\documentclass[12pt,a4paper,leqno]{amsart}

\usepackage{amsthm}
\usepackage{amsmath}
\usepackage{amssymb}
\usepackage[latin1]{inputenc}
\usepackage[T1]{fontenc}
\usepackage{tabularx}
\usepackage{graphicx}
\usepackage{enumerate}
\usepackage[all,2cell,ps]{xy}
\UseAllTwocells
\UsePSspecials{}

\theoremstyle{plain}
\newtheorem{thm}{Theorem}[subsection]
\newtheorem{prop}[thm]{Proposition}
\newtheorem{lem}[thm]{Lemma}
\newtheorem{sublem}[thm]{Sublemma}
\newtheorem{defn}[thm]{Definition}
\newtheorem{cor}[thm]{Corollary}

\theoremstyle{definition}
\newtheorem{nota}[thm]{Notation}
\newtheorem{ex}[thm]{Example}
\newtheorem{rem}[thm]{Remark}

\numberwithin{equation}{section}
\providecommand{\zeroset}[1]{\ensuremath{{\mathscr Z}\left(#1\right)}}
\renewcommand{\tan}[1]{\ensuremath{{\mathscr T}\left(#1\right)}}
\providecommand{\cat}[1]{\ensuremath{\mathcal{#1}}}
\providecommand{\kn}[1]{\ensuremath{\mathcal {#1}}}

\providecommand{\liouville}{\mathop{\omega}\nolimits}
\newcommand{\lio}[1][\X]{\liouville_{#1}}
\newcommand{\Sph}{\operatorname{S}\!}
\newcommand{\codim}{\operatorname{codim}}
\newcommand{\normcone}{{\operatorname{C}}}
\newcommand{\id}{\operatorname{id}}
\newcommand{\ktimes}{\operatorname{\!\times\!}}
\providecommand{\fp}[1]{\ensuremath{{{#1}_{\pi}}}}
\providecommand{\fd}[1]{\ensuremath{{{#1}_{d}}}}
\def\muHom{\mathop{\mu hom}\nolimits}
\def\nuHom{\mathop{\nu hom}\nolimits}
\def\Hom{\mathop{\rm Hom}\nolimits}
\def\Homo{\mathop{\mathscr Hom}\nolimits}

\def\IndHomo{\mathop{\mathscr {IH}om}\nolimits}

\def\Dlimind%
{\mathop{\textrm{\normalfont\textquotedblleft}%
\underrightarrow{\rm lim\raisebox{- 2pt}{}}
\textrm{\normalfont\textquotedblright}}\limits}
\providecommand{\f}[1]{\ensuremath{\mathscr{#1}}}
\def\df[#1]{%
  \ar@{}[#1]|(.35)*{}="A" \ar@{}[#1]|(.65)*{}="B" 
  \ar@{=>}"A";"B" }
\def\tdf[#1]{%
  \ar@{}[#1]|(.25)*{}="A" \ar@{}[#1]|(.55)*{}="B" 
  \ar@{=>}"A";"B" }

\newcommand{\on}{\operatorname}

\newcommand{\ctimes}[1]{\underset{#1}{\times}}

\newcommand{\ol}{\overline}

\newcommand{\geqs}{\geqslant}
\newcommand{\leqs}{\leqslant}
\newcommand{\mc}{\mathcal}

\newcommand{\me}{\matheus}
\newcommand{\ra}{\rightarrow}
\newcommand{\la}{\leftarrow}
\newcommand{\lra}[1][\,]{\xrightarrow{{\,\,#1\,\,}}}

\newcommand{\Z}{\mathbb{Z}}
\newcommand{\A}{\mathbb{A}}
\newcommand{\R}{\mathbb{R}}
\newcommand{\C}{\on{\mathbb{C}}}

\newcommand{\field}{\on{\mathbb{K}}}

\newcommand{\Db}{\operatorname{D}^{\operatorname{b}}}
\renewcommand{\SS}{\operatorname{SS}}
\newcommand{\supp}{\operatorname{Supp}}
\newcommand{\Supp}{\operatorname{Supp}}

\newcommand{\Dr}{\operatorname{R}\!}
\newcommand{\I}{\operatorname{I}}

\DeclareMathAlphabet{\matheus}{U}{eus}{m}{n}
\newdir{ (}{{}*!/-5pt/@^{(}}
\newdir{ )}{{}*!/-5pt/@_{(}}

\newcommand{\eq}{\begin{eqnarray}}
\newcommand{\eneq}{\end{eqnarray}}
\newcommand{\eqn}{\begin{eqnarray*}}
\newcommand{\eneqn}{\end{eqnarray*}}
\newcommand{\subsubset}{\subset\kern-.2ex\subset}
\newcommand{\cl}{\colon}
\newcommand{\DI}[1]{\Db(\on{I}(\field_{#1}))}
\newcommand{\op}{{\operatorname{op}}}

\newcommand{\db}[1]{*+!<0pc,.8ex>{#1}}
\newcommand{\ddb}[1]{*+!<0pc,1ex>{#1}}

\newcommand{\wCC}{\on{{\widetilde{\C}}}}

\newcommand{\bi}{\begin{enumerate}[{\rm(i)}]}
\newcommand{\ei}{\end{enumerate}}
\newcommand{\ba}{\begin{array}}
\newcommand{\ea}{\end{array}}
\newcommand{\seteq}{\on{{:=}}}
\newcommand{\bal}{\begin{align}}
\newcommand{\eal}{\end{align}}
\newcommand{\baln}{\begin{align*}}
\newcommand{\enaln}{\end{align*}}
\newcommand{\bdot}{\dot}

\newcommand{\isoto}{\lra[\,\sim\,]}
\newcommand{\hs}{\hspace}
\newcommand{\X}{{\mathfrak{X}}}
\newcommand{\Y}{{\mathfrak{Y}}}
\renewcommand{\d}{d}
\newcommand{\wC}{\on{\widetilde{\field}}}
\newcommand{\Ker}{\on{Ker}}
\newcommand{\K}{\on{K}}
\newcommand{\Mod}{\on{Mod}}

\newcommand{\comp}[1]{\mathop\circ\limits_{#1}}
\newcommand{\ot}{\leftarrow}
\newcommand{\bl}{\bigl}
\newcommand{\br}{\bigr}
\newcommand{\Bl}{\Bigl}
\newcommand{\Br}{\Bigr}

\DeclareMathAlphabet{\matheus}{U}{eus}{m}{n}
\newcommand{\mathscr}{\matheus}

\begin{document}
\title[Microlocalization of Ind-sheaves]%
{Microlocalization of Ind-sheaves}
\author{M. Kashiwara}
\address{Research Institute for Mathematical Sciences,
Kyoto University, Kyoto 606--8502, Japan
}
\email{masaki@kurims.kyoto-u.ac.jp}
\thanks{The author M.~K is partially supported by
Grant-in-Aid for Scientific Research (B1)13440006,
Japan Society for the Promotion of Science and
the 21-st century COE program
``Formation of an International Center of Excellence 
in the Frontier of Mathematics 
and Fostering of Researchers in Future Generations''.}

\author{P. Schapira}
\address{Universit{\'e} Pierre et Marie Curie,
Institut de Math{\'e}matiques, 175 rue du Chevaleret, 75013 Paris, 
France}
\email{schapira@math.jussieu.fr}

\author{F. Ivorra}
\address{Universit{\'e} Pierre et Marie Curie,
Institut de Math{\'e}matiques, 175 rue du Chevaleret, 75013 Paris, 
France}
\email{fivorra@math.jussieu.fr}

\author{I. Waschkies}
\address{
Universit{\'e} de Nice - Sophia Antipolis
Laboratoire J.A. Dieudonne
Parc Valrose
06108 Nice, France}
\email{ingo@math.nice.fr}
\thanks{The author I.~W is partially supported by the
same 21-st century COE program .}

\date{}

\maketitle

\begin{abstract}
Let $X$ be a $\on{C}^{\infty}$-manifold
and $T^*X$ its cotangent bundle. We construct a microlocalization 
functor
$ \mu_X\cl\DI{X} \lra \DI{T^*X}$,
where $\DI{X}$ denotes the bounded derived category of ind-sheaves of 
vector spaces on $X$ over a field $\field$. This functor satisfies
$ \Dr \me{H}om(\mu_X(\me{F}),\mu_X(\me{G}))\simeq \mu hom(\me{F},\me{G})$
for any $\me{F},\me{G}\in\Db(\field_X)$, thus generalizing the 
classical theory of microlocalization. Then we discuss the functoriality of 
$\mu_X$. The main result is the existence of a microlocal convolution morphism
$$ \mu_{X\times Y}(\me{K}_1)\overset{a}{\circ}\mu_{Y\times Z}(\me{K}_2) \lra 
  \mu_{X\times Z}(\me{K}_1\circ\me{K}_2)$$
which is an isomorphism under suitable non-characteristic 
conditions on $\me{K}_1$
and $\me{K}_2$.
\end{abstract}

\tableofcontents

\section*{Introduction}

This paper is based on ideas of the authors M.K.\ and P.S.\ announced in 
\cite{KS5} and developed in a preliminary manuscript of M.K.

The idea of microlocalization goes back to M.~Sato \cite{S} in 1969
who invented the functor of microlocalization of sheaves (along a
smooth submanifold of a real manifold) in order to analyze  the
singularities of hyperfunction solutions of systems of differential
equations in the cotangent bundle. 
This microlocalization procedure 
then  allowed Sato, Kashiwara and Kawai \cite{SKK} to define 
functorially the sheaf of rings of microdifferential operators
on the cotangent bundle $T^*X$ of a complex manifold $X$, a sheaf 
whose direct image is the sheaf of 
differential operators on $X$.

Then in the 80's, M.K.\ and P.S.\  (cf.\ \cite{KS2}, \cite{KS3}) 
developed a microlocal theory of sheaves 
on a $\on{C}^\infty$- manifold $X$, based on the notion of
microsupport (a conic involutive closed subset of the cotangent bundle
to $X$) and 
introduced in particular the functor $\mu hom$.
This is roughly speaking a functor which associate to a pair of
sheaves on $X$
the sheaf of microlocal morphisms between them.

On the other hand, the Riemann-Hilbert problem, solved by M.K.,
tells us that there is a one-to-one correspondence
between the regular holonomic modules over the ring of differential operators
and the perverse sheaves.
The notion of regular holonomic modules 
over the ring of differential operators can be easily microlocalized
to the notion of regular holonomic modules
over the ring of microdifferential operators and it is a natural
question to ask if there is a natural notion of microlocalization of
perverse sheaves, or, more generally a functor $\mu$  of microlocalization for
sheaves, the microsupport of a sheaf being the support of its
microlocalization and the functor $\mu hom$ being the internal $hom$ 
applied to the microlocalization. 
This is indeed what we do in this paper. 

As an application of the new functor $\mu$,  
 the author I.W.\ \cite{W} has recently constructed the stack
of microlocal perverse sheaves on a homogeneous symplectic manifold, 
after M.K.\ \cite{K} 
had constructed the stack of microdifferential modules.

The paper consists of two parts.
The first part is the technical heart of the paper. 
We define kernels on a $\on{C}^{\infty}$-manifold $X$, 
attached to the data of 
a closed submanifold $Z$ and a $1$-form $\sigma$ vanishing on $Z$.
Then we study its functorial properties. 
These kernels can be seen as 
``general'' microlocalization kernels, though their only role in this paper is 
to provide us with the tools for the proofs of the functorial properties of
$\mu$. 

In the second part we introduce the functor $\mu$,
which is the integral transform with respect to
the kernel $\on{K}_{T^*X}$ on $T^*X\times T^*X$ associated with
the fundamental $1$-form. 
 We discuss the functorial properties of $\mu$ 
deduced from the corresponding properties of the kernels studied in the 
first part. 
We then show how some classical microlocal properties can be generalized
to ind-sheaves. 
We give a comparison theorem 
between the micro-support of ind-sheaves $\me{F}$
and the support of its microlocalization $\mu(\me{F})$.

As an application,
we prove that, on a complex manifold $X$, 
$\mu hom$ induces a well-defined functor
$$\mu hom(\,\bullet\,,\me{O}_X)
\cl \Db(\C_X)^\op\lra\Db(\me{E}_X),$$
where $\me{E}_X$ is the ring of microdifferential operators.

The authors would like to thank A. D'Agnolo for many helpful comments.

\section{Microlocal kernels}
In all this paper, $\field$ denotes a field. 

        \subsection{Review on Ind-sheaves on manifolds}
In this section we shall give a short overview on the theory of 
ind-sheaves of \cite{KS1}. 

Let $X$ be a locally compact topological space 
with finite cohomological dimension,
$\on{Mod}(\field_X)$ the category of sheaves of $\field$-vector spaces on $X$,
and $\on{Mod}^{\on{c}}(\field_X)$ 
its full subcategory of sheaves with compact supports.

  We denote by $\I(\field_X)$ the category of ind-sheaves,
which is by definition the category of ind-objects of
$\on{Mod}^{\on{c}}(\field_X)$.
Then, $\I(\field_X)$ is an abelian category,
and its bounded derived category is denoted by 
  $\Db(\I(\field_X))$.  

There is a fully faithful exact functor
$$ \iota_X\colon\on{Mod}(\field_X)\lra \I(\field_X) \quad \mbox{given by} 
\quad {F}\mapsto 
      \Dlimind_{U\subsubset  X}{{F}_U}, $$
where the direct limit on the right is taken over the family of 
relatively compact open
subsets $U$ of $X$. In the sequel, 
we will regard $\on{Mod}(\field_X)$ as a full subcategory of
$\I(\field_X)$.

The functor $\iota_X$ admits an exact left adjoint functor
$$ \alpha_X\cl\I(\field_X) \lra \on{Mod}(\field_X), \quad 
  \Dlimind_{i\in I}{{F}_i}\mapsto \varinjlim_{i\in I}{{F}_i}. $$
Since $\iota_X$ is fully faithful, we have $\alpha_X\circ\iota_X\simeq 
\on{Id}_{\on{Mod}(\field_X)}$.\\
The functor $\alpha_X$ admits an exact fully faithful left adjoint 
$$ \beta_X\cl\on{Mod}(\field_X)\lra \I(\field_X). $$
Since $\beta_X$ is fully faithful, we get $\alpha_X\circ\beta_X\simeq
\on{Id}_{\on{Mod}(\field_X)}$. 
The functor $\beta_X$ is less easy to define than $\alpha_X$ and 
$\iota_X$. However, for a locally closed subset $S\subset X$,
$$\wC_S\seteq\beta_X(\field_S)$$
is described as follows.
Let $Z$ be a closed subset, then we have
$$ \wC_Z\simeq 
   \Dlimind_{Z\subset W}\field_{\ol{W}}, $$
where $W$ runs through the open subsets
containing $Z$.
If $U\subset X$ is an open subset 
then 
$$\wC_U\simeq \Dlimind_{V\subsubset U}\field_V, $$
where $V$ runs through the family of relatively compact open subsets of $U$. 
If $S\subset X$ 
is locally closed, then we can write $S=Z\cap U$ where $U$ is open 
and $Z$ is closed, and
$$\wC_S\simeq\wC_U\otimes\wC_W\simeq
 \Dlimind_{V\subsubset U,\, Z\subset W}
      \field_{V\cap \ol{W}}. $$
Therefore $\field_{V\cap\ol{W}}\ra \field_S$ induces a 
morphism $\wC_S\lra \field_S$ which is not an isomorphism in general.

Note that if $Z$ is closed and $S\subset Z$ is a locally closed subset, 
then
$$ \field_{S}\otimes\wC_Z\simeq \field_{S}. $$

The machinery of Grothendieck's six operations
is also applied to this context. We have the functors:
\eqn
f^{-1},\,f^!&\cl& \DI{Y}\to\DI{X},\\[1ex]
\Dr f_*,\,\Dr f_{!!}&\cl&\DI{X}\to\DI{Y},\\[1ex]
\Dr\me{IH}om&\cl&\DI{X}^\op\times\DI{X}\to\on{D}^+(\I(\field_X)),\\[1ex]
\otimes\ \ &\cl&\DI{X}\times\DI{X}\to\DI{X},
\eneqn
(here,  $f\cl X\to Y$ is a continuous map)
and we have the stack-theoretical hom
$$\Dr\me{H}om\cl\DI{X}^\op\times\DI{X}\to\on{D}^+(\field_X).$$
Note that the functor $\Dr\me{IH}om$ sends $\Db(\field_X)^\op\times\DI{X}$ to
$\DI{X}$ and $\Dr\me{H}om$ sends $\Db(\field_X)^\op\times\DI{X}$ to
$\Db(\field_X)$.

The inverse image functor $f^{-1}$ is a left adjoint of the direct image 
functor $\Dr f_*$.
The functor of direct image with proper support
$\Dr f_{!!}$ has a right adjoint functor $f^!$.
Most formulas of sheaves have their 
  counterpart in the theory of ind-sheaves, but some formulas are new. We
 shall not repeat them here and refer to \cite{KS1}. As an example we state 
  the following propositions:

\begin{prop}
  Consider a cartesian square
      $$ \xymatrix{
              X' \ar[r]^{f'}\ar[d]_{g'} & Y' \ar[d]^g \\
              X \ar[r]_f &  Y. }$$
  Then we have canonical isomorphisms
  $$ \Dr f'_{!!}{g'}^{-1}\simeq g^{-1}\Dr f_{!!}, \qquad
     \Dr f_*'g^{\prime!}\simeq g^!\Dr f_* ,\qquad
\Dr f'_{!!}g^{\prime!}\simeq g^!\Dr f_{!!}. $$
 \end{prop}
 Note that the last isomorphism has no counterpart in sheaf theory. 

\begin{prop}\label{prop:tenshom}
For a morphism $f\cl X\to Y$ and for $K\in\Db(\field_Y)$,
$\f F\in \DI{X}$, we have
\eqn
\ba{ll}
\Dr f_{!!}\Dr \me{IH}om(f^{-1}K,\f F)\simeq\Dr \me{IH}om(K,\Dr f_{!!}\f F)
&\text{in $\DI{Y}$,}\\[1ex]
\Dr f_{!}\Dr \me{H}om(f^{-1}K,\f F)\simeq\Dr \me{H}om(K,\Dr f_{!!}\f F)
&\text{in $\Db(\field_Y)$.}
\ea
\eneqn
\end{prop}

\begin{rem}
Let $Z$ be a closed subset of $X$ and
let $i\cl Z\to X$, $j\cl X\setminus Z\to Z$ be
the inclusion morphisms.
Then for $\f F$, $\f F'\in\DI{X}$, we have
\eq
&&\ba{ll}
\Dr j_{!!}j^{-1}\f F\simeq\wC_{X\setminus Z}\otimes \f F,&
\Dr i_*i^{-1}\f F\simeq\field_Z\otimes\f F,\\[.5ex]
\Dr j_{*}j^{-1}\f F\simeq\Dr\me{IH}om(\wC_{X\setminus Z},\f F),&
\Dr i_*i^{!}\f F\simeq\Dr\me{IH}om(\field_Z,\f F),\\[.5ex]
\Dr j_{*}j^{-1}\Dr\me{H}om(\f F',\f F)\simeq\Dr\me{H}om(\wC_{X\setminus Z}\otimes\f F',\f F).&
\ea
\eneq
Hence there are {\em not}
distinguished triangles
$$\Dr j_{!!}j^{-1}\f F\to\f F\to\Dr i_*i^{-1}\f F\lra[+1]\ \text{nor}
\ \Dr i_*i^{!}\f F\to\f F\to\Dr j_{*}j^{-1}\f F\lra[+1],$$
and instead there are distinguished triangles
\eq
&&\ba{l}
\Dr j_{!!}j^{-1}\f F\to\f F\to\f F\otimes\wC_Z\lra[+1]\ \text{and}
\ \Dr \me{IH}om(\wC_Z,\f F)\to\f F\to\Dr j_{*}j^{-1}\f F\lra[+1].
\ea
\eneq
\end{rem}

\medskip
The functor $\beta$ satisfies the following properties.

\begin{subequations}
\renewcommand{\theequation}{\theparentequation\kern.5ex\alph{equation}}
\begin{equation}
\text{$\beta_X({F})\otimes\beta_X({G})\simeq\beta_X({F}\otimes{G})$
for ${F}$, ${G}\in\Db(\field_X)$.}\label{eq:betashriek}
\end{equation}

For $f\cl X\to Y$ and ${G}\in\Db(\field_Y)$ and $\f G\in\DI{X}$, we have
\begin{equation}
f^{-1}\beta_Y({G})\simeq\beta_X(f^{-1}{G})\quad\text{and}\quad
f^!(\f G\otimes\beta_Y(G))\simeq f^!\f G\otimes\beta_X(f^{-1}G).
\end{equation}

For $\f F\in \DI{X}$ and $K$, $K'\in\Db(\field_X)$, we have
\begin{equation}
\ba{rcll}
\Dr\me{IH}om(K,\f F)\otimes\beta_X(K')&\simeq&
\Dr\me{IH}om\left(K,\f F\otimes\beta_X(K')\right)
\quad&\text{in $\DI{X}$,}\\[1ex]
\Dr\Homo( K,\f F)\otimes K'&\simeq&
\Dr\Homo\left( K,\f F\otimes\beta_X( K')\right)\quad&\text{in $\Db(\field_X)$.}
\ea\label{eq:betahom}\end{equation}

\end{subequations}


In general $\beta$ does not commute with direct image.
\begin{lem}\label{RemTrun}
  Consider a closed embedding $i\cl Z\hookrightarrow X$ 
and ${F}\in\Db(\field_Z)$.
  Then we have an isomorphism
   $$ \beta_X(\Dr i_*{F})\otimes \field_Z \simeq \Dr i_*\beta_Z({F}). $$ 
\end{lem}
\begin{proof}
We have
  $$ \beta_X(\Dr i_*{F})\otimes \field_Z\simeq
     \Dr i_*i^{-1}\beta_X(\Dr i_*{F})\simeq
     \Dr i_*\beta_Z(i^{-1}\Dr i_*{F})\simeq\Dr i_*\beta_Z({F}). $$
\end{proof}

The following fact will be used frequently in the paper:

\eq
&&\parbox{30em}{A morphism $u\cl\f F\to \f G$ in $\DI{X}$ is an isomorphism
if and only if
$\f F\otimes\wC_x\to \f G\otimes\wC_x$ is an isomorphism for
all $x\in X$.}
\eneq

\newcommand{\cross}{\times}
We list the commutativity of 
various functors. Here, ``{\Large{$\circ$}}''
means that the functors commute, and ``$\cross$'' that they do not.

\newcommand{\Yes}{{\Large$\circ$}}
\newcommand{\No}{$\cross$}
\begin{center}
\renewcommand{\arraystretch}{1.5}
\begin{tabular}{|c|c|c|c|c|c}
\cline{1-5}
& $\iota$ &\ $\alpha$\ &\ $\beta$\ & $\varinjlim$ & \\
\cline{1-5}
$\otimes$ &\Yes & \Yes & \Yes &  \Yes  & \\
\cline{1-5}
$f^{-1} $ & \Yes & \Yes & \Yes & \Yes & \\
\cline{1-5}
$\Dr f_* $ & \Yes & \Yes & \No  & \No  & \\
\cline{1-5}
$\Dr f_{!!}$ & \No  & \Yes & \No  & \Yes &\\
\cline{1-5}
$f^!$ & \Yes & \No  & \No  & \Yes &\multicolumn{1}{c}{} \\
\cline{1-5}
$\varinjlim$ & \No  & \Yes & \Yes & \multicolumn{2}{c}{} \\
\cline{1-4}
\end{tabular}
\end{center}

In the table, $\varinjlim$ means filtrant inductive limits.
For example, the commutativity of $\Dr f_{!!}$ and
$\varinjlim$ should be understood as in
Proposition \ref{functorJandcomposition} (i) below.
\medskip
\begin{nota}
For a continuous map $f\cl X\to Y$, we denote
by $\omega_{X/Y}$ the topological dualizing sheaf
$f^!\field_Y$, and $\omega_X=\omega_{X/\{\on{pt}\}}$. If $X$ and $Y$ are 
manifolds,
$\omega_{X/Y}\simeq\omega_X\otimes f^{-1}\omega_Y^{\otimes-1}$.
\end{nota}

\medskip
For three manifolds $X_i$ ($i=1,2,3$)
and for kernels $K\in\DI{X_1\times X_2}$ and $K'\in\DI{X_2\times X_3}$,
we define their convolution by
\eq\label{eq:int}
&&K\comp{X_2} K'=\Dr p_{13}{}_{!!}(p^{-1}_{12}K\otimes p^{-1}_{23}K'),
\eneq
where $p_{ij}$ is the projection from 
$X_1\times X_2\times X_3$ to $X_i\times X_j$.
We sometimes denote it simply by $K\circ K'$
when there is no risk of confusion.
This product of kernels satisfies the associative law:
$$(K\circ K')\circ K''\simeq K\circ (K'\circ K'')$$
for $K\in\DI{X_1\times X_2}$, $K'\in\DI{X_2\times X_3}$ and
$K''\in\DI{X_3\times X_4}$.
By taking $\{{\mathrm{pt}}\}$ as $X_3$
in \eqref{eq:int}, we obtain the integral transform functor:
$$K\circ\ \cl\DI{X_2}\to\DI{X_1}.$$
The following lemma is frequently used in \S \ref{sec:mic}.
\begin{lem}\label{lem:conv}
Let $f_k\cl X_k\to Y_k$ {\rm($k=1,2,3$)} be morphisms 
and $\f K_{ij}\in\DI{X_i\times X_j}$ and $\f L_{ij}\in\DI{Y_i\times Y_j}$.
\bi
\item
$\bl((f_1\times \id_{Y_2})^{-1}\f L_{12}\br)
\comp{Y_2}\bl((\id_{Y_2}\times f_3)^{-1}\f L_{23}\br)
\simeq (f_1\times f_3)^{-1}(\f L_{12}\comp{Y_2}\f L_{23})$ 
in $\DI{{X_1\times X_3}}$,\label{item}
\item
$\bl((f_1\times \id_{X_2})_{!!}\f K_{12}\br)
\comp{Y_2}\bl((\id_{X_2}\times f_3)_{!!}\f K_{23}\br)
\simeq (f_1\times f_3)_{!!}(\f K_{12}\comp{X_2}\f K_{23})$ 
in $\DI{{Y_1\times Y_3}}$,
\item
$\bl((\id_{Y_1}\times f_2)^{-1}\f L_{12}\br)
\comp{X_2}\f K_{23}
\simeq \f L_{12}\comp{Y_2}\Dr (f_2\times\id_{X_3})_{!!}\f K_{23}$ 
in $\DI{{Y_1\times X_3}}$.
\ei
\end{lem}

        \subsection{Kernels attached to $1$-forms}

Let us denote by $\pi_X\colon T^*X\to X$ the cotangent bundle to $X$.
For a closed submanifold $Z$ of $X$, we denote by $T^*_ZX$ its conormal bundle.
In particular, $T^*_XX$ is the zero section of
$T^*X$.
To a differentiable map  $f\cl X\ra Y$, we associate the diagram
$$ \xymatrix{
    {T^*X} &{T^*Y\smash{\ctimes{Y}}X} \ar[r]_{f_\pi} \ar[l]^{f_d} & {T^*Y.} 
    }$$

\begin{nota}
For a vector bundle $p\cl E\to X$, we denote by
$\bdot {E}$ the space $E$ with the zero section removed,
and by $\bdot p$ the projection $\bdot E\to X$.
For example, we use the notations
$\bdot\pi_X\cl\bdot{  T}^*X\to X$, $\bdot {T}^*_ZX$, etc.
\end{nota}
\begin{defn}
A kernel data is a triple $(X,Z,\sigma)$, 
where $X$ is a manifold, 
$Z$ is a closed submanifold of $X$ and 
$\sigma$ is a section of $T^*X\ctimes XZ\to Z$.
\end{defn}
We set $\tan{\sigma}=\sigma^{-1}(T^*_ZX)$ and 
$\zeroset{\sigma}=\sigma^{-1}(T^*_XX)$.
We have therefore
$$\zeroset{\sigma}\subset\tan{\sigma}\subset Z.$$
Each kernel data $(X,Z,\sigma)$ defines a closed cone $P_{\sigma}$ in 
$T_ZX\ctimes{X}\tan{\sigma}$ by
$$P_{\sigma}=\left\{(x,v)\in T_ZX; \text{$x\in\tan{\sigma}$ and
$\langle v,\sigma(x)\rangle\geqs 0$}\right\}.$$
Consider the deformation of the normal bundle to $Z$ in $X$ which will be denoted by 
$\widetilde X_Z$ or simply by $\widetilde X$ (see e.g. \cite{KS2}).
 We have the following commutative diagram 
where the squares marked by $\square$ are cartesian:
\begin{equation}\label{eq:deform}
\ba{c}
\xymatrix{&{\{0\}}\ar@<-1pt>@{ (->}[r]\ar@{}[rd]|{\square}&\R
\ar@{}[rd]|{\square}&{\{t\in\R;t>0\}}\ar@<1pt>@{ )->}[l]\\
{P_{\sigma}}\ar@<-1pt>@{ (->}[r] & {T_ZX}\ar[u]
        \ar@<-1pt>@{ (->}[r]^{s}\ar[d]^{\tau_Z} &       
          {\widetilde X_Z}\ar[d]^{p} \ar[u]_t& {\Omega}
           \ar@<1pt>@{ )->}[l]_{j}\ar[u]
           \ar[dl]^{\widetilde p}\\
        {} & {Z}\ar@<-1pt>@{ (->}[r]^{i} & {X.} & {}}
\ea
\end{equation}
Here $\Omega$ is the open subset defined by $\Omega=\{t>0\}$ 
for the natural smooth map $t\cl\widetilde X_Z\ra\R$. 
The normal bundle $T_ZX$ is identified with the inverse image of $0\in\R$
by $t$.
With a local coordinate system
$(x,z)=(x_1,\ldots,x_n,z_1,\ldots,z_m)$ of $X$ such that $Z$ is given by $x=0$,
$\widetilde X_Z$ has the coordinates
$(t,\tilde x,z)=(t,\tilde x_1,\ldots,\tilde x_n,z_1,\ldots,z_m)$
and $p$ is given by $p(t,\tilde x,z)=(t\tilde x,z)$.

Recall that the normal cone $C_Z(A)$ of a subset $A$ of $X$
is a closed cone of $T_ZX$ defined by
\eq
&&C_Z(A)=T_ZX\cap\overline{p^{-1}(A)\cap \Omega}.
\eneq

Note that $p$ is not smooth but the relative 
dualizing complex $\omega_{\widetilde{X}/X}$ is isomorphic to $\field_{\widetilde{X}}[1]$. 
In the sequel we will usually regard $P_{\sigma}$ as a closed subset of 
$\widetilde{X}_Z$ by $P_{\sigma}\subset T_ZX\subset \widetilde{X}_Z$.

\begin{defn}
(i) Let $(X,Z,\sigma)$ be a kernel data. We define the kernel 
  $\kn L_{\sigma}(Z,X)\in\Db(\I(\field_X))$ by
$$\kn L_{\sigma}(Z,X)=\Dr p_{!!}(\field_{\ol{\Omega}}\otimes
          \wC_{P_{\sigma}})\otimes\beta_X(\Dr i_*
      \omega_{Z/X}^{\otimes -1}).$$

\noindent
(ii) A morphism of kernel data
$f\cl (X_1,Z_1,\sigma_1)\to(X_2,Z_2,\sigma_2)$ 
is a morphism of manifolds $f\cl X_1\ra X_2$ 
satisfying
\bi
        \item{$f(Z_1)\subset Z_2$,}
        \item{$\sigma_1=f^*\sigma_2$.}
\ei
\end{defn}
\begin{rem}
Note that $\kn L_{\sigma}(Z,X)$ is supported on $\tan{\sigma}$, \emph{i.e.}
  $$ \kn L_{\sigma}(Z,X)\isoto\kn L_{\sigma}(Z,X)\otimes
      \wC_{\tan{\sigma}}.$$
\end{rem}

\medskip
This kernel behaves differently on
$\zeroset{\sigma}$ and outside.
We have
$$\kn L_{\sigma}(Z,X)\otimes \wC_{\zeroset{\sigma}}\simeq
\field_Z\otimes \wC_{\zeroset{\sigma}}$$
and $\kn L_{\sigma}(Z,X)\vert_{X\setminus \zeroset{\sigma}}$
is concentrated in degree $-\codim Z$ (see Corollary \ref{AffConExp}).

\medskip
In order to prove these facts,
we shall start by the following vanishing lemma.
\begin{lem}\label{lem:van}
\bi
\item
$\Dr p_{!!}(\field_{\Omega}\otimes\wC_{T_ZX})\simeq0$ and
$\Dr p_{!!}(\field_{\overline{\Omega}}\otimes\wC_{T_ZX})\simeq\Dr i_*\omega_{Z/X}$.
\item
Regarding $Z$ as the zero section of $T_ZX\subset \widetilde{X}_Z$,
we have
$$\Dr p_{!!}\left(\field_{\overline\Omega}\otimes\wC_{Z}\right)\simeq\wC_Z.
$$
\item
$\bl(\Dr {p}_{!!}(\field_{T_ZX}\otimes\wC_{P_{\sigma}})\br)\otimes\wC_{Z\setminus\zeroset{\sigma}}
\simeq 0$.
\ei
\end{lem}
\begin{proof}
(i)\quad
Since the problem is local, we may assume 
that $X$ is affine endowed with a system of 
global coordinates $(x,z)$ such that $Z=\{x=0\}$, $\widetilde{X}_Z=(t,\widetilde{x},z)$
and $p(t,\widetilde{x},z)=(t\widetilde{x},z)$.
We have then for all integer $j$
$$\Dr^jp_{!!}\big(\field_{\Omega}\otimes\wC_{T_ZX}\big)\kern-0.1em\simeq
  \kern-0.1em\Dr^jp_{!!}\Big(\Dlimind_{R>0,\,\varepsilon>0} 
\field_{\{ 0<t\leqs\varepsilon,\,|\widetilde x|<R\}}\Big)
\kern-0.1em\simeq\kern-0.1em
\Dlimind_{R>0,\,\varepsilon>0}\Dr^jp_!\field_{\{ 0<t\leqs\varepsilon,\,|\widetilde x|<R\}}
\simeq0,$$
which implies the first statement.
The last one follows from the distinguished triangle
$$
\Dr p_{!!}(\field_{\Omega}\otimes\wC_{T_ZX})\lra
\Dr p_{!!}(\field_{\overline{\Omega}}\otimes\wC_{T_ZX})
\lra\Dr p_{!!}(\field_{T_ZX})\lra[+1]$$
and $\Dr p_{!!}(\field_{T_ZX})\simeq
\Dr i_*\omega_{Z/X}$.

\medskip
\noindent
(ii)\quad
We have a chain of morphisms
$$\Dr p_{!!}\left(\field_{\overline\Omega}\otimes\wC_Z\right)
  \ra\Dr p_{!!}\left(\field_{\overline\Omega} \otimes\field_Z\right)
   \simeq\Dr p_{!!}\field_Z\simeq\field_Z.$$
which allows us 
to prove the isomorphism locally on $X$. 
With the coordinate system as above, we get for all integer $j$
\begin{eqnarray*}
 \Dr^jp_{!!}\left(\field_{\overline\Omega}\otimes\wC_{Z}\right) 
& \simeq & \Dr^jp_{!!}       
\Bigl(\Dlimind_{\varepsilon>0} \field_{\lbrace 0\leqs
  t\leqs\varepsilon, 
|\widetilde x|\leqs\varepsilon\rbrace}\Bigr)
\simeq\Dlimind_{\varepsilon>0}\Dr^jp_!
\field_{\lbrace 0\leqs t\leqs\varepsilon,|\widetilde x|
\leqs\varepsilon\rbrace}\\
        & \simeq &      
\begin{cases}
\Dlimind_{\varepsilon>0}\field_{\lbrace |x|\leqs\varepsilon^2\rbrace}
\simeq\wC_Z &                   \textrm{if $j=0$,}\\
0 & \textrm{if $j\neq 0$.}

                        \end{cases}
\end{eqnarray*}

\noindent
(iii)\quad
For $z_0\in\tan{\sigma}\setminus\zeroset{\sigma}$, we have
$$\left(\Dr {p}_{!!}(\field_{T_ZX}\otimes\wC_{P_{\sigma}})\right)
\otimes\wC_{z_0}
\simeq\Dr {p}_{!!} \left(\field_{T_ZX}
\otimes\wC_{P_{\sigma}\cap p^{-1}(z_0)}\right).$$
Set $\sigma(z_0)=\langle \xi_0,dx\rangle\not=0$.
Then we have
$$\field_{T_ZX}\otimes\wC_{P_{\sigma}\cap p^{-1}(z_0)}\simeq
\Dlimind_{R>0,\varepsilon>0}\field_{\{t=0,\,\,-\varepsilon
\le\langle \xi_0,\tilde x\rangle,\,|\tilde x|<R\}},$$
and for all integer $j$ 
$$\bl(\Dr^j{p}_{!!}(\field_{T_ZX}\otimes\wC_{P_{\sigma}})\br)\otimes\wC_{z_0}
\simeq\wC_{z_0}\otimes\Dlimind_{R>0,\,\varepsilon>0}\Dr^j{p}_!
\bl(\field_{\{t=0,\,\,-\varepsilon
\le\langle \xi_0,\tilde x\rangle,\,|\tilde x|<R\}}\br)\simeq0.$$
\end{proof}

\begin{lem}\label{lem:migi}
  There is a natural morphism
 $$\kn L_{\sigma}(Z,X)\lra
        \wC_{\tan{\sigma}}\otimes
\beta_X\left(\Dr i_*\omega^{\otimes -1}_{Z/X}\right).$$
\end{lem}
\begin{proof}
Regard $\tan{\sigma}$
as a subset of $\widetilde{X}_Z$ by $\tan{\sigma}\subset Z\subset T_ZX\subset 
\widetilde{X}_Z$. Then we get a natural morphism
\begin{eqnarray*}
&&\kn L_{\sigma}(Z,X)\ra
\Dr p_{!!}\left(\field_{\overline\Omega}\otimes\wC_{\tan{\sigma}}\right)
\otimes\beta_X\left(\Dr i_*\omega^{\otimes -1}_{Z/X}\right).
\end{eqnarray*}
Hence the desired morphism is obtained by Lemma \ref{lem:van} (ii).
\end{proof}

The following lemma provides a useful distinguished triangle to study some 
properties of the kernel $\kn L_{\sigma}(Z,X)$.

\begin{lem}\label{DT}
        There is a natural distinguished triangle
$$\Dr p_{!!}\bl(\field_{\Omega}\otimes\wC_{P_{\sigma}}\br)
\otimes\beta_{X}\bl(\Dr i_*\omega_{Z/X}^{\otimes-1}\br)\lra
     \kn L_{\sigma}(Z,X)\lra\Dr p_{!!}\bl(\field_{T_ZX}\otimes\wC_{P_{\sigma}}\br)
     \otimes\Dr i_*\omega_{Z/X}^{\otimes-1}\lra[+1].$$
\end{lem}

\begin{proof}
It is enough to apply the triangulated functor 
$\Dr p_{!!}\bl(\,\raisebox{-.2ex}{\Large{$\cdot$}}\,\otimes\wC_{P_{\sigma}}\br)
\otimes\beta_{X}(\Dr i_*\omega_{Z/X}^{\otimes-1})$ 
to the distinguished triangle 
\eq\label{eq:Omega}
&&\field_{\Omega}\lra\field_{\overline\Omega}
\lra\field_{T^*_ZX}\lra[+1],
\eneq
and to use $\field_Z\otimes \beta_{X}(\Dr i_*\omega_{Z/X}^{\otimes-1})
\simeq
\Dr i_*\omega_{Z/X}^{\otimes-1}$.
\end{proof}

Recall that $\zeroset{\sigma}$ is the set  of zeroes
of $\sigma$, i.e. $\zeroset{\sigma}=\sigma^{-1}(T^*_XX)\subset Z$.
\begin{prop}\label{DG1}
We have
$$\kn L_{\sigma}(Z,X)\otimes\wC_{\zeroset{\sigma}}
\simeq \field_Z\otimes\wC_{\zeroset{\sigma}}.$$
  In particular, if $\sigma=0$, then
  $ \kn L_{\sigma}(Z,X) \simeq \field_Z.$ 
\end{prop}
\begin{proof}
By the definition of $\zeroset{\sigma}$, 
the cone $P_{\sigma}\ctimes{Z}\zeroset{\sigma}$ coincides with
$T_ZX\ctimes{Z}\zeroset{\sigma}$. 
Hence we have
$\field_{\overline{\Omega}}\otimes \wC_{P\sigma}
\otimes p^{-1}\wC_{\zeroset{\sigma}}\simeq \field_{\overline{\Omega}}
\otimes p^{-1}\wC_{\zeroset{\sigma}}$,
which implies
\baln
\kn L_{\sigma}(Z,X)\otimes\wC_{\zeroset{\sigma}}
&\simeq\Dr p_{!!}(\field_{\overline{\Omega}}\otimes\wC_{T_ZX})
\otimes \widetilde\field_{\zeroset{\sigma}}\otimes
\beta_X(\Dr i_*\omega_{Z/X}^{\otimes-1}).
\end{align*}
Hence the result follows from
Lemma \ref{lem:van} (i).
\end{proof}

\begin{prop}
Let $(X,Z,\sigma)$ be a kernel data, and set
$X_0=X\setminus\zeroset{\sigma}$ and
$Z_0=Z\setminus\zeroset{\sigma}$.
Then there is a natural distinguished triangle
$$\Dr j_{!!}\kn L_{\sigma_0}\left(Z_0,X_0\right)\lra\kn L_{\sigma}(Z,X)
   \lra\field_Z\otimes \wC_{\zeroset{\sigma}}\lra[+1],$$
where $\sigma_0$ is the restriction of $\sigma$ to $Z_0$
and $j$ denotes the open immersion $X_0\hookrightarrow X$.
\end{prop}
\begin{proof}
We have the distinguished triangle
$$\kn L_{\sigma}(Z,X)\otimes\wC_{X_0}
\lra\kn L_{\sigma}(Z,X)\lra
\kn L_{\sigma}(Z,X)\otimes\wC_{\zeroset{\sigma}}
\lra[+1].$$
The first term is isomorphic to 
$\Dr j_{!!}\kn L_{\sigma_0}\left(Z_0,X_0\right)$,
and the last term is isomorphic to
$\field_Z\otimes \wC_{\zeroset{\sigma}}$
by Lemma \ref{DG1}.
\end{proof}

\begin{cor}\label{TrivMor}
        There are natural morphisms 
        $$\field_Z\lra\kn L_{\sigma}(Z,X)\lra
        \wC_{\tan{\sigma}}\otimes
        \beta_X\left(\Dr i_*\omega^{\otimes -1}_{Z/X}\right).$$
\end{cor}

\begin{proof}
The first arrow is constructed as
an immediate consequence of the preceding proposition and the obvious inclusion $P_{\sigma}\subset P_0=T_ZX$.
The last arrow follows from
Lemma \ref{lem:migi}.
\end{proof}

\begin{prop}\label{DG3}
   Assume the section ${\sigma}$ never vanishes. Then
  $$\kn L_{\sigma}(Z,X)
\simeq\Dr p_{!!}\left(\field_{\Omega}\otimes\wC_{P_{\sigma}}\right)
\otimes\beta_X
         \big(\Dr i_*\omega_{Z/X}^{\otimes{-1}}\big)
\simeq
\Dlimind_U\field_U\otimes\beta_X
         \big(\Dr i_*\omega_{Z/X}^{\otimes{-1}}\big)
            \otimes\wC_{\tan{\sigma}},$$
  where the inductive limit is taken over the family of open subsets 
$U$ of $X$ such that
  $$ P_\sigma\cap\normcone_Z(U)\subset Z.$$
Here, $Z$ is regarded as the zero section of $T_ZX$.
\end{prop}
Remark that the set of such $U$'s is a filtrant ordered set 
by the inclusion order.
\begin{proof} 

By Lemma \ref{DT} and Lemma \ref{lem:van} (iii), we have
$$\kn L_{\sigma}(Z,X)\simeq\Dr p_{!!}\left(\field_{\Omega}
\otimes\wC_{P_{\sigma}}\right)\otimes
\beta_{X}\left(\Dr i_*\omega_{Z/X}^{\otimes-1}\right).$$
Hence it is enough to show
$$\Dr p_{!!}\left(\field_{\Omega}\otimes\wC_{P_{\sigma}}\right)\simeq
\Dlimind_U\field_U\otimes\wC_{\tan{\sigma}}.$$
Since we have $Z\cap U=\varnothing$
on a neighborhood of $\tan{\sigma}$,
$p^{-1}(U)\cap\Omega=p^{-1}(U)\cap\ol\Omega$
is a closed subset of
$\Omega$ and we get the 
following chain of natural morphisms~:
$$p^{-1}\field_U\simeq\field_{p^{-1}(U)}\lra\field_{p^{-1}(U)\cap\Omega}
\lra\field_{\Omega}\lra\field_{\Omega}\otimes \wC_{P_{\sigma}}.$$
Since $\overline{p^{-1}(U)\cap\Omega}\cap P_{\sigma}
=\normcone_Z(U)\cap P_{\sigma}$ 
is contained in the zero section of $T_ZX$,
$\Supp(p^{-1}\field_U\otimes\wC_{P_{\sigma}})$ 
is proper over $Z$.
Hence we have a chain of morphisms
$$\field_U\lra p_*(p^{-1}\field_U\otimes\wC_{P_{\sigma}})
\simeq p_{!!}(p^{-1}\field_U\otimes\wC_{P_{\sigma}})
\lra p_{!!}\left(\field_{\Omega}\otimes\wC_{P_{\sigma}}\right),$$
which provides a natural morphism
$$\Dlimind_U\field_U\lra
\Dr p_{!!}\left(\field_{\Omega}\otimes\wC_{P_{\sigma}}\right).$$
By tensorisation we get the morphism
\eq\label{eq:UW}
&&
\Dlimind_U\field_U\otimes\wC_{\tan{\sigma}}\lra
\Dr p_{!!}\left(\field_{\Omega}\otimes\wC_{P_{\sigma}}\right).\eneq
We shall now show that this morphism is an isomorphism. 
It is enough to show that
\eqref{eq:UW} is an isomorphism after tensoring by $\wC_{x_0}$
for any $x_0\in \tan{\sigma}$.
Let us take  local coordinate system  
$(x,z)$ of $X$ such that $Z=\left\{ x=0\right\}$.
We may assume $x_0=(0,0)$, and we set 
${\sigma}(x_0)=\langle\xi_0,\d x\rangle$. 
We then have 
$$\Dr p_{!!}\left(\field_{\Omega}\otimes\wC_{P_{\sigma}}\right)
\otimes\wC_{x_0}\simeq
\Dr p_{!!}\left(\field_{\Omega}\otimes\wC_{P_{\sigma}}
\otimes\wC_{p^{-1}(x_0)}\right)
\simeq\Dr p_{!!}\left(\field_{\Omega}\otimes
\wC_{P_{\sigma}\cap p^{-1}(x_0)}\right),$$
and 
$$\field_{\Omega}\otimes\wC_{P_{\sigma}\cap p^{-1}(x_0)}
\simeq\Dlimind_{\substack{V\subsubset\widetilde X_Z,\,\,
P_{\sigma}\cap p^{-1}(x_0)\subset V'}}
\field_{\Omega\cap V\cap\overline{V'}}
\simeq \wC_{x_0}\otimes
\Dlimind_{R>0,\,\varepsilon_1>0,\,\varepsilon_2>0}
\field_{A_{R,\varepsilon_1,\varepsilon_2}},$$
where we have set 
$$A_{R,\varepsilon_1,\varepsilon_2}
=\left\{(t,\widetilde x,z)\in\widetilde X_Z\ ;\ 
0<t\leqs\varepsilon_1,\,-\varepsilon_2\leqs\langle\xi_0,
\,\widetilde x\rangle,|\widetilde x|<R\right\}.$$ 
Hence for all integer $j$, we have
$$ \Dr^jp_{!!}\left(\field_{\Omega}\otimes\wC_{P_{\sigma}}\right)
\otimes\wC_{x_0}\simeq \wC_{x_0}\otimes
\Dlimind_{R>0,\,\varepsilon_1>0,\,\varepsilon_2>0}
\Dr^jp_! \field_{A_{R,\varepsilon_1,\varepsilon_2}}.$$
We have 
\eqn
p^{-1}((x,z))&\simeq&
\{t\in\R;0<t\le\varepsilon_1,\,-\varepsilon_2\le \langle\xi_0,t^{-1}x\rangle,
\,|t^{-1}x|<R\}\\
&\simeq&\{t\in\R;R^{-1}|x|<t\le\varepsilon_1,
\,-\varepsilon_2^{-1}\langle\xi_0,x\rangle\le t\},
\eneqn
and hence
$$\Dr p_! (\field_{A_{R,\varepsilon_1,\varepsilon_2}})\simeq
\field_{\bigl\{R^{-1}|x|<-\varepsilon_2^{-1}\langle x,\,\xi_0\rangle
\leqs\varepsilon_1\bigr\}}.$$
Taking the limit we can use a cofinality argument to get
$$\Dr p_{!!}\left(\field_{\Omega}\otimes
\wC_{P_{\sigma}}\right)\otimes\wC_{x_0}\simeq
\wC_{x_0}\otimes\Dlimind_{\varepsilon>0}
\field_{\left\{ (x,z)\in X\ ;
\ -\langle\xi_0,x\rangle>\varepsilon |x|\right\}}.$$
Then the theorem follows from the following easy sublemma.
\end{proof}
\begin{sublem}
\bi
\item
Let $U=\{(x,z)\in X; \varepsilon|x|<-\langle \xi_0,x\rangle\}$.
Then 
  $P_\sigma\cap\normcone_Z(U)\subset Z$. 
   \item
Let $U\subset X$ be an open subset such that 
$P_\sigma\cap\normcone_Z(U)\subset Z$.
Then there exist $\varepsilon>0$ and $\delta>0$ such that 
    $$ U\cap \{|(x,z)|\leqs \delta\} \subset
          \{(x,z)\in X ; -\langle x,\xi_0\rangle >\varepsilon|x|\}.$$ 
\ei
\end{sublem}

\begin{cor}\label{AffConExp}
Let $(X,Z,\sigma)$ be a kernel data. 
Assume that $X$ is endowed with a 
local coordinate system $(x,z)$ such that $Z=\{x=0\}$ 
and $\sigma$ is a nowhere vanishing section. Then,
writing $\sigma(z)
=\langle\sigma_1(z),\d x\rangle+\langle\sigma_2(z),\d z\rangle$,
we have
$$\kn L_{\sigma}\left(Z,X\right)\simeq
\wC_{\{x=0,\,\sigma_2(z)=0\}}\otimes
\Dlimind_{\varepsilon>0}
\field_{\bigl\{(x,z) ; -\langle\sigma_1(z),x\rangle >\varepsilon |x|\bigr\}}\,
[\codim Z].$$

\end{cor}

\begin{rem}
\bi
\item
We have $$\alpha_X\bl(\kn{L}_\sigma(Z,X)\br)\simeq \field_{\zeroset{\sigma}}.$$
\item
Let $(X,Z,\sigma_1)$ and $(X,Z,\sigma_2)$ be kernel data,
and let $W$ be a closed subset of $Z$ such that 
$\sigma_1(x)=\sigma_2(x)$ for all $x\in W$.
Since $P_{\sigma_1}\cap\tau_Z^{-1}W=P_{\sigma_2}\cap\tau_Z^{-1}W$, we have
$$\kn L_{\sigma_1}(X,Z)\otimes\wC_W
\simeq \kn L_{\sigma_2}(X,Z)\otimes\wC_W.$$
\ei
\end{rem}

        \subsection{Functorial Properties}

In this subsection, we will investigate the 
behavior of microlocal kernels $\f L_\sigma(Z,X)$ under inverse and proper 
direct images, and under convolution.

Let $f\cl (X_1,Z_1,\sigma_1)\to(X_2,Z_2,\sigma_2)$ be morphism of 
kernel data. We have the  diagrams of manifolds
$$  
\ba{c}
\xymatrix{
     T^*_{Z_1}X_1&  \db{T^*_{Z_2}X_2\ctimes{Z_2}Z_1}\ar[l]_(.55){f_d}\ar[r]
^(.6){f_\pi}
       &T^*_{Z_2}X_2 \\
\tan{\sigma_1}\ar[u]^{\sigma_1} 
&\db{\tan{\sigma_2}\ctimes{Z_2}Z_1}\ar[u]\ar[r]\ar@<1pt>@{ )->}[l]&  
\tan{\sigma_2}\,\ar[u]_(.45){\sigma_2}
}\ea
\qquad\text{and}\quad\ba{c}
\xymatrix{
{\widetilde{X}_1}\ar[r]^{\tilde{f}}\ar[d]_{p_1}&{\widetilde{X}_2}\ar[r]^t
\ar[d]^{p_2}&\R\\
X_1\ar[r]^f&X_2
}\ea
$$
where $\widetilde{X}_k=\widetilde{X_k}_{Z_k}$ ($k=1,2$).
We denote by $i_k\cl Z_k\hookrightarrow X_k$ the inclusion map.
We have 
\eq\label{eq:tan12}
P_{\sigma_1}\ctimes{X_2}\tan{\sigma_2}=\tilde f^{-1}(P_{\sigma_2}).
\eneq
\begin{prop}\label{pro:cleanmorpkernel}
Let $f\cl(X_1,Z_1,\sigma_1)\ra (X_2,Z_2,\sigma_2)$ be a morphism of
kernel data. Assume 
that $Z_1=f^{-1}(Z_2)$
and the morphism $f\cl X_1\ra X_2$ is clean with respect to $Z_2$
{\rm(}i.e.\ $(T_{Z_1}X_1)_x\to (T_{Z_2}X_2)_{f(x)}$ is injective 
for any $x\in Z_1${\rm)}. 
Then there exists a natural 
morphism
  $$f^{-1}\kn L_{{\sigma}_2}(Z_2,X_2)\lra\kn L_{{\sigma}_1}(Z_1,X_1) 
    \otimes \beta_{X_1}(\Dr     i_{1*}\omega_{Z_1/Z_2})\otimes 
     \omega_{X_1/X_2}^{\otimes -1}\otimes\wC_{f^{-1}\tan{\sigma_2}}.$$ 
\end{prop}
\begin{proof}
Since $f$ is clean, $\widetilde X_1\to \widetilde X_2\ctimes{X_2}X_1$ is
a closed embedding
and there is a morphism of functors
$f^{-1}\Dr p_2{}_{!!}\to \Dr p_1{}_{!!}\tilde f^{-1}$ which induces
a natural morphism
\begin{eqnarray}
        f^{-1}\kn L_{{\sigma}_2}(Z_2,X_2) & \simeq & 
f^{-1}\Dr p_{2!!}\left(\field_{\ol{\Omega}_2}\otimes
\wC_{P_{{\sigma}_2}}\right)
\otimes f^{-1}\beta_{X_2}(\Dr i_{2*}\omega_{Z_2/X_2}^{\otimes-1})\nonumber\\
& \ra & \Dr p_{1!!}\widetilde f^{-1}
\left(\field_{\ol{\Omega}_2}\otimes\wC_{P_{{\sigma}_2}}\right) 
\otimes f^{-1}\beta_{X_2}(\Dr i_{2*}\omega_{Z_2/X_2}^{\otimes -1})
\label{eq:clean}\\
        & \simeq & \Dr p_{1!!}\left(\field_{\ol{\Omega}_1}\otimes\wC_{\widetilde f^{-1}(P_{{\sigma}_2})}\right)     \otimes\beta_{X_1}(f^{-1}\Dr i_{2*}\omega_{Z_2/X_2}^{\otimes-1}).\nonumber
\end{eqnarray} 
By \eqref{eq:tan12}, we have a morphism
\eq\label{eq:cl2}
&&f^{-1}\kn L_{{\sigma}_2}(Z_2,X_2)\lra
\Dr p_{1!!}\left(\field_{\ol{\Omega}_1}\otimes\wC_{P_{\sigma_1}} \right)
\otimes f^{-1}(\wC_{\tan{\sigma_2}})
\otimes\beta_{X_1}(f^{-1}\Dr i_{2*}\omega_{Z_2/X_2}^{\otimes-1}).
\eneq
Hence, to get the desired morphism, it is enough to remark that 
$$
f^{-1}\Dr i_{2*}\omega_{Z_2/X_2}^{\otimes -1} 
\simeq
\Dr i_{1*}\left(\omega_{Z_1/X_1}^{\otimes -1}
\otimes \omega_{Z_1/Z_2}\otimes i_1^{-1}
\omega_{X_1/X_2}^{\otimes -1}\right)
 \simeq  \Dr i_{1*}\omega_{Z_1/X_1}^{\otimes -1}
\otimes\Dr i_{1*}\omega_{Z_1/Z_2}\otimes\omega_{X_1/X_2}^{\otimes -1}.
$$
\end{proof}
By adjunction, we obtain:
\begin{cor}
 Under the hypothesis of the Proposition \ref{pro:cleanmorpkernel}, 
we have a natural morphism
 $$\kn L_{{\sigma}_2}(Z_2,X_2)\lra\Dr f_*\left(\kn L_{{\sigma}_1}
   (Z_1,X_1) \otimes\beta_{X_1}(\Dr i_{1*}\omega_{Z_1/Z_2})\otimes 
    \omega_{X_1/X_2}^{\otimes -1}\otimes f^{-1}\wC_{\tan{\sigma_2}}\right).$$ 
\end{cor}
\begin{prop}\label{InvIm}
Let $f\cl(X_1,Z_1,\sigma_1)\ra (X_2,Z_2,\sigma_2)$ be a morphism of
kernel data.
Assume that $f^{-1}(Z_2)=Z_1$ and $f$ is transversal to $Z_2$. 
Then we have a natural isomorphism
$$f^{-1}\kn L_{{\sigma}_2}(Z_2,X_2)\isoto
    \kn L_{{\sigma}_1}(Z_1,X_1).$$ 
\end{prop}
\begin{proof}
Indeed if $f$ is transversal,
$\widetilde{X}_1\to \widetilde{X}_2\ctimes{X_2}X_1$ is an isomorphism
and
$Z_1\cap f^{-1}(\tan{\sigma_2})=\tan{\sigma_1}$,
which implies that
the morphism \eqref{eq:clean} 
as well as \eqref{eq:cl2} is an 
isomorphism. We have furthermore 
$\omega_{Z_1/Z_2}\simeq i_1^{-1}\omega_{X_1/X_2}$.
\end{proof}

\begin{prop}
Let $f\cl(X_1,Z_1,\sigma_1)\ra (X_2,Z_2,\sigma_2)$ be a morphism of
kernel data.
Then there is a natural morphism
$$\Dr f_{!!}\Big(\kn L_{{\sigma}_1}(Z_1,X_1)\otimes 
   \beta_{X_1}(\Dr i_{1*}\omega_{Z_1/Z_2})\Big)
   \lra\kn L_{{\sigma}_2}(Z_2,X_2).$$
\end{prop}
\begin{proof}
The left hand side is isomorphic to
\begin{eqnarray}\label{eq:123}
&&\ba{rcl}
\lefteqn{\kern-2em\Dr f_{!!}\Bl(\Dr p_{1!!}\bl(\field_{\ol{\Omega}_1}\otimes
\wC_{P_{{\sigma}_1}}\br)\otimes\beta_{X_1}(\Dr i_{1*}
\omega_{Z_1/X_1}^{\otimes -1})\otimes\beta_{X_1}(\Dr i_{1*}
\omega_{Z_1/Z_2})\Br) }\\
& \kern 3em\simeq  & \Dr f_{!!}\Bl(\Dr p_{1!!}\bl(\field_{\ol{\Omega}_1}
\otimes\wC_{P_{{\sigma}_1}}\br)\otimes\omega_{X_1/X_2}\otimes
\beta_{X_1}(f^{-1}\Dr i_{2*}\omega^{\otimes-1}_{Z_2/X_2})\Br)\\
        & \kern 3em\simeq  & \Dr f_{!!}\Dr p_{1!!}\Bl(\bl(\field_{\ol{\Omega}_1}
\otimes\wC_{P_{{\sigma}_1}}\bl)\otimes p_1^{-1}
\omega_{X_1/X_2}\Br)
\otimes\beta_{X_2}(\Dr i_{2*}\omega^{\otimes-1}_{Z_2/X_2})\\
        & \kern 3em\simeq  & \Dr p_{2!!} \Dr \widetilde f_{!!}
\Bl(\widetilde f^{-1}\field_{\ol{\Omega}_2}\otimes\widetilde {\field}_{P_{{\sigma}_1}}
\otimes p_1^{-1}\omega_{X_1/X_2}\Br)\otimes\beta_{X_2}(\Dr i_{2*}
\omega^{\otimes-1}_{Z_2/X_2})\\
        & \kern 3em\simeq  & \Dr p_{2!!} \Bl(\field_{\ol{\Omega}_2}\otimes
\Dr \widetilde  f_{!!}\bl(\wC_{P_{{\sigma}_1}}\otimes 
\omega_{\widetilde{X}_1/\widetilde{X}_2}\br)\Br)\otimes
\beta_{X_2}(\Dr i_{2*}\omega^{\otimes-1}_{Z_2/X_2}). 
\ea
\end{eqnarray} 
Hence, it is enough to construct a morphism
\eq\label{eq:124}
&&\Dr\widetilde f_{!!}\left(\widetilde {\field}_{P_{{\sigma}_1}}
\otimes\omega_{\widetilde X_1/\widetilde X_2}\right)
\lra \wC_{P_{{\sigma}_2}}.
\eneq
By adjunction it is enough to construct a morphism
$\widetilde {\field}_{P_{{\sigma}_1}}
\otimes\omega_{\widetilde X_1/\widetilde X_2}
\lra \widetilde f^{!}\wC_{P_{{\sigma}_2}}$.
However by \eqref{eq:tan12},
we have
$$\widetilde {\field}_{P_{{\sigma}_1}}
\otimes\omega_{\widetilde X_1/\widetilde X_2}
\lra\widetilde {\field}_{P_{{\sigma}_1}\ctimes{X_2}\tan{\sigma_2}}
\otimes\omega_{\widetilde X_1/\widetilde X_2}
\simeq\widetilde f^{-1}\wC_{P_{{\sigma}_2}}
\otimes\omega_{\widetilde X_1/\widetilde X_2}
\simeq\widetilde f^!\wC_{P_{{\sigma}_2}},$$
where the last isomorphism follows from \eqref{eq:betashriek}.
\end{proof}

\begin{cor}\label{DirIm}
Let $f\cl(X_1,Z_1,\sigma_1)\ra (X_2,Z_2,\sigma_2)$ be a morphism,
and assume that 
$f$ is smooth and induces an isomorphism from $Z_1$ to $Z_2$.
Then we have a natural isomorphism
$$\Dr f_{!!}\kn L_{{\sigma}_1}(Z_1,X_1)\isoto
\kn L_{{\sigma}_2}(Z_2,X_2).$$
\end{cor}
\begin{proof}
By the assumption, we have
$\tan{\sigma_2}\ctimes{Z_2}Z_1=\tan{\sigma_1}$.
By \eqref{eq:123}, it is enough to prove that
\eqref{eq:124} is an isomorphism.
Since 
$P_{\sigma_1}=\widetilde f^{-1}(P_{\sigma_2})$,
we have
\eqn
\Dr\widetilde f_{!!}\left(\wC_{P_{\sigma_1}}
\otimes\omega_{\widetilde X_1/\widetilde X_2}\right)
&\simeq&\wC_{P_{{\sigma}_2}}
\otimes\Dr\widetilde f_{!!}
\left(\wC_{T_{Z_1}X_1}\otimes\omega_{\widetilde X_1/\widetilde X_2}
\right).
\eneqn
Hence we have reduced the problem to
$$\Dr\widetilde f_{!!}
\left(\wC_{T_{Z_1}X_1}\otimes\omega_{\widetilde X_1/\widetilde X_2}\right)
\simeq \wC_{T_{Z_2}X_2}.$$
Since $f$ is smooth,
we can take local coordinate systems $(x,z)$ on $X_2$ 
and $(x,y,z)$ on $X_1$ such that
$Z_2=\{x=0\}$, $Z_1=\{x=0,\,y=0\}$
and $f$ is given by the projection.
We then take a coordinate system
$(t,\tilde x,z)$ on $\widetilde {X_2}$ and 
$(t,\tilde x,\tilde y, z)$ on $\widetilde {X_1}$.
The associated morphism
$\widetilde f\cl\widetilde {X_1}\to\widetilde {X_2}$
is given by
$(t,\tilde x,\tilde y, z)\to(t,\tilde x,z)$.
Then we can check easily
$\Dr\widetilde f_{!!}(\wC_{T_{Z_1}X_1}
\otimes\omega_{\widetilde X_1/\widetilde X_2})\simeq
\Dr\widetilde f_{!!}(\wC_{\{t=0\}}
\otimes\omega_{\widetilde X_1/\widetilde X_2})\simeq 
\wC_{\{t=0\}}$.
\end{proof}

\begin{lem}\label{lem:vanDir}
Let $(X,Z,\sigma)$ be a kernel data on $X$, and
let $f\cl X\to Y$ be a smooth morphism which induces
a closed embedding
$Z\hookrightarrow Y$.
Assume that
$\sigma(x)\notin T^*_{f(x)}Y$ for any $x\in\tan{\sigma}$.
Then we have
$$\Dr f_{!!}\kn L_\sigma(Z,X)\simeq0.$$
\end{lem}

\begin{proof}
For any $x_0\in \tan{\sigma}$, take a local coordinate system
$(y,z)=(y_1,\ldots,y_n,z_1,\ldots,z_m)$ of $Y$ in a neighborhood of $f(x_0)$
such that $f(Z)$ is given by $y=0$.
Then we can take a local coordinate system $(t,x,y,z)$ of 
$X$ in a neighborhood of $x_0$ such that $Z$ is given by $\{t=0,\,x=0,\,y=0\}$,
and $\sigma(x_0)=-\d t(x_0)$.
Then we have
$$\kn L_\sigma(Z,X)\otimes\wC_{x_0}\simeq
\bigl(\Dlimind_{\delta>0,\,\varepsilon>0}\field_{F_{\delta,\varepsilon}}\bigr)
\otimes\beta_X(\Dr i_*\omega^{\otimes-1}_{Z/X})\otimes\wC_{x_0},$$
where
$$F_{\delta,\varepsilon}=\{(t,x,y,z);\delta\ge t> \varepsilon(|x|+|y|)\}.$$
Hence, $\bl(\Dr f_{!!}(\kn L_{\sigma}(Z,X)\br)\otimes\wC_{f(x_0)}\simeq
\Dr f_{!!}(\kn L_{\sigma}(Z,X)\otimes\wC_{x_0})\simeq0$ follows from
$$\on{R}^jf_!(\field_{F_{\delta,\varepsilon}})
\simeq0\quad\text{for any $j\in\Z$.}
$$
\end{proof}

\begin{prop}\label{CloRes}
Let $f\cl(X_1,Z_1,\sigma_1)\ra (X_2,Z_2,\sigma_2)$ be a morphism of
kernel data,
and assume that $f$ is a closed immersion
which induces an isomorphism $Z_1\isoto Z_2$.
Then there is a natural isomorphism
$$\kn L_{\sigma_1}(Z_1,X_1)\isoto f^!\kn L_{\sigma_2}(Z_2,X_2).$$
\end{prop}
\begin{proof}
Since $f$ is a closed immersion, 
we get the commutative diagrams
$$
\ba{c}
\xymatrix{
{Z_1}\ar@<-1pt>@{ (->}[r]^{i_1}\ar[d]_{\sim}
     \ar@{}[rd]|{\square} & {X_1}\ar@<-1pt>@{ (->}[d]^f\ar@{}[rd]|{} & 
     {\widetilde {X_1}}\ar@<-1pt>@{ (->} [d]^{\widetilde f}
      \ar[l]_{p_1}\ar@{}[rd]|{\square} & {\Omega_1}\ar@<1pt>@{ )->}_{j_1}
      [l] \ar@<1pt>@{ (->}[d] \\
        {Z_2}\ar@<-1pt>@{ (->}[r]_{i_2} & {X_2} & {\widetilde X_2} \ar[l]^{p_2} & 
         {\Omega_2}\ar@<1pt>@{ )->}[l]^{j_2}
}
\ea\quad\text{and}\quad
\ba{c}
\xymatrix{T_{Z_1}X_1 \ar@{ (->}[r]^(.6){s_1}\ar[d]_{T_Zf}\ar@{}[rd]|{\square} 
& {\widetilde {X_1}} \ar[d]^{\widetilde f} \\
        T_{Z_2}X_2 \ar@{ (->}[r]_(.6){s_2} & {\widetilde X_2\,,}}
\ea
$$ 
in which the squares marked by $\square$ are cartesian.
Recall the adjunction isomorphism $f^!\Dr f_{!!}\simeq\id$. Hence it is enough
to construct an isomorphism
$$\Dr f_{!!}\kn L_{{\sigma_1}}(Z_1,X_2)\isoto\Dr f_{!!}f^!
   \kn L_{\sigma_2}(Z_2,X_2).$$
Next recall that
$$\Dr f_{!!}f^!\kn L_{\sigma_2}(Z_2,X_2)\simeq\Dr\IndHomo\left(\field_{X_1},
  \kn L_{\sigma_2}(Z_2,X_2)\right).$$
Therefore we may write:
\begin{eqnarray*}
        \Dr f_{!!}f^!\kn L_{\sigma_2}(Z_2,X_2) & \simeq & 
\Dr\IndHomo\bigg(\field_{X_1},\Dr p_2{}_{!!}\Big(\field_{\ol{\Omega_2}}\otimes 
        \beta_{\widetilde X_2}\big(\field_{P_{\sigma_2}}\otimes
p_2^{-1}\Dr i_2{}_*\omega_{Z_2/X_2}^{\otimes-1}\big)\Big)\bigg)\\
        & \simeq &\Dr p_2{}_{!!}\Dr\IndHomo\Big(p_2{}^{-1}\field_{X_1},
\field_{\ol{\Omega_2}}\otimes\beta_{\widetilde{X_2}}\big(\field_{P_{\sigma_2}}
\otimes p_2^{-1}\Dr i_2{}_*\omega_{Z_2/X_2}^{\otimes-1}\big)\Big)\\
        & \simeq & \Dr p_2{}_{!!}
\Big(\Dr\IndHomo\big(p_2^{-1}\field_{X_1},\field_{\ol{\Omega_2}}\big)\otimes
\beta_{\widetilde{X_2}}(\field_{P_{\sigma_2}}\otimes 
p_2^{-1}\Dr i_2{}_*\omega_{Z_2/X_2}^{\otimes-1})\Big).
\end{eqnarray*}
On the other hand, $P_{\sigma_1}=\widetilde f^{-1}P_{\sigma_2}$ implies
\begin{align*}
\Dr f_{!!}&\kn L_{\sigma_1}(Z_1,X_1)  \simeq  
\Dr f_{!!}\Dr p_1{}_{!!}\Big(\field_{\ol{\Omega_1}}\otimes\beta_{\widetilde{X_1}}
\big(\field_{P_{{\sigma_1}}}\otimes p_1^{-1}\Dr i_1{}_*\omega_{Z_1/X_1}^{\otimes-1}
\big)\Big)\\
& \simeq  \Dr p_2{}_{!!}\Dr\widetilde f_{!!}
\Big(\field_{\widetilde f^{-1}(\ol{\Omega_2})}\otimes\beta_{\widetilde{X_1}}
\big(\widetilde f^{-1}\field_{P_{\sigma_1}}\otimes
\widetilde f^{-1}p_2^{-1}
\Dr s_2{}_*\omega_{Z_2/X_2}^{\otimes-1}
\otimes p_1{}^{-1}\omega_{X_1/X_2}\big)\Big)\\
& \simeq  \Dr p_2{}_{!!}\Dr\widetilde f_{!!}
\bigg(\widetilde f^{-1}\Big(\field_{\ol{\Omega_2}}\otimes\beta_{\widetilde{X_2}}
\big(\field_{P_{\sigma_2}}\otimes p_2^{-1}\Dr i_2{}_*
\omega_{Z_2/X_2}^{\otimes-1}\big)\Big)
\otimes \omega_{\widetilde X_1/\widetilde X_2}\bigg)\\
& \simeq  
\Dr p_2{}_{!!}\Big(\field_{\ol{\Omega_2}}\otimes
\beta_{\widetilde {X_2}}\big(\field_{P_{\sigma_2}}
\otimes p_2^{-1}\Dr i_2{}_*
\omega_{Z_2/X_2}^{\otimes-1}\big)\otimes
\Dr\widetilde f_{!!}\omega_{\widetilde X_1/\widetilde X_2}\Big),
\end{align*}
and it is enough to show that
$$\Dr\IndHomo(p_2^{-1}\field_{ X_1},\field_{\ol{\Omega_2}})\simeq 
\field_{\ol{\Omega_2}}\otimes
\Dr\widetilde f_{!!}\omega_{\widetilde X_1/\widetilde X_2}.$$
However we have the natural chain of isomorphisms
\begin{eqnarray*}
  \Dr\IndHomo(p_2^{-1}\field_{ X_1},\field_{\ol{\Omega_2}})& \simeq &
   \Dr\IndHomo(p_2^{-1}\field_{ X_1},\Dr j_2{}_*\field_{\Omega_2})\\
   &\simeq&\Dr j_2{}_*
      \Dr\IndHomo(j_2{}^{-1}p_2^{-1}\field_{ X_1},\field_{\Omega_2})
         \simeq  \Dr j_2{}_*\Dr\IndHomo(\field_{\Omega_1},\field_{\Omega_2}).
\end{eqnarray*}
On the other hand, we have, as an object of $\DI{\Omega_2}$,
$$\Dr\IndHomo(\field_{\Omega_1},\field_{\Omega_2})
\simeq j_2^{-1}\,\Dr\widetilde{f}_*\omega_{\widetilde X_1/\widetilde X_2},$$
and hence
\begin{eqnarray*}
\Dr\IndHomo(p_2^{-1}\field_{ X_1},\field_{\ol{\Omega_2}})& \simeq &
\Dr j_2{}_*j_2^{-1}\Dr\widetilde{f}_*\omega_{\widetilde X_1/\widetilde X_2}\\
&\simeq&   \Dr j_2{}_*\field_{\Omega_1}\otimes
\Dr\widetilde{f}_*\omega_{\widetilde X_1/\widetilde X_2}
\simeq 
\field_{\ol{\Omega_1}}\otimes
\Dr\widetilde f_{!!}\omega_{\widetilde X_1/\widetilde X_2}.
\end{eqnarray*}
\end{proof}

\begin{prop}\label{pro:kerndatsamemnf}
Let $(X,Z_1,\sigma_1)$ and $(X,Z_2,\sigma_2)$ be kernel data on the same 
base manifold $X$. 
Assume that $Z_1$, $Z_2$ are transversal submanifolds. 
Then there is a natural morphism
 $$\kn L_{\sigma_1}(Z_1,X)\otimes\kn L_{\sigma_2}(Z_2,X)\lra
     \kn L_{\sigma_1+\sigma_2}\left(Z_1\cap Z_2,X\right)\otimes
\wC_{\tan{\sigma_1}\cap\tan{\sigma_2}}.$$
\end{prop}
\begin{proof}
Set $Z=Z_1\cap Z_2$,
${\sigma}={\sigma}_1+{\sigma}_2$
and $N=\tan{\sigma_1}\cap\tan{\sigma_2}\subset \tan{\sigma}\subset Z$.

(i) Assume first that ${\sigma}_1(x)$ and ${\sigma}_2(x)$ are 
linearly independent 
vectors of $T^*X$ for every $x\in Z$.
Then we have
$$\kn L_{{\sigma}_k}\left(Z_k,X\right)
\otimes\wC_N\simeq\Dlimind_{U_k}\field_{U_k}\otimes
\wC_N
\otimes\beta_X\left(\Dr i_{k*}\omega^{\otimes -1}_{Z_k/X}\right),$$
where the inductive limits is taken over the family of 
open subsets $U_k$ of $X$ 
such that $\normcone_{Z_k}\left(U_k\right)\cap P_{{\sigma}_k}
\subset Z_k$. For such open subsets 
$U_1$, $U_2$,
we have 
$$\normcone_Z\left(U_1\cap U_2\right)\cap \bl(P_{\sigma}\ctimes{Z}N\br)
\subset Z,$$
since $P_\sigma\ctimes{Z}N\subset P_{\sigma_1}\cup P_{\sigma_2}$.
Hence we get a natural morphism
\eqn
&&\kn L_{{\sigma}_1}\left(Z_1,X\right)\otimes
\kn L_{{\sigma}_2}\left(Z_2,X\right)\otimes\wC_N\\
&&\hs{10ex}\simeq \Bigl(\Dlimind_{U_1}\field_{U_1}
\otimes\beta\left(\Dr i_{1*}\omega^{\otimes -1}_{Z_1/X}\right)
\Bigr)
\otimes
\Bigl(\Dlimind_{U_2}\field_{U_2}
\otimes\beta\left(\Dr i_{2*}\omega^{\otimes-1}_{Z_2/X}\right)\Bigr)
\otimes\wC_N\\
&&\hs{10ex}\lra \Bigl(\Dlimind_U\field_U\Bigr)\otimes
\beta\left(\Dr i_{1*}\omega^{\otimes-1}_{Z_1/X}\right)\otimes 
\beta\left(\Dr i_{2*}\omega^{\otimes -1}_{Z_2/X}\right)
\otimes\wC_N,
\eneqn
where $U$ ranges over the family of open subsets of $X$ such that 
$\normcone_Z\left(U\right)\cap\bl(P_{\sigma}\ctimes{Z}N\br)\subset Z$.\\
Since $Z_1$ and $Z_2$ are transversal submanifolds of $X$, we have
$\omega^{\otimes -1}_{Z/X}\simeq
\bl(\omega^{\otimes -1}_{Z_1/X}\vert_Z\br)\otimes 
\bl(\omega^{\otimes -1}_{Z_2/X}\vert_Z\br)$.
Hence we obtain
$$\Dlimind_U\field_U\otimes
\beta_X\left(\Dr i_{1*}\omega^{\otimes -1}_{Z_1/X}\right)\otimes
\beta_X\left(\Dr i_{2*}\omega^{\otimes -1}_{Z_2/X}\right)\otimes\wC_N
\simeq\kn L_{\sigma}(Z,X)\otimes\wC_N,$$
which provides the desired morphism.

\medskip
(ii) Consider the general case.
We set
$\A^n_X=X\times \R^n$ for $n=1,2$.
We use coordinates $(x,t_1,t_2)$ on $\A^2_X$.
We regard the manifold $\A^1_{Z_k}$ as a submanifold 
of $\A^2_X$ by
$$\A^1_{Z_k}\seteq\left\{(x,t_1,t_2)\,;\,x\in Z_k,\,t_k=0\right\},$$
and $\A^1_X$ as the submanifold $\{t_2=0\}$ of $\A^2_X$.
We identify $Z$ with
$$\A^1_{Z_1}\cap\A^1_{Z_2}
=\left\{(x,t_1,t_2)\,;\,x\in Z,\,t_1=t_2=0\right\}.$$
Thus we obtain the following commutative diagrams
\eqn
&&\ba{c}
\xymatrix{&{X}\ar@<-1pt>@{ (->}[r]^{i} & {\A^1_X}\ar@<-1pt>@{ (->}[r]^{i'}
   \ar@{}[rd]|{\textrm{tr}}  & {\A^2_X}\\
Z\ar@<-1pt>@{ (->}[r]&{Z_1}\ar@<-1pt>@{ ->}^{\sim}[r]\ar@<-1pt>@{ (->}[u] & 
      {Z_1}\ar@<-1pt>@{ (->}[r]\ar@<-1pt>@{ (->}[u] & {\A^1_{Z_1}}
      \ar@<-1pt>@{ (->}[u]_{j_1}}
\ea\qquad\textrm{and}\qquad
\ba{c}
        \xymatrix{&{X}\ar@<-1pt>@{ (->}[r]^{i}\ar@{}[rd]|{\textrm{tr}} & 
       {\A^1_X}\ar@<-1pt>@{ (->}[r]^{i'}  & {\A^2_X}\\
Z\ar@<-1pt>@{ (->}[r]&{Z_2}\ar@<-1pt>@{ (->}[r]\ar@<-1pt>@{ (->}[u] & 
        {\A^1_{Z_2}}\ar@<-1pt>@{->}[r]^{\sim}\ar@<-1pt>@{ (->}[u] & 
       {\A^1_{Z_2}}\ar@<-1pt>@{ (->}[u]_{j_2}}
\ea
\eneqn
where $j_1(z_1,t)=(z_1,0,t)$ and $j_2(z_2,t)=(z_2,t,0)$. Note that the squares marked
with $\on{tr}$ are transversal.
Define the sections 
\eqn
&&\ba{rcccl}
\tilde\sigma_1=\sigma_1+dt_1&\cl& \A^1_{Z_1}&\lra
&T^*\A^2_X,\\[1ex]
\tilde\sigma_2=\sigma_2+dt_2&\cl& \A^1_{Z_2}&\lra
& T^*\A^2_X,\\[1ex]
\tilde\sigma=\sigma_1+\sigma_2+dt_1+dt_2&\cl& Z&\lra& T^*\A^2_X.
\ea
\eneqn
Clearly $\tilde\sigma_1$ and $\tilde\sigma_2$
are linearly independent at each point, and the result in the first part
gives a morphism
$$\kn L_{\tilde{\sigma}_1}\left(\A^1_{Z_1},\A^2_X\right)\otimes
  \kn L_{\tilde{\sigma}_2}\left(\A^1_{Z_2},\A^2_X\right)\lra
   \kn L_{\tilde\sigma}(Z,\A^2_X)\otimes\wC_N.$$
We then deduce morphisms with the help of
Proposition \ref{InvIm} and Proposition \ref{CloRes}
\begin{eqnarray*}
\kn L_{{\sigma}_1}\left(Z_1,X\right)\otimes
\kn L_{{\sigma}_2}\left(Z_2,X\right) 
& \simeq & i^!\kn L_{\tilde\sigma_1}\left(Z_1,\A^1_X\right)
\otimes i^{-1}\kn L_{\sigma_2}\left(\A^1_{Z_2},\A^1_X\right)\\
& \ra &  i^!\left(\kn L_{\tilde\sigma_1}\left(Z_1,\A^1_X\right)\otimes
\kn L_{\sigma_2}\left(\A^1_{Z_2},\A^1_X\right)\right)\\
& \simeq & i^!\left(
{i'}^{-1}\kn L_{\tilde\sigma_1}\left(\A^1_{Z_1},\A^2_X\right)
\otimes{i'}^!\kn L_{\tilde\sigma_2}\left(\A^1_{Z_2},\A^2_X\right)\right)\\
&\ra&i^!{i'}^!\left(\kn L_{\tilde\sigma_1}\left(\A^1_{Z_1},\A^2_X\right)
\otimes\kn L_{\tilde\sigma_2}\left(\A^1_{Z_2},\A^2_X\right)\right)\\
&\ra& i^!{i'}^!\bigl(\kn L_{\tilde\sigma}(Z,\A^2_X)\otimes\wC_N\bigr)
\simeq\kn L_{\sigma}(Z,X)\otimes\wC_N,
\end{eqnarray*}
which completes the proof.
\end{proof}
\begin{rem}
Although we do not give proofs,
the following two facts hold.
\bi
\item
If $\sigma_1$ and $\sigma_2$ are linearly independent, the two morphisms 
constructed in the parts (i) and (ii) of the proof of Proposition
\ref{pro:kerndatsamemnf}  coincide.
\item
If $(X,Z_3,\sigma_3)$ is a third kernel data such that 
$(Z_1, Z_2)$,  $(Z_1,Z_3)$ and $(Z_2, Z_3)$  are transversal in $X$ and
that $(Z_1\cap Z_3,Z_2\cap Z_3)$ is transversal in $Z_3$,
then
the following diagram is commutative
where $N=\tan{\sigma_1}\cap\tan{\sigma_2}\cap\tan{\sigma_3}$:
$$
\xymatrix@C=.45cm{
\kn L_{\sigma_1}(Z_1,X)\otimes\kn L_{\sigma_2}(Z_2,X)\otimes
      \kn L_{\sigma_3}(Z_3,X)\ar[r]\ar[d]&
\kn L_{\sigma_1+\sigma_2}(Z_1\cap Z_2,X)\otimes
   \kn L_{\sigma_3}(Z_3,X)\otimes\wC_N\ar[d]\\
\kn L_{\sigma_1}(Z_1,X)\otimes\kn L_{\sigma_2+\sigma_3}
 \left(Z_2\cap Z_3,X\right)\otimes\wC_N\ar[r]&
\kn L_{\sigma_1+\sigma_2+\sigma_3}\left(Z_1\cap Z_2\cap Z_3,X\right)
 \otimes\wC_N,}
$$
    \emph{i.e.} the composition morphisms are associative.
\ei
\end{rem}

\begin{lem}\label{lem:vandir}
Let $(X,Z_1,\sigma_1)$, $(X,Z_2,\sigma_2)$ be kernel data on $X$ and
assume that $Z_1$, $Z_2$ are transversal submanifolds of $X$ and 
that $\sigma_1$ and $\sigma_2$ never vanish.
Let $f\cl X\ra Y$ be a smooth morphism which
induces a closed embedding $Z_1\cap Z_2\hookrightarrow Y$.
Assume the following condition:
$$\text{$\bigl(\R_{\ge0}\sigma_1(x)+\R_{\ge0}\sigma_2(x)\bigr)\cap
T_{f(x)}^*Y=\{0\}$ for every 
$x\in\tan{\sigma_1}\cap\tan{\sigma_2}$.}$$
Here $T^*_{f(x)}Y$ is regarded as a subspace of $T^*_xX$ by $f_d$.
Then we have
$$\Dr f_{!!}\left(\kn L_{\sigma_1}(Z_1,X)\otimes
   \kn L_{\sigma_2}(Z_2,X)\right)\simeq0.$$
\end{lem}
\begin{proof}
Let us show that 
$$\Dr f_{!!}\left(\kn L_{\sigma_1}(Z_1,X)\otimes
   \kn L_{\sigma_2}(Z_2,X)\otimes\wC_{x_0}\right)\simeq0$$
for any $x_0\in\tan{\sigma_1}\cap\tan{\sigma_2}$.
We first reduce the proof to the case 
where $X$ is of relative dimension one over $Y$. 
Assume the assertion to be true for relative one-dimensional morphisms. 
Set $E=T_{x_0}(f^{-1}f(x_0))$.
Then by the assumption, $E$ satisfies
$\bigl(\R_{\ge0}\sigma_1(x_0)+\R_{\ge0}\sigma_2(x_0)\bigr)\cap
E^\perp=\{0\}$.
Hence there exists a line
$\ell\subset E$ such that
$\bigl(\R_{\ge0}\sigma_1(x_0)+\R_{\ge0}\sigma_2(x_0)\bigr)\cap
\ell^\perp=\{0\}$.
Decompose $f$ into the composition of smooth morphisms
$X\lra[g] Y'\lra[h] Y$ on a neighborhood of $x_0$
such that $g$ and $h$ are smooth
and $T_{x_0}(g^{-1}g(x_0))=\ell$.
Then $g$ satisfies the conditions in the lemma.
Hence applying to $g$ the relative one-dimensional morphism
case, we obtain
$\Dr g_{!!}\left(\kn L_{\sigma_1}(Z_1,X)\otimes
   \kn L_{\sigma_2}(Z_2,X)\otimes\wC_{x_0}\right)\simeq0$,
which implies the desired result.

\medskip
Now assume that 
$f$ has relative dimension one.
Since $\sigma_k(x_0)\notin T^*_{f(x_0)}Y$,
the map $Z_k\to Y$ is a (local) embedding, and
$T_{x_0}Z_k=f_*^{-1}\bl(T_{f(x_0)}Z'_k\br)\cap\sigma_k(x_0)^{-1}(0)$, where
$Z'_k\seteq f(Z_k)\subset Y$.
Then $Z'_1$ and $Z'_2$ are transversal, and
$f(Z_1\cap Z_2)$ is a hypersurface of $Z'_1\cap Z_2'$ since
\begin{align*}
\codim_Y(f(Z_1\cap Z_2))&=\codim_X(Z_1\cap Z_2)-1
=\codim_X(Z_1)+\codim_X(Z_2)-1\\
&=\codim_Y(Z'_1)+\codim_Y(Z'_2)+1=\codim_Y(Z'_1\cap Z'_2)+1.
\end{align*}
Since $T_{x_0}(Z_1\cap Z_2)=f_*^{-1}\bl(T_{f(x_0)}(Z'_1\cap Z'_2)\br)
\cap \sigma_1(x_0)^{-1}(0)\cap \sigma_2(x_0)^{-1}(0)$, the vectors
$\sigma_1(x_0)$ and $\sigma_2(x_0)$ are linearly independent.
By multiplying by a positive constant, we may therefore assume
that 
$$\sigma_1(x_0)-\sigma_2(x_0)\in T^*_{f(x_0)}Y\setminus\{0\}.$$
Take a local coordinate system $(t,y_1,y_2,z)$ of $Y$
such that
$$\text{$Z'_k=\{y_k=0\}$ 
and $\sigma_2(x_0)-\sigma_1(x_0)=\d t$.}$$
Then take a local coordinate system $(x,t,y_1,y_2,z)$ of $X$ such that
$\sigma_1(x_0)=-\d x$ (and hence $\sigma_2(x_0)=\d t-\d x$),
and $Z_1=\{y_1=0,\, x=0\}$ and $f$ is given by forgetting $x$. 
Set $Z_2=\{y_2=0,\,x=\varphi(t,y_1,z)\}$. Then replacing
$\varphi(t,y_1,z)$ with $t$, we may assume from the beginning that 
$$\text{$Z_2=\{y_2=0,\,x=t\}$, $Z_1\cap Z_2=\{y_1=0,\,y_2=0,\,x=t=0\}$.}$$
Then we have
$$
\kn L_{{\sigma}_1}(Z_1,X)\otimes
\kn L_{{\sigma}_2}(Z_2,X)\otimes\wC_{x_0}
\simeq
\Dlimind_{\delta>0,\,\varepsilon>0}\bigl(\field_{U^1_{\delta,\,\varepsilon}}
\otimes\field_{U^2_{\delta,\,\varepsilon}}\bigr)
\otimes
\beta_X\big(\Dr i_{1*}\omega^{\otimes-1}_{Z_1/X}
\otimes\Dr i_{2*}\omega^{\otimes-1}_{Z_2/X}\big)
\otimes\wC_{x_0},
$$
where the open sets $U^k_{\delta,\,\varepsilon}$ are given by
\eqn
&&U^1_{\delta,\,\varepsilon}=\left\{
\varepsilon |y_1|<x\le\delta\right\}\quad\text{and}\quad
U^2_{\delta,\,\varepsilon}=\left\{
\varepsilon |y_2|<x-t\le\delta\right\}.
\eneqn
Hence
we have
$$
U^1_{\delta,\,\varepsilon}\cap
U^2_{\delta,\,\varepsilon}=\left\{
\max(\varepsilon |y_1|,\,\varepsilon |y_2|+t)
<x\le\min(\delta,\,\delta+t)\right\}.
$$
Then the result follows from
$$\Dr f_!(\field_{U^1_{\delta,\,\varepsilon}\cap
U^2_{\delta,\,\varepsilon}})\simeq0.$$
\end{proof}

\begin{prop}\label{DirImTens}
Let $(X,Z_1,\sigma_1)$, $(X,Z_2,\sigma_2)$ be kernel data on $X$ and
$(Y,Z,\sigma)$ a kernel data on $Y$. Assume that $Z_1$, $Z_2$ are transversal
submanifolds of $X$.
Let $f\cl X\ra Y$ be a smooth morphism which
induces an isomorphism $Z_1\cap Z_2\isoto Z$.
Let $N$ be a closed subset of $\tan{\sigma_1}\cap\tan{\sigma_2}$
satisfying the following conditions:
\bi
\item
$\zeroset{\sigma_1}\cap\zeroset{\sigma_2}\subset N$,
\item
$f^*\sigma(x)=\sigma_1(x)+\sigma_2(x)$ for every 
$x\in N$,
\item
$\sigma_1(x)\not\in T^*_{f(x)}Y$ for 
any $x\in N\setminus\bigl(\zeroset{\sigma_1}\cup\zeroset{\sigma_2}\bigr)$,
\item
$\bigl(\R_{\ge0}\sigma_1(x)+\R_{\ge0}\sigma_2(x)\bigr)\cap
T_{f(x)}^*Y=\{0\}$ for every 
$x\in \bigl(\tan{\sigma_1}\cap\tan{\sigma_2}\bigr)\setminus N$,
\item
the morphism $Z_k\ra Y$ is smooth 
at each point of $\zeroset{\sigma_k}$ for $k=1,2$.
%
\ei
Then there is a natural isomorphism
$$\Dr f_{!!}\left(\kn L_{\sigma_1}(Z_1,X)\otimes
   \kn L_{\sigma_2}(Z_2,X)\right)\isoto\kn L_{\sigma}(Z,Y)
   \otimes\wC_{f(N)}.$$
\end{prop}

\begin{proof}
The morphism is obtained as the composition
\eqn
\Dr f_{!!}\left(\kn L_{\sigma_1}(Z_1,X)\otimes
   \kn L_{\sigma_2}(Z_2,X)\right)
&\lra&\Dr f_{!!}\bigl(\kn L_{\sigma_1+\sigma_2}(Z_1\cap Z_2,X)
\otimes\wC_N\bigr)\\
&\simeq&\Dr f_{!!}\bigl(\kn L_{f^*\sigma}(Z_1\cap Z_2,X)\otimes\wC_N\bigr)
\lra\kn L_{\sigma}(Z,Y)\otimes\wC_{f(N)}.\eneqn
In order to see that it is an isomorphism, 
it is enough to prove the isomorphism
$$\Dr f_{!!}\left(\kn L_{\sigma_1}(Z_1,X)\otimes
\kn L_{\sigma_2}(Z_2,X)\otimes\wC_{x_0} \right)
\isoto \kn L_{\sigma}(Z,Y)\otimes\wC_{f(N)}\otimes\wC_{f(x_0)}$$
for any $x_0\in \tan{\sigma_1}\cap\tan{\sigma_2}$.

\medskip
\noindent
(a)\quad
Assume first that 
$\sigma_1(x_0)=\sigma_2(x_0)=0$.  Then,
(i) implies $x_0\in N$, and we have $\sigma(f(x_0))=0$ by (ii).
Hence Proposition \ref{DG1} implies
\begin{eqnarray*}
\Dr f_{!!}\left(\kn L_{{\sigma}_1}(Z_1,X)\otimes
\kn L_{{\sigma}_2}(Z_2,X)\otimes \wC_{x_0} \right)
   & \simeq & 
\Dr f_{!!}\left(\field_{Z_1}\otimes\field_{Z_2}\otimes \wC_{x_0}
  \right)\\
        & \simeq & \field_Z\otimes\wC_{f(x_0)}\simeq\kn L_{\sigma}(Z,Y)
\otimes\wC_{f(N)}\otimes\wC_{f(x_0)}.
\end{eqnarray*}

\medskip
\noindent
(b)\quad
Assume ${\sigma}_1(x_0)=0$ and ${\sigma}_2(x_0)\not=0$. Then we have
\begin{eqnarray*}
\Dr f_{!!}\left(\kn L_{{\sigma}_1}(Z_1,X)\otimes
     \kn L_{{\sigma}_2}(Z_2,X)\otimes\wC_{x_0} \right) 
         &\simeq&  \Dr f_{!!}\left(\field_{Z_1}\otimes
    \kn L_{{\sigma}_2}(Z_2,X)\otimes \wC_{x_0}\right)\\
 &\simeq&  \Dr f_{!!}i_{1!!}i^{-1}_1\kn L_{{\sigma}_2}(Z_2,X)
\otimes\wC_{f(x_0)},
\end{eqnarray*}
where $i_1\cl Z_1\lra X$ is the inclusion.
Proposition \ref{InvIm} implies
$i^{-1}_1\kn L_{{\sigma}_2}(Z_2,X)\simeq
\kn L_{{\sigma}_2}(Z_1\cap Z_2,Z_1)$.
Note that $Z_1\ra Y$  is smooth at $x_0$ by the assumption (v).
If $x_0\in N$, then
Corollary \ref{DirIm}, along with by the hypothesis (ii), implies
$\Dr f_{!!}i_{1!!}\kn L_{{\sigma}_2}(Z_1\cap Z_2,Z_1)
\simeq\kn L_{\sigma}(Z,Y)$.
Assume $x\in\bigl(\tan{\sigma_1}\cap\tan{\sigma_2}\bigr)\setminus N$.
Then (iv) implies that $\sigma_2(x_0)\not\in T^*{_f(x_0)}Y$,
and hence Lemma \ref{lem:vanDir} implies 
$\Dr f_{!!}i_{1!!}\kn L_{{\sigma}_2}(Z_1\cap Z_2,Z_1)\simeq0$.

\medskip
\noindent
(c)\quad
Therefore we may assume that 
${\sigma}_1(x_0)\neq 0$ and ${\sigma}_2(x_0)\neq 0$.
If $x_0\not\in N$, then the result follows from (iv)
and Lemma \ref{lem:vandir}.
We may assume  therefore $x_0\in N$.
Similarly to the proof of  Lemma \ref{lem:vandir},
we first reduce the proof to the case 
where $X$ is of relative dimension one over $Y$. 
Assume the theorem to be true in the relative one-dimensional morphism
case. 
Set $E=T_{x_0}(f^{-1}f(x_0))$.
Let us choose a line
$\ell\subset E$ such that $\sigma_1(x_0)\vert_\ell\not=0$, and then
decompose $f$ into
$X\lra[g] Y'\lra[h] Y$ on a neighborhood of $x_0$
such that $g$ and $h$ are smooth,
and $T_{x_0}(g^{-1}g(x_0))=\ell$.
Then $g$ satisfies the conditions (i)--(iv),
and applying the relative dimension one case to $g$,
we obtain
$$\Dr f_{!!}\left(\kn L_{\sigma_1}(Z_1,X)\otimes
  \kn L_{\sigma_2}(Z_2,X)\right)\simeq
  \Dr h_{!!}\kn L_{h^*\sigma}\left(g(Z_1\cap Z_2),Y'\right)
  \simeq\kn L_{\sigma}(Z,Y),$$
where the last isomorphism is deduced from Corollary \ref{DirIm}.

Hence me may assume that the relative dimension of $X$ over $Y$ is one. 
By the assumption (iii), $Z_k\ra Y$ is a (local) embedding
and $T_{x_0}Z_k=f_*^{-1}\bl(T_{f(x_0)}Z'_k\br)\cap \sigma_k(x_0)^{-1}(0)$
where $Z'_k\seteq f\left(Z_k\right)$. 
Then $Z'_1$ and $Z'_2$ are transversal submanifolds of $Y$ and 
$Z$ is a one-codimensional submanifold of $Z'\seteq Z'_1\cap Z'_2$.
We have 
$$\sigma(f(x_0))\notin T^*_{Z'}Y.$$
Indeed, we have
\begin{align*}
T_{x_0}(Z_1\cap Z_2)&=f_*^{-1}\bl(T_{f(x_0)}Z'\br)
\cap \sigma_1(x_0)^{-1}(0)\cap
\sigma_2(x_0)^{-1}(0)\\
&=f_*^{-1}\Bl(T_{f(x_0)}Z'\cap\sigma(f(x_0))^{-1}(0)\Br)
\cap\sigma_1(x_0)^{-1}(0),
\end{align*}
which implies
$T_{f(x_0)}Z=T_{f(x_0)}Z'\cap\sigma(f(x_0))^{-1}(0)\not=T_{f(x_0)}Z'$.

Hence we can take local coordinates
$(t,y_1,y_2,z)\in\R\times \R^{m_1}\times\R^{m_2}\times \R^{n}$ of $Y$
such that ${\sigma}(f(x_0))=-\d t(f(x_0))$ and
$Z'_k=\{y_k=0\}$ ($k=1,2$).
Then we can choose a system of coordinates $(x,t,y_1,y_2,z)$ on $X$ 
such that $f$ is given by forgetting $x$,
 ${\sigma}_1(x_0)=-\d x(x_0)$ by (iii)
(and hence $\sigma_2(x_0)=\d x(x_0)-\d t(x_0)$) and that
$Z_1=\{y_1=0,\,x=0\}$.
Set $Z_2=\{y_2=0,\,x=\varphi(t,y_1,z)\}$.
Replacing $\varphi(t,y_1,z)$ with $t$,
we may assume from the beginning
that 
$$\text{$Z_2=\{y_2=0,\,x=t\}$ and $Z=\{y_1=0,\,y_2=0,\,t=0\}$.}$$
We have then using Corollary \ref{AffConExp}
\eqn
\kn L_{{\sigma}_1}(Z_1,X)\otimes\wC_{x_0} 
        & \simeq & \wC_{x_0}\otimes
\Dlimind_{\varepsilon>0}\field_{U^1_{\varepsilon}}\otimes
      \beta_X\big(\Dr i_{1*}\omega^{\otimes-1}_{Z_1/X}\big),\\[.5ex]
\kn L_{{\sigma}_2}(Z_2,X)\otimes\wC_{x_0}
&\simeq&\wC_{x_0}\otimes\Dlimind_{\varepsilon>0} 
\field_{U^2_{\varepsilon}}
\otimes\beta_X\big(\Dr i_{2*}\omega^{\otimes-1}_{Z_2/X}\big),
\eneqn
where the open sets $U^k_{\varepsilon}$ are given by
\eqn
&&U^1_{\varepsilon}=\left\{
\varepsilon |y_1|<x\right\}\quad\text{and}\quad
U^2_{\varepsilon}=\{\varepsilon |y_2|<t-x\}.
\eneqn
We may therefore write
\begin{eqnarray*}
&&\kn L_{{\sigma}_1}(Z_1,X)\otimes\kn L_{{\sigma}_2}(Z_2,X)\otimes
\wC_{x_0}\\
& &\hs{7em} \simeq 
\wC_{x_0} \otimes
\Dlimind_{\varepsilon>0}\field_{U^1_{\varepsilon}\cap U^2_{\varepsilon}}
\otimes\beta_X\big(\Dr i_{1*}
\omega^{\otimes-1}_{Z_1/X}\big)\otimes\beta_X\big(\Dr i_{2*}
\omega^{\otimes-1}_{Z_2/X}\big)\\
&&\hs{7em} \simeq  \wC_{x_0} \otimes
\Dlimind_{\varepsilon>0}\field_{U^1_{\varepsilon}
\cap U^2_{\varepsilon}}\otimes\beta_X\big(f^{-1}\Dr i_{*}
\omega_{Z/Y}^{\otimes-1}\big) \otimes\omega_{X/Y}^{\otimes-1}.
\end{eqnarray*}
Since the relative dimension of $X$ over $Y$ is one, we have 
$\omega_{X/Y}^{\otimes-1}\otimes\wC_{x_0}\simeq\wC_{x_0}\left[1\right],$
and we deduce an isomorphism
\begin{eqnarray*}
&&\Dr f_{!!}\Big(\kn L_{{\sigma}_1}(Z_1,X)\otimes
\kn L_{{\sigma}_2}(Z_2,X)\otimes\wC_{x_0}\Big)\\
&&\hs{7em} \simeq 
\Dr f_{!!}\Big(\wC_{x_0} \otimes\Dlimind_{\varepsilon>0}
\field_{U^1_{\varepsilon}\cap U^2_{\varepsilon}}
\otimes\omega_{X/Y}\Big) \otimes\beta_Y
\Dr i_{*}\omega^{\otimes-1}_{Z/Y}\\
&&\hs{7em}\simeq 
\Dr f_{!!}\Big(\Dlimind_{\varepsilon>0}
\field_{U^1_{\varepsilon}\cap U^2_{\varepsilon}}\Big)
\left[1\right]\otimes\wC_{f(x_0)} \otimes\beta_Y\Dr i_{*}\omega^{\otimes-1}_{Z/Y}.
\end{eqnarray*}
Since
$U^1_{\varepsilon}\cap U^2_{\varepsilon}
=\{\varepsilon|y_1|<x<t-\varepsilon|y_2|\}$,
we have
\begin{eqnarray*}
        \Dr f_{!}\bl(\field_{U^1_{\varepsilon}\cap U^2_{\varepsilon}}\br) &
    \simeq & \field_{\{\varepsilon(|y_1|+|y_2|)<t\}}[-1].
\eneqn
Hence we finally deduce that
\begin{eqnarray*}
\Dr f_{!!}\Big(\kn L_{{\sigma}_1}(Z_1,X)\otimes
\kn L_{{\sigma}_2}(Z_2,X)\otimes\wC_{x_0}\Big)& \simeq &
\Big(\Dlimind_{\varepsilon>0} 
\field_{\{\varepsilon (|y_1|+|y_2|)<t\}}\Big)
\otimes\beta_Y\Dr i_*\omega^{\otimes-1}_{Z/Y}\otimes \wC_{f(x_0)}\\
 & \simeq & \kn L_{\sigma}(Z,Y)\otimes\wC_{f(x_0)}.
\end{eqnarray*}
\end{proof}

\begin{prop}\label{prop:kernel}
Let $(X_1,X_2,X_3)$ be a triplet of manifolds
and $(X_i\times X_j,Z_{ij},\sigma_{ij})$
be a kernel data for $1\le i<j\le 3$.
Assume that 
$Z_{12}\times X_3$ and $X_1\times Z_{23}$
are transversal in $X_1\times X_2\times X_3$
and that the projections 
$p_{ij}\cl X_1\times X_2\times X_3\to X_i\times X_j$
induce an isomorphism
$Z_{12}\ctimes{X_2}Z_{23}\isoto Z_{13}$.
Let us denote by $p_2\cl X_1\times X_2\times X_3\to X_2$ the second projection
and by $p_2{}_*\cl T^*(X_1\times X_2\times X_3)\to T^*X_2$ the induced projection.
Let $N\subset\tan{\sigma_{12}}\ctimes{X_2}\tan{\sigma_{23}}$
be a closed subset satisfying the following conditions:
\bi
\item
$\zeroset{\sigma_{12}}\ctimes{X_2}\zeroset{\sigma_{23}}\subset N$,
\item
$p_{13}^*\sigma_{13}(x)=p_{12}^*\sigma_{12}(x)+p_{23}^*\sigma_{23}(x)$ 
for every $x\in N$,
\item
$p_2{}_*\sigma_{12}(x)\not\in T^*_{X_2}X_2$ for 
any $x\in N\setminus\bigl(\zeroset{\sigma_{12}}\times X_3\cup
X_1\times\zeroset{\sigma_{23}}\bigr)$,
\item
$\R_{\ge0}p_2{}_*\sigma_{12}(x)\not=
\R_{\le0}p_2{}_*\sigma_{23}(x)$ for every 
$x\in \bigl(\tan{\sigma_{12}}\ctimes{X_2}\tan{\sigma_{23}}\bigr)\setminus N$,
\item
the morphism $Z_{12}\ra X_1$ is smooth 
at each point of $\zeroset{\sigma_{12}}$ 
and
the morphism $Z_{23}\ra X_3$ is smooth 
at each point of $\zeroset{\sigma_{23}}$.
\ei
Then we have an isomorphism
$$\kn L_{\sigma_{12}}(Z_{12},X_1\times X_2)
\circ \kn L_{\sigma_{23}}(Z_{23},X_2\times X_3)
\isoto \kn L_{\sigma_{13}}(Z_{13},X_1\times X_3)
\otimes\wC_{f(N)}.$$
\end{prop}
\begin{proof}
By Proposition \ref{InvIm},
we have
\eqn
&&p_{12}{}^{-1}\kn L_{\sigma_{12}}(Z_{12},X_1\times X_2)
\simeq \kn L_{p_{12}^*\sigma_{12}}(Z_{12}\times X_3,X_1\times X_2\times X_3),\\
&&p_{23}{}^{-1}\kn L_{\sigma_{23}}(Z_{23},X_2\times X_3)
\simeq \kn L_{p_{23}^*\sigma_{23}}(X_1\times Z_{23},X_1\times X_2\times X_3),
\eneqn
and Proposition \ref{DirImTens} implies
\eqn
&&\Dr p_{13}{}_{!!}
\left(\kn L_{p_{12}^*\sigma_{12}}(Z_{12}\times X_3,X_1\times X_2\times X_3)
\otimes\kn L_{p_{23}^*\sigma_{23}}(X_1\times Z_{23},X_1\times X_2\times X_3)
\right)\\
&&\hs{50ex}\simeq
\kn L_{\sigma_{13}}(Z_{13},X_1\times X_3)
\otimes\wC_{f(N)}.
\eneqn
\end{proof}

\section{Microlocalization of ind-sheaves}\label{sec:mic}

\subsection{The kernel $\on{K}_\X$ of ind-microlocalization}\unboldmath


We shall construct the kernel of microlocalization by 
the methods of the preceding section using the fundamental $1$-form 
$\liouville_X$ of $T^*X$. 
Since the construction uses only a $1$-form,
we shall discuss it on homogeneous symplectic manifolds.
A {\em homogeneous symplectic manifold} is a manifold $\X$ 
of even dimension 
endowed with a $1$-form $\lio$
such that $(d\lio)^{\dim\X/2}$ never vanishes.
It is a classical result that there locally exists 
a coordinate system $(x_1,\ldots,x_n;\xi_1,\ldots,\xi_n)$
where $\lio$ does not vanish and 
\eq
\lio=\sum_{i=1}^n\xi_idx_i.
\eneq

Let $p_i\cl \X\times\X\to\X$ ($i=1,2$) be the projection and let 
$\Delta_{\X}$ denote the diagonal of $\X\times \X$.
Then $\sigma_\X=p_1^*\lio-p_2^*\lio$
gives a section of $T^*_{\Delta_\X}(\X\times\X)\to\Delta_\X$.

\begin{defn}
The  microlocalization kernel is the kernel defined 
on $\X\times \X$ by: 
   $$\on{K}_\X=\kn L_{\sigma_\X}(\Delta_{\X},\X\times \X)
\in\DI{\X\times\X}.$$
\end{defn}

\begin{lem}\label{damorphism}
There is a natural morphism
 $$\varepsilon_\X\cl \field_{\Delta_{\X}}\lra \on{K}_\X $$
 such that the compositions
\begin{eqnarray*}
  \on{K}_\X\simeq \on{K}_\X\circ \field_{\Delta_{\X}} &
 \lra[{\on{K}_\X\circ \varepsilon_\X}]& 
     \on{K}_\X\circ \on{K}_\X,\\
   \on{K}_\X\simeq\field_{\Delta_{\X}}\circ \on{K}_\X & 
   \lra[{\varepsilon_\X\circ \on{K}_\X}]& 
      \on{K}_\X\circ \on{K}_\X  
\end{eqnarray*}
are isomorphisms, and these two isomorphisms coincide.
\end{lem}
\begin{proof}
We have constructed the morphism $\varepsilon_\X$ in Corollary \ref{TrivMor}. 
The second statement easily follows from
Proposition \ref{prop:kernel}.
The last statement follows from
Lemma \ref{lem:monoid} below
\end{proof}

\begin{lem}\label{lem:monoid}
Let $F\cl\cat{C}\ra\cat{C}$ be a functor and 
$\alpha\cl\id_{\cat{C}}\ra F$ a morphism of functors. 
Assume that for any object $X\in\on{Ob}(\cat{C})$
the morphisms 
$$\alpha_{F(X)}\cl F(X)\ra F(F(X))\qquad\qquad F(\alpha_X)\cl F(X)\ra F(F(X))$$
        are isomorphisms. Then
\bi
\item{ For any two objects $X,Y\in\on{Ob}(\mc{C})$, the composition with 
$\alpha_X$ defines a bijection
   $$\Hom_{\cat{C}}(F(X),F(Y))\isoto\Hom_{\cat{C}}(X,F(Y)),$$}
\item{$\alpha_{F(X)}=F(\alpha_X)$ for any $X\in\on{Ob}(\cat{C})$.}
\ei
\end{lem}

\begin{lem}\label{lem:boxtens}
For two homogeneous symplectic manifolds
$\X$ and $\Y$, we have
$$\text{$\on{K}_{\X\times\Y}\circ(\on{K}_\X\boxtimes\on{K}_\Y)
\simeq \on{K}_{\X\times\Y}$\quad and \quad
$\on{K}_{\X\times\Y}\circ \on{K}_\X\simeq \on{K}_{\X\times\Y}$.}
$$
\end{lem}
\begin{proof}
The last isomorphism is obtained by
applying Proposition \ref{prop:kernel}
to $(\X\times\Y\times \Y,\X,\X)$, and the first isomorphism follows 
from the second since
$$\on{K}_{\X\times\Y}\circ(\on{K}_\X\boxtimes\on{K}_\Y)
\simeq (\on{K}_{\X\times\Y}\circ \on{K}_\X)\circ \on{K}_\Y.$$
\end{proof}

Now let $X$ be a manifold and set $\X\seteq T^*X$.
Then $\X$ has a canonical structure of a homogeneous symplectic manifold.
The microlocalization functor is defined by:
$$\mu_X\cl\Db(\I(\field_X))\lra\Db(\I(\field_\X))\quad 
;\quad\f F\mapsto\mu_X\f F\seteq
  \on{K}_\X\circ\pi^{-1}_X\f F.$$

The microlocalization functor $\mu_X$ may 
also be obtained as an integral transform associated with a kernel 
$\on{L}_X\in\DI{T^*X\times X}$ which is often easier to manipulate
 than $\on{K}_\X$. 

\begin{defn}\label{sigma}
The kernel $\on{L}_X\in\DI{\field_{T^*X\times X}}$ is given by 
$$\on{L}_X=\kn L_{\sigma_X}\bigl(T^*X\ctimes{X}X,T^*X\times X\bigr),$$
where $\sigma_X$ is induced by $\liouville_X$ on the first factor and 
$-id$ 
on the second factor.
\end{defn}
\begin{rem}
Let $(x;\xi)$ be a local coordinate system on $\X=T^*X$ and let 
$(x,\xi;\eta,y)$ denote the associated coordinates on $T^*\X$. 
Then $\sigma_X$ is defined by
$$ \sigma_X(x;\xi)=((x,\xi;\xi,0),(x;-\xi))\in T^*\X\times T^*X. $$
Therefore $\tan{\sigma_X}=T^*X\ctimes{X}X$.
\end{rem}

\begin{prop}
        Let $\f F\in\Db(\I(\field_X))$. There is a canonical isomorphism 
        $$\mu_X\f F\simeq L_X\circ\f F.$$
\end{prop}

\begin{proof}
Consider the following diagram
\eq\label{eq:dg}
&&\ba{c}
\xymatrix@R=.6cm{{} & {T^*X\times T^*X}\ar[d]^{q}\ar[ldd]_{p_1}
   \ar[rd]^{p_2} & {}\\
   {} & {T^*X\times X}\ar[ld]^{p'_1}\ar[rd]_{p'_2} & T^*X\ar[d]^{\pi_X}\\
    {T^*X} & {} & {X.}}\ea\eneq
Since $q$ satisfies the assumptions of Corollary \ref{DirIm}, we have the 
isomorphism
$ \Dr q_{!!}\on{K}_\X\simeq \on{L}_X$,
which implies
\begin{align*} \on{L}_X\circ\f F \simeq \Dr p'_{1!!}\left(
\Dr q_{!!}\on{K}_\X\otimes{p'}^{-1}_2
\f F\right)&\simeq \Dr p'_{1!!}\Dr q_{!!}\left(\on{K}_\X\otimes 
q^{-1}p^{\prime-1}_2\f F\right)\\
&\simeq \Dr p_{1!!}\left(\on{K}_\X\otimes 
p_2^{-1}\pi_X^{-1}\f F\right)\simeq
\on{K}_\X\circ\pi^{-1}_X\f F\simeq\mu_X\f F.
\end{align*}
\end{proof}

The next lemma immediately follows from Lemma \ref{damorphism}.

\begin{lem}\label{lem:Kmu}
For $\f F\in \DI{X}$, we have
$$\on{K}_{T^*X}\circ\mu_X\f F\simeq\mu_X\f F.$$
\end{lem}

\begin{ex}
Let $Z\subset X$ be a closed submanifold. Then
$$ \mu_X(\field_Z) \simeq
      \kn L_{\omega_X}(T^*X\ctimes{X}Z,T^*X). $$
Indeed, noticing that $\field_Z\simeq\kn L_{0}(Z,X)$,
it is enough to apply Proposition \ref{prop:kernel}
to the triplet $(T^*X,X,\on{pt})$ with $N=T^*X\ctimes{X}Z$.

Note that the support of $\mu_X(\field_Z)$ is $T^*_ZX$.
Let us take a local coordinate  system $(x,z)$ on $X$ such that $Z=\{x=0\}$.
Let $(x,z;\xi,\zeta)$ be the corresponding coordinates on $T^*X$. 
Then on $\bdot T^*X$, we have
$$ \mu_X(\field_Z)\simeq
    \Dlimind_{\varepsilon>0}\field_{\{-\langle\xi,x\rangle >
     \varepsilon |x|\}}\otimes \wC_{\{x=0,\,\zeta=0\}}[\codim Z].$$
Note that
\eq
&&\mu_X(\wC_Z)\simeq\wC_{T^*_XX\ctimes{X}Z}.
\label{eq:mutild}
\eneq
\end{ex}

\begin{lem}\label{lem:muzero}
Let $\f F\in\Db\left(\I\left(\field_{T^*X}\right)\right)$. Then 
$$ (\on{K}_{T^*X}\circ\me{F})\otimes\wC_{T^*_XX}\simeq 
   \me{F}\otimes\wC_{T^*_XX}, $$
In particular if $\f F\in\Db\left(\I\left(\field_{X}\right)\right)$ then
$$\mu_X\f F\otimes\wC_{T^*_XX}\simeq\pi^{-1}_X\f F\otimes
\wC_{T^*_XX}.$$
\end{lem}
\begin{proof}
With the notations in \eqref{eq:dg},
we have an isomorphism by Proposition \ref{DG1}:
$$\on{K}_{T^*X}  \otimes p^{-1}_1\wC_{T^*_XX}
\simeq\field_{\Delta_{T^*X}}\otimes p^{-1}_1\wC_{T^*_XX}.$$
Therefore we have for $\f F\in\Db\left(\I\left(\field_{T^*X}\right)\right)$
\begin{align*}
\bigl(\on{K}_{T^*X}\circ\f F\bigr)&
\otimes\wC_{T^*_XX} \simeq \Dr p_{1!!}\left(\on{K}_{T^*X}
    \otimes p^{-1}_2\f F\right)\otimes \wC_{T^*_XX}=
    \Dr p_{1!!}\left(\on{K}_{T^*X}\otimes p^{-1}_1\wC_{T^*_XX}
    \otimes p^{-1}_2\f F\right)\\
        & \simeq \Dr p_{1!!}\left(\field_{\Delta_{T^*X}}\otimes p^{-1}_1
    \wC_{T^*_XX}\otimes p^{-1}_2\f F\right)=
     \Dr p_{1!!}\left(\field_{\Delta_{T^*X}}\otimes p^{-1}_2\f F\right)
     \otimes\wC_{T^*_XX}\\
        & \simeq \f F\otimes\wC_{T^*_XX}.
\end{align*}
\end{proof}
\begin{rem}\label{Conical}
  The ind-sheaf $\mu_X\me{F}$ is conical in the sense that
it is equivariant with respect to the $\R_{>0}$-action
on $T^*X$. We will not develop here 
the theory of conic ind-sheaves but simply give some  consequences
sufficient for our purpose.
  Let $\bdot T^*X$ be the cotangent bundle with its zero section removed,
  and $S^*X$ the associated sphere bundle. 
Let $\gamma\cl\bdot T^*X\ra S^*X$ be the natural projection and
  $\me{F}\in\Db(\I(\field_X))$. Then we have the following isomorphism:
$$ \mu_X\me{F}|_{\bdot T^*X}\simeq 
\gamma^{-1}\Dr\gamma_*\mu_X\me{F}|_{\bdot T^*X}. $$
  Indeed, the kernel $\on{L}_X$ satisfies a similar property.
\end{rem}

\begin{lem}\label{BundLem}
Let $X$ be a real manifold and $\pi_E\cl E\ra X$ a real vector bundle over 
$X$. Denote by $\Sph E$ the spherical bundle associated with $E$ and by
$$j\cl\bdot E\hookrightarrow E\qquad p\cl\bdot E\ra\Sph E$$
the natural morphisms. 
Assume that $\f F\in\DI{E}$ satisfies
$j^{-1}\f F\simeq p^{-1}\f G$ for some
$\f G\in\DI{SE}$. Then 
\bi
\item
$\Dr{\pi_E}_*\Dr j_{!!}j^{-1}\f F\simeq0$,
\item
$\Dr\pi_E{}_*(\f F)\isoto\Dr\pi_E{}_*(\wC_{X}\otimes\f F)$,
where $X$ is identified to the zero section of $E$,
\item
there is a natural distinguished triangle
$$\Dr\bdot{\pi}_{E}{}_{!!}j^{-1}\f F\lra
\Dr \pi_{E}{}_{!!}\f F\lra\Dr {\pi_{E}}{}_*\f F\lra[+1].$$
\ei
\end{lem}

\begin{proof}
(a)\quad
Let $E_X$ denote the real blow up of $E$ along $X$ 
identified with the zero section, {\em i.e.}
$E_X=\bl(\bdot{E}\times\R_{\ge0}\br)/\R_{>0}$,
hence $E_X=\bdot E\sqcup SE$ as a set.
We have the following commutative diagram
$$\xymatrix{{\bdot E}\ar@<-1pt>@{ (->}[r]^{i}
    \ar@{ (->}[rd]_{j}  & {E_X}\ar[r]^{q}\ar[d]^{\pi_{E_X}} &   
   {\Sph E}\ar[d]^{\pi_{\Sph E}}\\
        {} & {E}\ar[r]_{\pi_E} & {X}}$$
where $\pi_{E_X}$ and $\pi_{SE}$ are proper.

\medskip
\noindent
(b)\quad
We shall first show 
$$\Dr q_*\Dr i_{!!}j^{-1}\f F\simeq0.$$
Since $q$ is locally trivial with fiber
$R_{\ge0}$,
we have $q^!\f G\simeq q^{-1}\f G\otimes q^!\field_{SE}
\simeq q^{-1}\f G\otimes\field_{i(\bdot E)}[1]$.
Therefore we have 
$$\Dr q_*(\field_{i(\bdot E)}\otimes q^{-1}\f G)
  \simeq \Dr q_*\Dr\IndHomo\left(\field_{E_X}[1],q^!\f G\right)\simeq
  \Dr\IndHomo\left(\Dr q_{!!}\field_{E_X}[1],\f G\right)\simeq 0$$
since $\Dr q_{!!}\field_{E_X}=0$. 
On the other hand, we have
$$\Dr q_*\left((\field_{i(\bdot E)}/\wC_{i(\bdot E)})\otimes
q^{-1}\f G\right)\simeq
\Dr q_{!!}\left((\field_{i(\bdot E)}/\wC_{i(\bdot E)})\otimes
q^{-1}\f G\right)\simeq
\Dr q_{!!}\left((\field_{i(\bdot E)}/\wC_{i(\bdot E)})\right)\otimes
\f G\simeq0.$$
Hence the desired result follows from the distinguished triangle:
$$\Dr q_*(\wC_{i(\bdot E)}\otimes q^{-1}\f G)
\lra\Dr q_*(\field_{i(\bdot E)}\otimes q^{-1}\f G)
\lra\Dr q_*\left((\field_{i(\bdot E)}/\wC_{i(\bdot E)})\otimes
q^{-1}\f G\right)\lra[+1],$$
in which the first term is isomorphic to $\Dr q_*\Dr i_{!!}j^{-1}\f F$.

\medskip
\noindent
(i)\quad 
We have a chain of isomorphisms
\begin{align*}
\Dr \pi_E{}_*\Dr j_{!!}j^{-1}\f F
&\simeq\Dr \pi_E{}_*\Dr \pi_{E_X}{}_{!!}\Dr i_{!!}j^{-1}\f F\\
&\simeq\Dr \pi_E{}_*\Dr \pi_{E_X}{}_{*}\Dr i_{!!}j^{-1}\f F
\simeq\Dr \pi_{SE}{}_*\Dr q_*\Dr i_{!!}j^{-1}\f F,
\end{align*}
which vanishes by (b).

\medskip\noindent
(ii)\quad
Applying the functor
$\Dr \pi_E{}_*(\,\bullet\,\otimes\f F)$ to the distinguished triangle
\eq\label{eq:dt1}
\wC_{\bdot{E}}\lra\field_E\lra\wC_{X}\lra[+1],
\eneq
we obtain the distinguished triangle
$$\Dr \pi_E{}_*(\wC_{\bdot{E}}\otimes\f F)\lra\Dr \pi_E{}_*\f F\lra
\Dr \pi_E{}_*(\wC_{X}\otimes\f F)\lra[+1],$$
in which the first term vanishes by (i).

\medskip\noindent
(iii)\quad
Applying the functor
$\Dr \pi_E{}_{!!}(\,\bullet\,\otimes\f F)$ to the distinguished triangle
\eqref{eq:dt1},
we obtain the distinguished triangle
$$\Dr \pi_E{}_{!!}(\wC_{\bdot{E}}\otimes\f F)\lra\Dr \pi_E{}_{!!}\f F\lra
\Dr \pi_E{}_{!!}(\wC_{X}\otimes\f F)\lra[+1],$$
in which the first term is isomorphic to
$\Dr\bdot{\pi}_{E}{}_{!!}j^{-1}\f F$
and the last term is isomorphic to $\Dr {\pi_E}{}_*\f F$ by (ii).
\end{proof}

\begin{prop}\label{prop:mudir}
Let $\f F\in\Db\left(\I\left(\field_X\right)\right)$. Then 
\bi
\item{$\Dr{\pi_X}_*\mu_X\f F\simeq\f F$,}
\item{$\Dr{\pi_X}_{!!}\mu_X\f F\simeq\wC_{\Delta_X}
                        \circ\f F$,}
\item{${\Dr\bdot{\pi}_X}_{!!}\left(\mu_X\f F|_{\bdot T^*X}\right)
                  \simeq \left(\field_{X\times X\setminus \Delta_X}\otimes
                       \wC_{\Delta_X}\right)\circ\f F$,}
                \item{there is a natural distinguished triangle
$${\Dr\bdot{\pi}_X}_{!!}\left(\mu_X\f F|_{\bdot T^*X}\right)
       \lra\Dr{\pi_X}_{!!}\mu_X\f F\lra \f F\lra[+1].$$}
\ei
\end{prop}

\begin{proof}
(i)\quad
By Lemma \ref{BundLem} (ii), we have
$$\Dr\pi_{X*}\mu_X\f F\simeq
\Dr\pi_{X*}\left(\mu_X\f F\otimes\wC_{T^*_XX}\right)\simeq
  \Dr\pi_{X!!}\left(\pi^{-1}_X\f F\otimes\wC_{T^*_XX}\right)
  \simeq\f F\otimes\Dr\pi_{X!!}\wC_{T^*_XX}\simeq \f F,$$
where the second isomorphism follows from Lemma \ref{lem:muzero}.

\medskip\noindent
(ii) and (iii)\quad
Let us denote by $p\cl T^*X\times X\to X\times X$ the canonical morphism.
Then we have isomorphisms:
\eqn
\Dr\pi_{X!!}\,\mu_X\f F&\simeq&\left(\Dr p_{!!}L_X\right)\circ\f F,\\
\Dr\bdot\pi_{X!!}(\mu_X\f F|_{\bdot T^*X})
&\simeq&\bl(\Dr p_{!!}(L_X\otimes\wC_{\bdot T^*X\times X})\br)\circ\f F.
\eneqn
Hence, it is enough to show the isomorphism 
\eq
\Dr p_{!!}L_X&\simeq& \wC_{\Delta_X},\label{eq:Ldir}\\
\Dr p_{!!}(L_X\otimes\wC_{\bdot T^*X\times X})
&\simeq&\field_{X\times X\setminus\Delta_X}
\otimes\wC_{\Delta_X}.\label{eq:Ldir2}
\eneq

The natural morphism given in Corollary \ref{TrivMor}
$$\on{L}_X\lra\wC_{T^*X\ctimes{X}X}\otimes\beta_{T^*X\ktimes X}
   \left(\omega^{\otimes-1}_{T^*X\ctimes{X} X/T^*X\ktimes X}\right)=
   p^!\wC_{\Delta_X}$$
provides a morphism $\Dr p_{!!}L_X\lra\wC_{\Delta_X}$.

We shall first show \eqref{eq:Ldir2}.
Take a local coordinate system $x=(x_1,\ldots,x_n)$
on $X$ and let 
$((x;\xi),x')$ be the associated local coordinates on $T^*X\ktimes X$. We have
\begin{eqnarray*}
\on{L}_X\otimes\wC_{\bdot T^*X\ktimes X}&\simeq & \Dlimind_{\varepsilon>0}
\field_{\big\{((x;\xi),x')\; ; \; \langle\xi,x'-x\rangle >
      \varepsilon |x'-x|\big\}}\otimes\wC_{\bdot T^*X\ctimes{X}X} 
\otimes\beta\left(\Dr i_*\omega^{\otimes-1}_{T^*X
       \ctimes{X}X/T^*X\ktimes X}\right) \\
& \simeq & \Dlimind_{\varepsilon>0}
\field_{\big\{((x;\xi),x')\; ; \; \langle \xi,x'-x\rangle >\varepsilon |x'-x|
      \big\}}\otimes p^{-1}\wC_{\Delta_X}[n].
\end{eqnarray*}
Hence 
$$\Dr p_{!!}\left(\on{L}_X\otimes\wC_{\bdot T^*X\ktimes X}\right)\simeq
\Dr p_{!!}\Big(\Dlimind_{\varepsilon>0,\,R>0} 
\field_{\big\{((x;\xi),x')\; ; \; 
\langle \xi,x'-x\rangle >\varepsilon |x'-x|,\,|\xi|<R\big\}}\Big)
\otimes \wC_{\Delta_X}[n].$$
For $0<\varepsilon<R$, we have
$$\Dr^k p_{!}\Big(\field_{\big\{((x;\xi),x')\; ; \; 
   \langle \xi,x'-x\rangle >\varepsilon |x'-x|,\,|\xi|<R\big\}}\Big)
   \simeq
\begin{cases}
        \field_{\{0<|x'-x|<\varepsilon^{-1}R\}} & \textrm{if $k=n$.}\\
        0 & \textrm{if $k\neq n$.}
\end{cases}$$
Hence we have shown that
\begin{eqnarray*}
\Dr p_{!!}\left(\on{L}_X\otimes\wC_{\bdot T^*X\ktimes X}\right) & \simeq & \field_{X\ktimes X\setminus\Delta_X}[-   n]\otimes\wC_{\Delta_X}[n]\simeq\field_{X\ktimes X\setminus\Delta_X}\otimes\wC_{\Delta_X},
\end{eqnarray*}
which proves \eqref{eq:Ldir2}.
In the morphism of distinguished triangles
$$\xymatrix{
\Dr p_{!!}\left(\on{L}_X\otimes\wC_{\bdot T^*X\ktimes X}\right)
\ar[r]\ar[d]^\sim&
\Dr p_{!!}\left(\on{L}_X\right)
\ar[r]\ar[d]&{}
\Dr p_{!!}(\left(\on{L}_X\otimes\wC_{T^*_XX\ktimes X}\right)
\ar[r]^(.8){+1}\ar[d]&{\strut}\\
\field_{X\ktimes X\setminus\Delta_X}\otimes\wC_{\Delta_X}
\ar[r]
&\wC_{\Delta_X}\ar[r]&
\field_{\Delta_X}\ar[r]^{+1}&{\strut},
}
$$
the left vertical arrow is an isomorphism by \eqref{eq:Ldir2}
and the right vertical arrow is an isomorphism since
$$\Dr p_{!!}\left(\on{L}_X\otimes\wC_{T^*_XX\ktimes X}\right)
\simeq \Dr p_{!!}\left(\field_{{T^*X\ctimes{X}X}}
\otimes\wC_{T^*_XX\ktimes X}\right)
\simeq
\field_{\Delta_X}\otimes\Dr p_{!!}(\wC_{T^*_XX\ktimes X})
\simeq \field_{\Delta_X}.$$
Hence we obtain \eqref{eq:Ldir}.

\medskip
\noindent
(iv) follows immediately from Lemma \ref{BundLem} and (i).
\end{proof}

\begin{prop}\label{boxcomposition}
For $\f F\in\Db\left(\I\left(\field_X\right)\right)$ 
and $\f G\in\Db\left(\I\left(\field_Y\right)\right)$,
we have an isomorphism
        $$\mu_{X\times Y}\left(\f F\boxtimes\f G\right)\simeq 
\K_{T^*(X\times Y)}\circ\left(\mu_X\f F\boxtimes\mu_Y\f         G\right).$$
\end{prop}

\begin{proof}
This follows immediately from Lemma \ref{lem:boxtens}.
\end{proof}

\boldmath\subsection{The link with $\muHom$ and classical microlocalization}
\unboldmath


\begin{prop}\label{FundProp}
Let $\sigma\in\Gamma(X,\Omega^1_X)$ and $\f F,\f G\in\Db(\field_X)$. 
Then we have an isomorphism
        $$\sigma^{-1}\muHom(\f F,\f G)\simeq
\Dr\Homo\Big(\f F,\kn L_{\widetilde\sigma}(\Delta_X,X\times X)
\circ\f G\Big),$$
where $\widetilde\sigma=q_1^*\sigma-q_2^*\sigma$ and 
$q_i\cl X\times X\to X$ is the $i$-th projection {\rm($i=1,2$)} .
\end{prop}

\begin{proof}
By definition we have
$$\muHom(\f F,\f G)\simeq\nuHom(\f F,\f G)^{\wedge},$$
where $\nuHom$ is the specialization of the functor $\Dr\Homo$ (see below),
and $(\cdot)^\wedge$ is the Fourier-Sato transform.
Setting
$$
P'=\left\{\left((x;\xi),(x;v)\right)\in T^*X\ctimes{X}TX;
\langle\xi,v\rangle\leqs 0\right\},$$
the Fourier-Sato transform is the integral transform with kernel $\field_{P'}$. 
Consider the following commutative diagram
$$\xymatrix{{TX} & {} \\
*+!<0pc,0.25pc>{T^*X\ctimes{X} TX} \ar[u]^{\pi_2}
\ar[d]_{\pi_1}\ar@{}[rd]|{\square} 
&{TX}\ar[l]^(.4){\sigma'}\ar[d]^{\tau_X}\ar[ul]_{\id_{T_X}} \\
{T^*X} & {X}.\ar[l]^{\sigma}} $$ 
Then 
$\muHom(\f F,\f G)\simeq \nuHom(\f F,\f G)^{\wedge}
\simeq\Dr \pi_{1!}\left(\pi_2^{-1}\nuHom(\f F,\f G)\otimes \field_{P'}\right)$. 
Hence 
\begin{eqnarray*} 
{\sigma}^{-1}\muHom(\f F,\f G) & \simeq & 
{\sigma}^{-1}\Dr\pi_{1!}\left(\pi_2^{-1}
\nuHom(\f F,\f  G)\otimes\field_{P'}\right)\\
&\simeq&\Dr\tau_{X!}{\sigma'}^{-1}
\left(\pi_2^{-1}\nuHom(\f F,\f G)\otimes\field_{P'}\right)
 \simeq  \Dr\tau_{X!} \left(\nuHom(\f F,\f G) \otimes \field_{P'_{\sigma}}\right),
\end{eqnarray*}
where we have set 
$P'_{\sigma}={\sigma'}^{-1}(P')
=\left\{ (x,v)\in TX ; \langle\sigma(x),v\rangle\leqs 0\right\}$. 
Consider the normal deformation of $\Delta_X$ in $X\times X$,
visualized by the diagram:
$$\xymatrix{
{TX}\ar[r]^(.3){\sim}
    &{T_{\Delta_X}(X\times X)}\ar[d]^{\tau_X}\ar@<-1pt>@{ (->}[r]^(.6){s} 
         & {\widetilde{X\times X}}\ar[d]_p\ar@<-2pt>@/_5ex/[dd]_(.4){p_1}   
                         \ar@<+2pt>@/^5ex/[dd]^{p_2} 
   &{\Omega}\ar@<1pt>@{ )->}[l]_(.3){j}|(.3){\circ}\ar[ld]^(.3){\widetilde p}\\
 & {X}\ar@<-1pt>@{ (->}[r]_{i}\ar[dr]_{\id_{_X}} 
  &{X\times X} \ar@<+2pt>[d]^(.4){q_2}\ar@<-2pt>[d]_(.4){q_1} & {} \\
        & {} & {X} &  }$$
Then $\nuHom(\f F,\f G)$ is by definition
$s^{-1}\Dr j_*{\widetilde p}^{-1}\Dr\Homo(q_2^{-1}\f F,q_1^!\f G)$.
Since $\widetilde p$ is smooth we have
$${\widetilde p}^{-1}\Dr\Homo(q_2^{-1}\f F,q_1^!\f G)
\simeq\Dr\Homo({\widetilde p}^{-1}q_2^{-1}\f F,{\widetilde p}^{-1}q_1^!\f G)
\simeq\Dr\Homo(j^{-1}p_2^{-1}\f F,j^{-1}p^{-1}q_1^!\f G).$$
Hence we have
\begin{align*}
\nuHom(\f F,\f G) 
& \simeq  
s^{-1}\Dr j_*\Dr\Homo(j^{-1}p_2^{-1}\f F,j^{-1}p^{-1}q_1^!\f G)
\simeq  s^{-1}\Dr\Homo(p_2^{-1}\f F,\Dr j_*j^{-1}p^{-1}q_1^{!}\f G) \\
& \simeq s^{-1}\Dr\Homo(p_2^{-1}\f F,\Dr j_*j^{-1}p_1^{-1}\f G)
\otimes\tau_X^{-1}\omega_X.
\end{align*}
Since $p_1$ is smooth, we have the estimate 
$$\SS(p^{-1}_1\f G)\cap\SS\field_{\Omega}
\subset (p_1){}_\pi^{-1}(T^*X)\cap
\left(T^*_{T_{\Delta_X(X\times X)}}\widetilde{X\times X}
\cup T^*_{\widetilde{X\times X}}\widetilde{X\times X}\right)
\subset T^*_{\widetilde{X\times X}}\widetilde{X\times X},$$
which implies
\begin{eqnarray*}
        \Dr j_*j^{-1}p_1^{-1}\f G & \simeq &
\Dr\Homo(\field_{\Omega},p_1^{-1}\f G)
 \simeq  \Dr\Homo(\field_{\Omega},\field_X)\otimes p_1^{-1}\f G\\
&\simeq&
\field_{\ol{\Omega}}\otimes p^{-1}_1\f G.
\end{eqnarray*}
Applying this result we obtain
\begin{eqnarray*}
\nuHom(\f F,\f G) & \simeq &
s^{-1}\Dr\Homo(p_2^{-1}\f F,p_1^{-1}\f  G\otimes\field_{\ol{\Omega}})
\otimes \tau_X^{-1}\omega_X\\
& \simeq & s^{-1}\Dr\Homo(p_2^{-1}\f F,p_1^{-1}\f G\otimes\field_{\ol{\Omega}})\otimes \tau_X^{-1}
      \omega_{\Delta_X/X\times X}^{\otimes -1},
\end{eqnarray*}
and finally
\begin{eqnarray*}
{\sigma}^{-1}\muHom(\f F,\f G) & \simeq & 
\Dr\tau_{X!}\left(s^{-1}\Dr\Homo\left(p_2^{-1}\f F,p_1^{-1}\f G\otimes
\field_{\ol{\Omega}}\right)\otimes
\tau_X^{-1}\omega_{\Delta_X/X\times X}^{\otimes -1}
\otimes\field_{P'_{\sigma}}\right)\\
& \simeq & \Dr p_{2}{}_!
\Dr s{}_!\left(s^{-1}\Dr\Homo\left(p_2^{-1}\f F,p_1^{-1}\f G
\otimes\field_{\ol{\Omega}}\right)\otimes
\tau_X^{-1}\omega_{\Delta_X/X\times X}^{\otimes -1}
\otimes\field_{P'_{\sigma}}\right)\\
& \simeq & \Dr p_{2!}\left(
\Dr\Homo\left(p_2^{-1}\f F,p_1^{-1}\f G\otimes\field_{\ol{\Omega}}\right)
\otimes \field_{P'_{\sigma}}\otimes
p^{-1}\Dr i_*\omega_{\Delta_X/X\times X}^{\otimes -1}
)\right).
\end{eqnarray*}
Note that this intermediate result is obtained by means of
 classical sheaf theory.
However, formulas in the derived category of ind-sheaves
allow us to continue the calculations.
Using the properties \eqref{eq:betahom} of the 
functor $\beta$ and Proposition \ref{prop:tenshom}, we have
\begin{eqnarray*}
{\sigma}^{-1}\muHom(\f F,\f G)\kern -0.5em & \simeq &
\kern -0.5em \Dr p_{2!}
\Dr\Homo\Bl(p_2^{-1}\f F,p_1^{-1}\f G\otimes\field_{\ol{\Omega}}
\otimes\beta_{\widetilde{X\times X}}
\bl(\field_{P'_{\sigma}}\otimes
p^{-1}\Dr i_*\omega_{\Delta_X/X\times X}^{\otimes -1}\br)\Br)\\
\kern -0.5em& \simeq &\kern -0.5em \Dr\Homo
\Bl(\f F,\Dr p_{2!!}\bl(p_1^{-1}\f G\otimes\field_{\ol{\Omega}}\otimes      
\wC_{P'_{\sigma}}\otimes
p^{-1}\beta_{X\times X}(\Dr i_*\omega_{\Delta_X/X\times X}^{\otimes -1})\br)
\Br).
\end{eqnarray*}
We have furthermore
\eqn
&&\Dr p_{2!!}\bl(p_1^{-1}\f G\otimes\field_{\ol{\Omega}}\otimes
\wC_{P'_{\sigma}}\otimes
p^{-1}\beta_{X\times X}(\Dr i_*\omega_{\Delta_X/X\times X}^{\otimes -1})\br)\\
&&\hs{15ex}\simeq
\Dr q_2{}_{!!}\Dr p_{!!}
\left(p^{-1}q_1^{-1}\f G\otimes\field_{\ol{\Omega}}\otimes
\wC_{P'_\sigma}\otimes
p^{-1}\beta_{X\times X}(\Dr i_*\omega_{\Delta_X/X\times X}^{\otimes -1})
\right)\\
&&\hs{15ex}\simeq 
\Dr q_{2!!}\left(q_1^{-1}\f G\otimes\Dr p_{!!}
\bl(\field_{\ol{\Omega}}\otimes\wC_{P'_{\sigma}}\br)\otimes
\beta_{{X\times X}}
(\Dr i_*\omega_{\Delta_X/X\times X}^{\otimes-1})\right)\\
&&\hs{15ex}\simeq 
\Dr q_{1!!}\left(q_2^{-1}\f G\otimes\Dr p_{!!}
\bl(\field_{\ol{\Omega}}\otimes\wC_{P_{\widetilde\sigma}}\br)\otimes
\beta_{{X\times X}}
(\Dr i_*\omega_{\Delta_X/X\times X}^{\otimes-1})\right)\\
&&\hs{15ex}\simeq
\Dr q_{1!!}\left(q_2^{-1}\f G\otimes
\kn L_{\widetilde\sigma}(\Delta_X,X\times X)\right)\simeq
\kn L_{\widetilde\sigma}(\Delta_X,X\times X)\circ\f G. 
\eneqn
\end{proof}

\begin{cor}\label{Homformula}
        Let $\f F,\f G\in\Db(\field_X)$. Then we have an isomorphism
  $$\muHom(\f F,\f G)\simeq\Dr\Homo(\pi_X^{-1}\f F,\mu_X\f G)\simeq
       \Dr\Homo(\mu_X\f F,\mu_X\f G).$$
\end{cor}

\begin{proof}
Consider the fundamental $1$-form $\liouville_X\in\Gamma(T^*X,\Omega^1_{T^*X})$ of the cotangent bundle of $X$. Then we have
$$\muHom(\f F,\f G)\simeq\liouville_X^{-1}\muHom(\pi^{-1}_X\f F,\pi^{-1}_X\f G)$$
and by Proposition \ref{FundProp} we get a natural isomorphism
$$\muHom(\f F,\f G)=\Dr\Homo(\pi_X^{-1}\f F,\K_{T^*X}\circ\pi^{-1}_X\f G)\simeq
\Dr\Homo(\pi_X^{-1}\f F,\mu_X\f G)$$
The last isomorphism is a consequence of
Lemma \ref{lem:monoid} 
and Lemma \ref{damorphism}.
\end{proof}

\begin{prop}
Let $\f F\in\Db(\field_X)$ and let $Z$ be a closed submanifold of $X$. 
Denote by $i$ the closed immersion $i\cl T^*X\ctimes{X}Z\hookrightarrow T^*X$. Then we have a natural isomorphism
$$\mu_Z(\f F)\simeq
\alpha_{{}_{{T^*X\ctimes{X}Z}}}\big(i^!\mu_X\f F)|_{T^*_ZX}\big)
\simeq\Dr\me{H}om(\field_{T^*X\ctimes{X}Z},\mu_X\f F)|_{T^*_ZX}.$$
 Here $\mu_Z(\me{F})$ denotes the classical functor of Sato's 
microlocalization 
\end{prop}
See \cite{KS2}, Chapter IV for definitions and a detailed study
for $\mu_Z$. We only remark here that 
$\mu_Z(\f F)\simeq \muHom(\field_Z,\f F)|_{T^*_ZX}$.

\begin{proof}
We have by Corollary \ref{Homformula}
\eqn
\mu_Z(\f F)&\simeq&\Dr\Homo(\pi_X^{-1}\field_{Z},\mu_X(\f F))\vert_{T^*_ZX}
\simeq \Dr\Homo\bigl(\Dr i_{!!}\field_{{T^*X\ctimes{X}Z}},\mu_X(\f F)\bigr)
\vert_{T^*_ZX}\\
&\simeq&\Dr\Homo\bigl(\field_{{T^*X\ctimes{X}Z}},i^{!}\mu_X(\f F)\bigr)
\vert_{T^*_ZX}
\simeq\bl(\alpha_{{}_{T^*X\ctimes{X}Z}}i^{!}\mu_X(\f F)\br)\vert_{T^*_ZX}.
\eneqn
\end{proof}


\subsection{Review on the microsupport of ind-sheaves}

In this section we shall give a short overview on the results of
\cite{KS4} on the 
microsupport of ind-sheaves .

The microsupport $\SS(\me{F})$ of an object $\me{F}\in\Db(\field_X)$ 
is a closed 
involutive cone in the cotangent bundle $T^*X$ which describes 
the codirections in 
which the cohomology of $\me{F}$ does not propagate 
(cf.\ \cite{KS2}, \cite{KS3}).
The corresponding notions for ind-sheaves 
are more intricate.

Let $\mc{C}$ be an abelian category, and consider the functor
$$ \on{J}\cl\Db(\on{Ind}(\mc{C})) \lra \Db(\mc{C})^{\wedge}
  \qquad\text{given by}\qquad \me{F}\mapsto 
\on{Hom}_{\Db(\on{Ind}(\mc{C}))}(\,\cdot\,,\me{F}). $$
Here, $\Db(\mc{C})^{\wedge}$ is the category of contravariant functors
from $\Db(\mc{C})$ to the category of sets.
Then it can be shown that $\on{J}$ factors through $\on{Ind}(\Db(\mc{C}))$. 
Note that $\on{J}$
is conservative, which is a consequence of the commutative diagram
$$ \xymatrix{
     \Db(\on{Ind}(\mc{C})) \ar[rr]^{\on{J}} \ar[dr]_{\on{H}^k} & & 
      \on{Ind}(\Db(\mc{C})) \ar[dl]^{\on{IH}^k} \\
     & \on{Ind}(\mc{C}) & } $$
Finally assume 
\eq\label{eq:fincoh}
&&\text{$\mc{C}$ has enough injectives and finite homological dimension.}
\eneq
Recall that 
in this case $\varphi\cl \me{F}\ra\me{G}$ is an isomorphism in 
$\on{Ind}(\Db(\mc{C}))$ 
if and only if $\on{IH}^k(\varphi)$ is an isomorphism for all $k$. 
Then we easily get the 
following result.

\begin{lem}\label{thefunctorJ}
Assume \eqref{eq:fincoh}.
  Let $\me{F}\in\Db(\on{Ind}(\mc{C}))$ and let
$\{\me{F}_i\to\f F\}_{i\in I}$ be a filtrant inductive system of morphisms
   in $\Db(\on{Ind}(\mc{C}))$. 
Then $\Dlimind \on{J}(\me{F}_i)\isoto\on{J}(\me{F})$
if and only if\/
$\smash{\Dlimind_{i\in I}} \on{H}^k(\me{F}_i) \isoto\on{H}^k(\me{F})$.
\end{lem}

In particular if $\Dlimind \on{J}(\me{F}_i)\isoto\on{J}(\me{F})$ the we have
$\Dlimind \on{J}(\tau^{\leqs n}\me{F}_i)\isoto\on{J}(\tau^{\leqs n}\me{F})$
for all $k$.

We shall apply the results above to the case of ind-sheaves,
by taking $\Mod^{\on{c}}(\field_X)$ as $\mc{C}$.
For a $\on{C}^\infty$-manifold $X$,
let
$$J_X\cl \DI{X}\to\bigl(\Db(\Mod^{\on{c}}(\field_X))\bigr)^\wedge$$
be the canonical functor.

\begin{prop}\label{functorJandcomposition}
Let $f\cl X\ra Y$ be a continuous map. Let $\{\me{F}_i\to\me{F}\}_{i\in I}$ 
be a filtrant inductive system of morphisms
in $\Db(\I(\field_X))$ and
  $\{\me{G}_j\to \me{G}\}_{j\in J}$ a filtrant inductive system in 
$\Db(\I(\field_Y))$
  such that 
 $$J_X(\me{F})\simeq  \Dlimind_{i\in I}J_X(\me{F}_i) \qquad \text{and} \qquad 
   J_Y(\me{G})\simeq  \Dlimind_{j\in J}J_Y(\me{G}_j).   $$
Then 
  \bi
 {\item 
     $$J_Y(\Dr f_{!!} \me{F})\simeq\Dlimind_{i\in I} J_Y(\Dr f_{!!}\me F_i),$$}
  {\item  For $\me{K}\in\Db(\I(\field_X))$, we have
      $$J_X(\me{K}\otimes\me{F})=\Dlimind_{i\in I}J_X
      (\me{K}\otimes \me{F}_i),$$ }
 {\item  
         $$J_X(f^{-1}\me{G})\simeq\Dlimind_{j\in J} J_X(f^{-1}\me{G}_j)
\quad\text{and}\quad
J_X(f^{!}\me{G})\simeq\Dlimind_{j\in J} J_X(f^{!}\me{G}_j)$$ }
 {\item 
$$J_{T^*X}(\mu_X\me{F}) \simeq \Dlimind_{i\in I} J_{T^*X}(\mu_X\me{F}_i).$$ }
  \ei
\end{prop}
\begin{proof}
By Lemma \ref{thefunctorJ}, we can reduce the situation by d{\'e}vissage
to usual ind-sheaves, where the formulas are obvious.
\end{proof}
\begin{defn}
\bi
\item  Let $\me{F}\in\Db(\I(\field_X))$. The micro-support of $\me{F}$,
denoted  $\SS(\me{F})$,
is the closed conic
  subset of $T^*X$ whose complementary is the set of points 
$p\in T^*X$ such that 
  there exist 
a conic open neighborhood $U$ of $p$ in $T^*X$,
an open neighborhood $W$ 
of $\pi_X(p)$ and a small filtrant 
  inductive system $\{\me{F}_i\}_{i\in I}$ of objects 
  $\me{F}_i\in\Db(\on{Mod}^{\on{c}}(\field_X))$
such that $\SS(\me{F}_i)\cap U=\varnothing$ and
$$J_X(\me{F}\otimes\field_W)\simeq 
\Dlimind_{i\in I}{\me{F}_i\otimes\field_W}. $$
\item
For $\f F\in \DI{X}$, one sets $\SS_0(\me{F})=\supp(\mu_X\me{F})$.
\ei
\end{defn}
\begin{rem}
  The micro-support defined above coincide with the classical definition
  for objects of $\Db(\field_X)$, it satisfies the triangular 
  inequality (in a distinguished triangle, the micro-support of an
  object is contained in the union of the micro-supports of the two
  others),  and we have
  \eqn&& 
\ba{c}
\supp(\me{F})=\SS(\me{F})\cap T^*_XX,\quad
\SS(\alpha_X(\f F))\subset\SS(\f F)
\ea
           \quad\text{for $\me{F}\in\Db(\I(\field_X))$.}
\eneqn
In general, it is no longer an involutive subset of $T^*X$.
\end{rem}

\begin{prop}
Let $\f F\in\Db(\I(\field_X))$. Then 
        $$\SS_0(\f F)\subset \SS(\f F).$$
If $\f F\in\Db(\field_X)$, then
        $$\SS_0(\f F)=\SS(\f F).$$
\end{prop}

\begin{proof}
The result for sheaves is actually an obvious consequence of Corollary 
\ref{Homformula} since 
$$\SS(\f F)=\supp(\muHom(\f F,\f F))=\supp(\Dr\Homo(\mu_X\f F,\mu_X\f F))=
  \supp(\mu_X\f F).$$
Now assume that $\f F\in\Db(\I(\field_X))$ 
and $p\not\in\SS \f F$. Consider a filtrant 
inductive system $\f F_i$ in 
$\Db({\rm Mod}^{c}(\field_X))$ and
an open neighborhood $W$ of $\pi_X(p)$, a neighborhood 
$U\subset\pi_{T^*X}{}^{-1}(W)$ of $p$
such that
$$J_X(\f F\otimes\field_W)\simeq\Dlimind_i(\f F_i\otimes\field_W) $$
and $\SS(\me{F}_i)\cap \ol{U}=\varnothing$. We have by Proposition 
\ref{functorJandcomposition}
$$J_X\bl(\mu_X(\f F\otimes\field_W)\br)\simeq 
\Dlimind_i J_X\bl(\mu_X(\f F_i\otimes\field_W)\br),$$
and we get $\mu_X\me{F}|_{U}\simeq 0$ since 
$\supp(\mu_X\me{F}_i)=\SS\me({F}_i)$.
\end{proof}

\begin{ex}
For a closed submanifold $Z$ of $X$, we have
\eqn
&&\SS_0(\field_Z)=\SS(\field_Z)=T^*_ZX\quad\text{and}\\[1ex]
&&\SS_0(\wC_Z)=T^*_XX\ctimes{X}Z,\quad\SS(\wC_Z)=T^*_ZX.
\eneqn
\end{ex}

\begin{lem}\label{LemMS1}
Let $\Omega$ be an open subset of $\bdot{T}^*X$ and 
let $\f F\in\Db(\field_{\Omega})$, $\f G\in\DI{\Omega}$. 
Assume that $\f F$ is cohomologically constructible
{\rm(}see {\rm \cite[Definition 3.4.1]{KS2})}.
Assume further
$$\liouville_X^{-1}\bigl(\SS(\f F)\bigr)\cap\Supp(\f G)
=\varnothing,$$
where $\omega_X$ is considered as a map $T^*X\to T^*(T^*X)$.
        Then we have an isomorphism
        $$\Dr\Homo(\f F,\field_{\Omega})\otimes(\K_{\Omega}\circ\f G)
\isoto\Dr\IndHomo(\f F,\K_{\Omega}\circ\f G)\quad\text{in $\DI{\Omega}$.}$$
\end{lem}

\begin{proof}
By shrinking $\Omega$, we may assume from the beginning that 
$\liouville_X^{-1}\bigl(\SS(\f F)\bigr)=\varnothing$.

(i)\quad Assume first that $\f G\in\Db(\field_\Omega)$.
For $p=(x_0,\xi_0)\in\Omega$, we shall prove that
$$\Dr\Homo(\f F,\field_{\Omega})
\otimes(\K_{\Omega}\circ\f G)\otimes\wC_{p}\isoto
\Dr\IndHomo(\f F,\K_{\Omega}\circ\f G)\otimes\wC_{p}.$$

Since $p\notin T^*_XX$, we have:
\begin{eqnarray}\label{eq:mup}
&&\begin{array}{rcl}
(\K_{\Omega}\circ\f G)\otimes\wC_{p} & \simeq & 
(\K_{\Omega}\otimes\wC_{(p,p)})\circ\f G
     \simeq 
\Dlimind_{\rho>0}\field_{K_\rho}
\otimes\Bl(\bigl(\Dlimind_{\delta>0,\varepsilon>0}
      \field_{F_{\delta,\varepsilon}}\bigr)[-n]\circ\f G\Br),
\end{array}
\end{eqnarray}
where 
$$ K_\rho=\{(x,\xi); |x-x_0|\le\rho,\,|\xi-\xi_0|\le\rho\}$$
and $$F_{\delta,\varepsilon}=
\big\{\delta\geqs\langle\xi_0,x'-x\rangle
>\varepsilon(|x'-x|+|\xi'-\xi|)\big\}.$$
Let $p_1\cl T^*\Omega\times T^*\Omega\to T^*\Omega$ be the first projection.
For sufficiently small $\varepsilon$, $\delta$ and $\rho$, 
$\pi_\X^{-1}K_\rho\cap
p_1\bl(\SS(\field_{F_{\delta,\varepsilon}})\br)$ is contained in 
a sufficiently small neighborhood of 
$\omega_X(p)$, and hence so is 
$\pi_\X^{-1}K_\rho\cap\SS(\field_{F_{\delta,\varepsilon}}\circ\f G)$.
Thus we obtain by assumption
$$\pi_\X^{-1}K_\rho\cap\SS\f F\cap\SS(\field_{F_{\delta,\varepsilon}}\circ\f G)
\subset T^*_{\X}\X.$$
Then by \cite[Corollary 6.4.3]{KS2}, we have an isomorphism
$$\field_{K_\rho}\otimes\Dr\Homo(\f F,\field_{\Omega})
\otimes(\field_{F_{\delta,\varepsilon}}\circ\f G)
\isoto\field_{K_\rho}\otimes\Dr\Homo(\f F,\field_{F_{\delta,\varepsilon}}\circ\f G)$$
in $\Db(\field_{\Omega})$.
Therefore we have
\begin{align*}
  J_{\Omega}\left(\Dr\Homo(\f F,\field_{\Omega})\right.&\left.\otimes(\K_{\Omega}\circ\f G)\otimes\wC_{p}\right) \\
   & \simeq  \hs{-3ex}
  \Dlimind_{\delta>0,\,\varepsilon>0,\,\rho>0}J_{\Omega}
\left(\field_{K_\rho}\otimes\Dr\Homo(\f F,\field_{\Omega})\otimes
   (\field_{F_{\delta,\varepsilon}}\circ\f G)[-n]\right)\\
        & \simeq  \hs{-2ex}
 \Dlimind_{\delta>0,\,\varepsilon>0,\,\rho>0}{J_{\Omega}\Bigl(
\field_{K_\rho}\otimes
\Dr\Homo\bigl(\f F,\field_{F_{\delta,\varepsilon}}\circ\f G[-n]\bigr)\Bigr)}\\
\simeq J_{\Omega}\left(\Dr\IndHomo(\f F,\K_{\Omega}\circ\f G)\otimes\wC_{p}
     \right),
\end{align*}
and the lemma is proved when $\f G\in\Db(\field_\Omega)$.

\medskip
In the general case, taking a filtrant inductive system
$\f{G}_k$ in $\Db(\field_{\Omega})$
such that $J_\Omega(\f{G})\simeq\Dlimind\f{G}_k$. 
we have
\begin{align*}
  J_\Omega\left(\Dr\Homo(\f F,\field_{\Omega})\right.&\left.\otimes(\K_{\Omega}\circ\f G)\right) 
  \simeq 
  \Dlimind_{k}J_\Omega\left(\Dr\Homo(\f F,\field_{\Omega})\otimes
(\K_{\Omega}\circ\f G_k)\right)\\
&\simeq  \Dlimind_{k}J_\Omega
\left(\Dr\Homo(\f F,K_{\Omega}\circ\f G_k)\right) \simeq 
J_\Omega\Bigl(\Dr\Homo(\f F,\K_{\Omega}\circ\f G)\bigr),
\end{align*}
which completes the proof.
\end{proof}

We prove now in the framework of ind-sheaves a well known result for sheaves.

\begin{prop}
        Let $\f F\in\Db(\field_X)$ and $\f G\in\Db(\I(\field_X))$. 
Assume that $\f F$ is cohomologically constructible.
Assume further the non-characteristic condition
        $$\SS(\f F)\cap\SS_0(\f G)\subset T^*_XX.$$
        Then we have an isomorphism
        $$\Dr\Homo(\f F,\field_X)\otimes\f G\isoto\Dr\IndHomo(\f F,\f G).$$
\end{prop}

\begin{proof}
Since $\liouville^{-1}_X\SS(\pi_X^{-1}\f F)=\SS\f F$,
the  non-characteristic condition may be rewritten as
$$\liouville^{-1}_X\SS(\pi_X^{-1}\f F)\cap\supp\mu_X\f G\cap\bdot T^*X=\varnothing,$$
and Lemma \ref{LemMS1} assures that
\eqn
\Big(\pi_X^{-1}\Dr\Homo(\f F,\field_{X})\otimes\mu_X\f G\Big)|_{\bdot T^*X}
&\simeq&
\Big(\Dr\Homo(\pi^{-1}_X\f F,\field_{T^*X})\otimes\mu_X\f G\Big)|_{\bdot T^*X}\\
&\simeq&\Dr\IndHomo(\pi^{-1}_X\f F,\mu_X\f G)|_{\bdot T^*X}.
\eneqn
Applying the functor $\Dr\bdot\pi_{X!!}$, we obtain
$$\Dr\Homo(\f F,\field_X)\otimes\Dr\bdot\pi_{X!!}\left(\mu_X\f G|_{\bdot T^*X}\right)\simeq\Dr\IndHomo\left(\f F,\Dr\bdot\pi_{X!!}\left(\mu_X\f G|_{\bdot T^*X}\right)\right).$$
Now, 
Proposition \ref{prop:mudir} gives
the following morphism of distinguished triangles
where $\f F^*=\Dr\Homo(\f F,\field_X)$:
$$\xymatrix{{\f F^*\otimes\Dr\bdot\pi_{X!!}
\left(\mu_X\f G|_{\bdot T^*X}\right)}\ar[r]\ar[d]^{\sim} 
& {\f F^*\otimes(\wC_{\Delta_X}\circ\f G})\ar[r]\ar[d]
& {\f F^*\otimes\f G}\ar[r]^(.6){+1}\ar[d] & {}\\
{\Dr\IndHomo\left(\f F,\Dr\bdot\pi_{X!!}(\mu_X\f G|_{\bdot T^*X})\right)}
\ar[r] & {\Dr\IndHomo(\f F,\wC_{\Delta_X}\circ\f G)}\ar[r] 
& {\Dr\IndHomo(\f F,\f G)}\ar[r]^(.75){+1} & {}}.$$
The middle vertical arrow is an isomorphism by 
the following lemma, and hence the right
arrow is an isomorphism.
\end{proof}

\begin{lem}
Let $\f F\in\Db(\field_X)$ and $\f G\in\Db(\I(\field_X))$. 
Assume that $\f F$ is cohomologically constructible.
Then we have an isomorphism
$$\Dr\Homo(\f F,\field_X)\otimes\bl(\wC_{\Delta_X}\circ\f G\br)
\isoto\Dr\IndHomo(\f F,\wC_{\Delta_X}\circ\f G).$$
\end{lem}
\begin{proof}
Let $p_k\cl X\times X\to X$ be the $k$-th projection ($k=1,2$).
Then we have
$$p_1^{-1}\Dr\Homo(\f F,\field_X)\otimes p_2^{-1}\f G\isoto
\Dr\Homo(p_1^{-1}\f F,p_2^{-1}\f G)\quad\text{for any $\f G\in\DI{X}$.}$$
Hence we have
\begin{align*}
\Dr\Homo(\f F,\field_X)\otimes\bl(\wC_{\Delta_X}\circ\f G\br)
&\simeq
\Dr p_1{}_{!!}\bl(p_1^{-1}\Dr\Homo(\f F,\field_X)
\otimes p_2^{-1}\f G\otimes\wC_{\Delta_X}\br)\\
&\simeq
\Dr p_1{}_{!!}\bl(
\Dr\IndHomo(p_1^{-1}\f F,p_2^{-1}\f G)\otimes\wC_{\Delta_X}\br)\\
&\simeq
\Dr p_1{}_{!!}\Dr\IndHomo(p_1^{-1}\f F,p_2^{-1}\f G\otimes\wC_{\Delta_X})\\
&\simeq\Dr\IndHomo\Bl(\f F,
\Dr p_1{}_{!!}\bl(p_2^{-1}\f G\otimes\wC_{\Delta_X}\br)\Br)
\simeq\Dr\IndHomo(\f F,\wC_{\Delta_X}\circ\f G).
\end{align*}
\end{proof}

\begin{cor}
        Assume that $i\cl Z\hookrightarrow X$ is a closed immersion 
and $\f F\in\Db(\I(\field_X))$ satisfies the condition
        $$\SS_0(\f F)\cap T^*_ZX\subset T^*_XX.$$
        Then we have an isomorphism
        $$i^{-1}\f F\otimes\omega_{Z/X}\isoto i^!\f F.$$
\end{cor}

\begin{proof}
We have 
$ i^{-1}\me{F}\otimes\omega_{Z/X}\simeq i^{-1}\me{F}\otimes
  i^{-1}\Dr\me{H}om(\field_Z,\field_X)\simeq i^{-1}\Dr\me{IH}om(\field_Z,\me{F})
  \simeq i^!\me{F}$. 
\end{proof}

\begin{lem}\label{LemMS2}
Let $\Omega\subset\bdot T^*X$ be an open subset
and $\f K\in\DI{Y\times\Omega}$.
Assume that
$$\SS(\f K)^a\cap \big(T^*Y\times\liouville_X(\Omega)\big)=\varnothing,$$ 
where $a$ denotes the antipodal map.
Then
$$(\f K\circ \on{K}_{T^*X})|_{Y\times\Omega}=0.$$
\end{lem}

\begin{proof}
We can easily reduce to the case where
$\me{K}\in \Db(\field_{Y\times\Omega})$.
In this case, let us prove that
$$ (\me{K}\circ\on{K}_{T^*X})\otimes\wC_p\simeq 0\quad
\text{for $p\in Y\times\Omega$.}$$
We may assume that $X$, $Y$ are affine and $p=(y_0,x_0;\xi_0)$. We have
$$ \on{K}_{T^*X}\otimes\wC_{(x_0,\xi_0)}\simeq
   \Dlimind_{\delta>0,\varepsilon>0}{\field_{F_{\delta,\varepsilon}}[2\dim X]},$$
where we have set 
$F_{\delta,\varepsilon}=\{\delta\geqs \langle \xi_0,x'-x\rangle
>\varepsilon(|x'-x|+|\xi'-\xi|)\}$.\\
Hence it is enough to show that there exists a
neighborhood $U$ of $p$ such that
$$ (\me{K}\circ\field_{F_{\delta,\varepsilon}})\vert_U\simeq 0$$
for $0<\delta\ll\varepsilon\ll1$.
Let $p_{ij}$ 
be the $(i,j)$-th projection
from $Y\times \Omega\times\Omega$
to $Y\times\Omega$ or $\Omega\times\Omega$.
Then we have
$$\me{K}\circ \field_{F_{\delta,\varepsilon}}
\simeq\Dr p_{13}{}_{!}(p_{12}^{-1}\me{K}\otimes p_{23}^{-1}
\field_{F_{\delta,\varepsilon}}).$$
For $\SS(F_{\delta,\varepsilon})$ contained in a sufficiently small
neighborhood of 
$(\omega_X(p),-\omega_X(p))$,
$\SS\bigl(p_{12}^{-1}\me{K}\otimes p_{23}^{-1}
\field_{F_{\delta,\varepsilon}}\bigr)$
does not intersect $T^*Y\times\{-\langle \xi_0,\d x\rangle\}\times T^*\Omega$.
Since the map $Y\times\Supp(\field_{F_{\delta,\,\varepsilon}})
\to Y\times\R\times T^*X$ induced by $\langle\xi_0,x\rangle$
is proper,
Proposition 5.4.17 in \cite{KS2} implies that
$\bl(\me{K}\circ \field_{F_{\delta,\varepsilon}}\br)\vert_U\simeq0$.
\end{proof}

\begin{prop}\label{Comphelp}
Let $\f K\in\DI{Y\times X}$ be a kernel and $\f F\in\Db(\I(\field_X))$. 
Assume that
$$\SS(\f K)^a\cap \bl(T^*Y\times\SS_0(\f F)\br)
\subset T^*Y\times T^*_XX.$$
        Then we have an isomorphism
        $$\f K\circ\wC_{\Delta_X}\circ\f F\isoto\f K\circ\f F.$$

\end{prop}

\begin{proof}
It is enough to show that
$\f K\circ\bl(\Ker(\wC_{\Delta_X}\to\field_{\Delta_X})\br)\circ\f F\simeq0$.

Let $p\cl Y\times T^*X\to Y\times X$ be the projection.
We have
$$\SS\left(p^{-1}\f K\right)\subset\Big\{((y;\eta),(x,\xi;\xi',0));
\, ((y;\eta),(x;\xi'))\in\SS(\f K)\Big\}.$$
Hence,
\begin{align*}
\SS\left(\pi_X^{-1}\f K\right)&\cap 
\Bl(T^*Y\times\liouville_X(\SS_0(\f F)\setminus T^*_XX)\Br)\\
\kern-0.2em & \subset\kern-0.2em 
\Big\{((y;\eta),(x,\xi;\xi,0));((y;\eta),(x;\xi))\in\SS(\f K),(x,\xi)\in
\SS_0(\f F)\setminus T^*_XX\Big\}
\end{align*}
is empty by assumption. Therefore Lemma \ref{LemMS2} assures that
$$\Supp(p^{-1}\f K\circ \on{K}_{T^*X})\cap
\bl(Y\times\SS_0(\f F)\br)\subset Y\times T^*_XX.$$
Let $p_1\cl Y\times T^*X\to Y$ and $p_2\cl Y\times T^*X\to T^*X$
be the projections.
Then
\eqn
p^{-1}\f K\circ \bl(\mu_X\f F\otimes\wC_{\bdot T^*X}\br)
&\simeq&
p^{-1}\f K\circ\on{K}_{T^*X}\circ \bl(\mu_X\f F\otimes\wC_{\bdot T^*X}\br)\\
&\simeq&\Dr p_1{}_!\Bl((p^{-1}\f K\circ \K_{T^*X})
\otimes p_2^{-1}(\mu_X\f F\otimes\wC_{\bdot T^*X})\Br)\simeq0.\eneqn
This proves the proposition since 
$p^{-1}\f K\circ (\mu_X\f F\otimes\wC_{\bdot T^*X})
\simeq \f K\circ \Dr\pi_{X!!}(\mu_X\f F\otimes\wC_{\bdot T^*X})$
by Lemma \ref{lem:conv} (iii),
and $\Dr\pi_{X!!}\bl(\mu_X\f F\otimes\wC_{\bdot T^*X}\br)\simeq
\Ker\bl(\wC_{\Delta_X}\ra \field_{\Delta_X}\br)\circ\f F$
by Proposition \ref{prop:mudir} (iii).
\end{proof}


        \subsection{Functorial properties of microlocalization}

To study the functorial behavior of the functor $\mu_X$, it is convenient to
introduce various transfer kernels. 
They will be used exclusively inside the proofs in 
  order to keep notations as simple as possible.
In the sequel, we frequently use Lemma \ref{lem:conv} without mentioning it.

  Let $f\cl X\ra Y$ be a morphism of manifolds. 
Let us recall the commutative
diagram:
$$\xymatrix{
T^*X\ar[d]_{\pi_X}&\db{T^*Y\ctimes{Y}X}\ar[l]_(.55){f_d}
\ar[r]^(.6){f_\pi}\ar[d]
&T^*Y\\
X\ar[r]^f&\ \,Y\,.
}$$
We have $f_d^*\omega_X=f_\pi^*\omega_Y$.
Consider the maps
 \begin{eqnarray*}
    (T^*Y\ctimes{Y}X)\times X & \lra[{f_d\times \id_X}] & T^*X\times X, \\
    (T^*Y\ctimes{Y}X)\times Y & \lra[f_\pi\times\on{\id}_Y] & T^*Y\times Y, \\
     T^*Y\times X & \lra[\id_{T^*Y}\times f] & T^*Y\times Y.
 \end{eqnarray*}
  They define morphisms 
 \begin{eqnarray*}
    \Gamma(T^*X\ctimes{X}X,\Omega^1_{T^*X\times X}) & \lra & 
            \Gamma(T^*Y\ctimes{Y}X,
          \Omega^1_{(T^*Y\ctimes{Y}X)\times X}), \\ 
     \Gamma(T^*Y,\Omega^1_{T^*Y\times Y}) & \lra & \Gamma(T^*Y\ctimes{Y}X,
          \Omega^1_{(T^*Y\ctimes{Y}X)\times Y}), \\ 
     \Gamma(T^*Y,\Omega^1_{T^*Y\times Y}) & \lra & \Gamma(T^*Y\ctimes{Y}X,
          \Omega^1_{T^*Y\times X}).
 \end{eqnarray*}
 We denote by $\sigma_{Y\ot X},\sigma_{X\to Y}$ and $\sigma_{X|Y}$ the images
of the section $\sigma_X$, $\sigma_Y$ and $\sigma_Y$ (defined in \ref{sigma}), 
respectively.
We set
\begin{eqnarray*}
   \on{L}_{Y\ot X}&=&\kn L_{\sigma_{Y\ot X}}((T^*Y\ctimes{Y}X)\ctimes{X}X,
               (T^*Y\ctimes{Y}X)\times X),         \\
    \on{L}_{X\to Y}&=&\kn L_{\sigma_{X\to Y}}
         ((T^*Y\ctimes{Y}X)\ctimes{Y}Y,(T^*Y\ctimes{Y}X)\times Y),       \\
   \on{L}_{X|Y}&=&\kn L_{\sigma_{X|Y}}(T^*Y\ctimes{Y}X,T^*Y\times X).    
\end{eqnarray*}
Note that if $f=\id_X\cl X\to X$, 
then these three kernels coincide and
are isomorphic to $\on{L}_X$.

\begin{lem}\label{lem:LL}
  Let $f\cl X\ra Y$ be a morphism of manifolds. There are natural
  isomorphisms 
\bi
\item $\on{L}_X\simeq \Dr\,\, (\on{id}_{T^*X}\times \pi_X)_{!!}\on{K}_{T^*X}$,
\item $(f_d\times\on{id}_X)^{-1}\on{L}_X\simeq \on{L}_{Y\ot X}$,
\item $\on{L}_{X|Y}\simeq\,(\on{id}_{T^*Y}\times f)^{-1}\on{L}_Y$,
\item $\on{K}_{T^*Y}\comp{T^*Y}\on{L}_{X|Y}\simeq \on{L}_{X|Y}$,
\item $\Dr\, (f_{\pi}\times\on{id}_X)_{!!}\on{L}_{Y\ot X}\lra \on{K}_{T^*Y}
\comp{T^*Y}
       \Dr\, (f_{\pi}\times\on{id}_X)_{!!}\on{L}_{Y\ot X} \isoto \on{L}_{X|Y}$,
\item $\Dr\,(f_{\pi}\times\on{id}_X)_{!!}\on{L}_{Y\ot X}\isoto\on{L}_{X|Y}$
      if $f$ is smooth,
\item $(f_{\pi}\times\on{id}_Y)^{-1}\on{L}_Y\simeq \on{L}_{X\to Y}$.
\item Moreover, there is  a morphism 
$\Dr\,(\on{id}_{T^*Y\ctimes{Y}X}\times f)_{!!}
\on{L}_{Y\ot X}\lra \on{L}_{X\to Y}$ which is an isomorphism if $f$ is smooth.
\ei
\end{lem}
The results easily follow from the first part of the paper.

\begin{thm}[proper direct image]\label{properthm}
        Let $f\cl X\ra Y$ be a morphism of manifolds and 
$\f F\in\Db\left(\I\left(\field_X\right)\right)$. Then
        \bi
                \item{we have a natural morphism and a natural isomorphism
$$\Dr\fp{f}_{!!}\fd{f}^{-1}\mu_X\f F\lra \K_{T^*Y}
\circ\Dr\fp{f}_{!!}\fd{f}^{-1}\mu_X\f F
\isoto\mu_Y\left(\Dr f_{!!}\f F\right),$$}
                \item{if $f$ is smooth we get an isomorphism
$$\Dr\fp{f}_{!!}\fd{f}^{-1}\mu_X\f F
\isoto\mu_Y\left(\Dr f_{!!}\f  F\right).$$}
        \ei
\end{thm}

\begin{proof}
We have $\fd{f}^{-1}\mu_X\f F\simeq \on{L}_{Y\ot X}\circ\f F$ 
by Lemma \ref{lem:LL} (ii),
and a natural morphism by Lemma \ref{lem:LL} (v),
$$\Dr\, (f_{\pi}\times\on{id}_X)_{!!}\on{L}_{Y\ot X}\lra \on{K}_{T^*Y}
\comp{T^*Y}
 \Dr\, (f_{\pi}\times\on{id}_X)_{!!}\on{L}_{Y\ot X} \isoto \on{L}_{X|Y}.$$
However $\left(\Dr\, (f_{\pi}\times\on{id}_X)_{!!}L_{Y\ot X}\right)\circ\f F
\simeq\Dr\fp{f}_{!!}\fd{f}^{-1}\mu_X\f F$ 
and $L_{X|Y}\circ\f F\simeq\mu_Y\left(\Dr f_{!!}\f F\right)$. Hence we get 
natural morphisms
$$\Dr\fp{f}_{!!}\fd{f}^{-1}\mu_X\f F\lra 
  \on{K}_{T^*Y}\circ\Dr\fp{f}_{!!}\fd{f}^{-1}\mu_X\f F\isoto
  \mu_Y\left(\Dr f_{!!}\f F\right),$$
which are isomorphisms if $f$ is smooth by Lemma \ref{lem:LL} (vi).
\end{proof}

\begin{prop}[inverse image]\label{inversethm}
        Let $f\cl X\ra Y$ be a morphism of manifolds and $\f G\in\Db\left(\I\left(\field_Y\right)\right)$. Then
        \bi
                \item{we have a natural morphism 
                $$\fd{f}^{-1}\mu_X(f^{-1}\f G)\lra \fp{f}^{-1}\mu_Y\f G,$$
                which is an isomorphism if $f$ is smooth,}
                \item{we have a natural morphism
                $$\mu_X(f^{-1}\f G)\lra \Dr\fd{f}_*\fp{f}^{-1}\mu_Y\f G.$$}
        \ei
\end{prop}

\begin{proof}
We have 
$$\fd{f}^{-1}\mu_X\left(f^{-1}\f G\right)\simeq \on{L}_{Y\ot X}\circ f^{-1}\f G\qquad
\textrm{and}\qquad\fp{f}^{-1}\mu_Y\f G\simeq \on{L}_{X\to Y}\circ\f G.$$
Since $\on{L}_{Y\ot X}\circ f^{-1}\f G\simeq 
\Bigl(\Dr\,\,(\on{id}_{T^*Y\ctimes{Y}X}\times f)_{!!} \on{L}_{Y\ot X}\Bigr)
\circ \f G$,  we deduce a morphism by Lemma \ref{lem:LL} (viii):
$$\fd{f}^{-1}\mu_X\left(f^{-1}\f G\right)\simeq
\Bigl(\Dr\,\,(\on{id}_{T^*Y\ctimes{Y}X}\times f)_{!!} \on{L}_{Y\ot X}\Bigr)
\circ\f G\lra \on{L}_{X\la Y}\circ\f G\simeq\fp{f}^{-1}\mu_Y\f G,$$
which is an isomorphism whenever $f$ is smooth. By adjunction we get then the inverse image morphism
$\mu_X(f^{-1}\f G)\lra \Dr\fd{f}_*\fp{f}^{-1}\mu_Y\f G$.
\end{proof}

\begin{thm}[embedding case]\label{immersionthm}
        Let $f\cl X\hookrightarrow Y$ be a closed embedding.
Then the following statements hold:
for $\f G\in\Db\left(\I\left(\field_Y\right)\right)$.
\bi
\item
we have a natural morphism 
$$\Dr\fd{f}{}_{!!}\fp{f}^{-1}\mu_Y(\f G)\lra \mu_X(f^{-1}\f G),$$
\item
if $X$ is non characteristic for $\f G$ 
{\rm(}\emph{i.e.} 
$\SS_0(\f G)\cap T^*_XY\subset T^*_YY${\rm)},
then the morphism in {\rm (i)} is an isomorphism
and $\SS_0(f^{-1}\f G)\subset \fd{f}\fp{f}^{-1}\SS_0(\f G)$.
        \ei
\end{thm}

\begin{proof}
(i)\quad Consider the following diagrams
$$
\ba[t]{c}
\xymatrix{
T^*X&\db{T^*Y\ctimes{Y}X}\ar[l]_(.6){f_d}\ar@{ (->}[d]^(.6){f_\pi}\\
&T^*Y}
\ea
\hs{4ex}\text{and}\hs{4ex}
\ba[t]{c}
\xymatrix
{
X&
\db{T^*Y\ctimes{Y}X}\ar@{=}[d]\ar[l]_(.6)p & 
\db{\bigl(T^*Y\ctimes{Y}X\bigr)\times X}\ar[r]^(.7){p_2}\ar[l]_(.6){p_1}
\ar@{^(->}[d]^(.55){f'} & {X}\ar@<-1pt>@{ (->}[d]^{f}\\
&\db{T^*Y\ctimes{Y}X} & 
\db{\bigl(T^*Y\ctimes{Y}X\bigr)\times Y}\ar[r]^(.7){p'_2}
\ar[l]_(.6){p'_1} & {Y.}}
\ea
$$
We have
$$\fd{f}^{-1}\mu_X\left(f^{-1}\f G\right)\simeq \on{L}_{Y\ot X}\circ f^{-1}\f G\qquad\textrm{and}\qquad\fp{f}^{-1}\mu_Y\f G\simeq \on{L}_{X\to Y}\circ\f G.$$
Since $f$ is a closed immersion, $\fd{f}$ is smooth and we get 
$$\fd{f}^!\mu_X\left(f^{-1}\f G\right)\simeq 
\bl(\on{L}_{Y\ot X}\circ f^{-1}\f G\br)\otimes\omega_{T^*Y\ctimes{Y}X/T^*X}.$$
The cotangent bundles being canonically orientable, we have 
$$\omega_{T^*Y\ctimes{Y}X/T^*X}\simeq p^{-1}
\omega_{X/Y}[2(\dim Y-\dim X)]\simeq p^{-1}\omega^{\otimes -1}_{X/Y},$$
where $p\cl T^*Y\ctimes{Y}X\to X$ is the projection.
Hence we get
$$\fd{f}^!\mu_X\left(f^{-1}\f G\right)
\simeq \bl(\on{L}_{Y\ot X}\circ f^{-1}\f G\br)
\otimes p^{-1}\omega^{\otimes -1}_{X/Y}.$$
Now since $f'$ is a closed immersion, 
$\on{L}_{Y\ot X}\simeq f'\,{}^{!}\on{L}_{X\to Y}$
using Proposition \ref{CloRes}, which induces a morphism
$${f'}^{-1}\on{L}_{X\to Y}\ra \on{L}_{Y\ot X}\otimes\omega^{\otimes-1}_{X\times(T^*Y\ctimes{Y}X)/Y\times(T^*Y\ctimes{Y}X)}\simeq
\on{L}_{Y\ot X}\otimes p^{-1}_2\omega^{\otimes-1}_{X/Y}
\simeq\on{L}_{Y\ot X}\otimes p^{-1}_1p^{-1}\omega^{\otimes-1}_{X/Y}.$$
Then the preceding morphism together with the adjunction morphism 
$\id\ra\Dr f'_*f^{\prime-1}\simeq\Dr f'_{!!}f^{\prime-1}$ provides a morphism
\begin{eqnarray*}
\fp{f}^{-1}\mu_Y\f G &\simeq& \on{L}_{X\to Y}\circ\f G
=\Dr p'_{1!!}(\on{L}_{X\to Y}\otimes {p'}^{-1}_2\f G)
\simeq\Dr p'_{1!!}\Dr f'_{!!}{f'}^{-1}(\on{L}_{X\to Y}
\otimes {p'}^{-1}_2\f G)\\
&\lra& \Dr p_{1!!}(\on{L}_{Y\ot X}\otimes p^{-1}_1p^{-1}\omega^{\otimes-1}_{X/Y}
\otimes p^{-1}_2f^{-1}\f G)
\simeq (\on{L}_{Y\ot X}\circ f^{-1}\f G)\otimes p^{-1}
\omega^{\otimes -1}_{X/Y}.
\end{eqnarray*}
Finally we obtain a morphism
$$\fp{f}^{-1}\mu_Y\f G\lra (\on{L}_{Y\ot X}\circ f^{-1}\f G)\otimes p^{-1}
\omega^{\otimes -1}_{X/Y}\simeq\fd{f}^{-1}\mu_X\left(f^{-1}\f G\right)
\otimes  p^{-1}\omega^{\otimes -1}_{X/Y}
\simeq\fd{f}^!\mu_X\left(f^{-1}\f G\right),$$
and by adjunction the desired morphism
$$\Dr\fd{f}_{!!}\fp{f}^{-1}\mu_Y(\f G)\lra \mu_X(f^{-1}\f G).$$

\medskip
\noindent
(ii)\quad
Assume now that $X$ is non characteristic for $\f G$. 
By induction we may assume that $X$ is a hypersurface in $Y$.
For $p\in T^*X$, let us show that
$\Dr\fd{f}_{!!}\fp{f}^{-1}\mu_Y(\f G)\otimes\wC_p
\isoto\mu_X(f^{-1}\f G)\otimes\wC_p$. 

Assume first that $p\in T^*_XX$.
Since $X$ is non characteristic for $\me{G}$ we get
\begin{align*}
  \Dr f_{d!!}f_{\pi}^{-1}\mu_Y\me{G}&\otimes \wC_p  \simeq
    \Dr f_{d!!}\big(f_{\pi}^{-1}\mu_Y\me{G}\otimes \wC_{T^*_XY}\big)
    \otimes \wC_p \simeq
    \Dr f_{d!!}\big(f_{\pi}^{-1}\big(\mu_Y\me{G}\otimes\wC_{T^*_YY}\big)
     \big)\otimes \wC_p \\
   & \simeq \Dr f_{d!!}\big( f_{\pi}^{-1}(\pi^{-1}_Y\me{G}\otimes 
       \wC_{T^*_YY})\big)\otimes\wC_p 
     \simeq \Dr f_{d!!}\big( f_d^{-1}\pi_X^{-1}f^{-1}\me{G}\otimes
      \wC_{T^*_YY\ctimes{Y}X}\big) \otimes \wC_p\\
  &\simeq \pi_X^{-1}f^{-1}\me{G}\otimes \Dr f_{d!!}\wC_{T^*_YY\ctimes{Y}X}
   \otimes\wC_p\simeq
    \pi_X^{-1}f^{-1}\me{G}\otimes \wC_p\\
  &\simeq \mu_Xf^{-1}\me{G}\otimes\wC_p. 
\end{align*}

Assume now that $p\not\in T^*_XX$. 
Consider the following diagram
$$\xymatrix{
T^*X\times X\ar@{ (->}[r]^(.55){f_1}&
T^*X\times Y\ar@/^20pt/[rr]^{q_2}\ar[d]_{q_1}\ar@{}[dr]|{\square}
&\db{(T^*Y\ctimes{Y}X)\times Y}\ar[r]_(.7){p'_2}\ar[d]^(.55){p'_1}
\ar[l]_(.55){r}
&Y\\
&T^*X&\db{T^*Y\ctimes{Y}X}\ar[l]^{f_d}
} 
$$
Note that
$$ \Dr f_{d!!} f_{\pi}^{-1}\mu_Y\me{G}\simeq
    (\Dr r_{!!}\on{L}_{X\to Y})\circ \me{G} \qquad
    \text{and} \qquad \mu_Xf^{-1}\me{G}\simeq
       \on{L}_X\circ f^{-1}\me{G}\simeq
    (\Dr f_1{}_{!!}\on{L}_X)\circ\me{G}. $$
Hence we have to prove that 
$$ \big(\Dr r_{!!}\on{L}_{X\to Y}\otimes \wC_p\big)\circ\me{G} \simeq
   \big(\Dr f_1{}_{!!}\on{L}_X\otimes \wC_p\big)\circ\me{G}. $$
Here we identify $p\in T^*X$ with $\bl(p,f(\pi_X(p))\br)\in T^*X\times Y$.
Take a local coordinate system $(t,x)=(t,x_1,\ldots,x_n)$
of $Y$ such that $X$ is given by $t=0$ and denote by 
 $(t,x,\tau,\xi)$  and
$(x,\xi)$ the associated coordinates on  $T^*Y$ and  $T^*X$, respectively. 
Set $p=(0,\xi_0)$.
Let $((x,\tau,\xi),(t',x'))$ be the coordinates of $(T^*Y\ctimes{Y}X)\times Y$.
Then $r((x,\tau,\xi),(t',x'))=((x,\xi),(t',x'))$.
We have
$$ \Dr r_{!!}\on{L}_{X\to Y}\otimes\wC_p \simeq \Dr r_{!!}\Bigl(
\Dlimind_{\varepsilon>0,\,R>0}{\field_{\{\tau t'+\langle\xi_0,x'-x\rangle>
    \varepsilon(|t'|+|x'-x|),\,|\tau|<R\}}}[\dim Y]\Br)
\otimes \wC_p.$$  
Since the fiber of
$\{\tau t'+\langle\xi_0,x'-x\rangle>
    \varepsilon(|t'|+|x'-x|),\,|\tau|<R\}$
over $((x,\xi),t',x')$ is a non-empty open interval 
if $R|t'|+\langle\xi_0,x'-x\rangle>\varepsilon(|t'|+|x'-x|)$,
and empty otherwise,
we obtain
$$ \Dr r_{!!}\on{L}_{X\to Y}\otimes\wC_p \simeq
\Bigl(
\Dlimind_{\varepsilon>0,\,R>0}{\field_{\{R|t'|+\langle\xi_0,x'-x\rangle>
    \varepsilon(|t'|+|x'-x|)\}}}[\dim Y-1]\Br)
\otimes \wC_p.$$
Therefore
\begin{align*}  
  \big(\Dr r_{!!}\on{L}_{X\to Y}\otimes\wC_p\big)\circ\me{G}
 & \simeq\bl(\Dlimind_{\varepsilon>0,\,R>0}
{\field_{\{R|t'|+\langle\xi_0,x'-x\rangle>\varepsilon|x'-x|\}}}[\dim X]\otimes 
        \wC_p\br)\circ\f G.
\end{align*}
On the other hand we have
\begin{align*}
 \big(\Dr f_1{}_{!!}\on{L}_X\otimes \wC_p\big)\circ\me{G} & \simeq
\Bl( \Dr f_1{}_{!!}\bl(\Dlimind_{\varepsilon>0}
{\field_{\{\langle\xi_0,x'-x>\varepsilon|x'-x|\}}
[\dim X]}\br)\otimes \wC_p\Br)\circ\f G \\
&\simeq \Bl(\Dlimind_{\varepsilon>0}
{\field_{\{\langle \xi_0,x'-x\rangle>\varepsilon|x'-x|,\,t'=0\}}[\dim X]}
\otimes \wC_p\Br)\circ\me{G}.
\end{align*}
Hence it is enough to show that
$$ \big(\Dlimind_{\varepsilon>0,\,R>0}{\field_{\{R|t'|+\langle\xi_0,x'-x\rangle>
 \varepsilon|x'-x|,\,0<t'\}}}\otimes\wC_p\big)\circ\me{G}
    \simeq 0. $$
Let us set $U_{\varepsilon,\,\delta,\,R}=
\{R\,t'+\langle\xi_0,x'-x\rangle>|x-x|,\,0<t'\le\delta\}$.
For 
$\varepsilon$, $\delta$ sufficiently small and $R$ sufficiently large,
$\SS(\field_{U_{\varepsilon,\,\delta,\,R}})$ is contained in a 
sufficiently small neighborhood of 
$-R\d t'+\langle \xi_0,\d(x-x')\rangle$ on a neighborhood of $p$.
Hence we obtain
$$ \SS(\field_{U_{\varepsilon,\,\delta,\,R}})^a
\cap T^*(T^*X)\times \SS_0(\me{G})\subset
      T^*(T^*X)\times T^*_YX\quad\text{on a neighborhood of $p$
for $R\gg0$.} $$
Therefore Proposition \ref{Comphelp} implies
$$\bl(\field_{U_{\varepsilon,\,\delta,\,R}}
\circ \wC_{\Delta_Y}\br)\circ\me{G} \simeq
   \field_{U_{\varepsilon,\,\delta,\,R}}\circ\me{G}
\quad\text{on a neighborhood of $p$ for $R\gg0$.}$$
Hence we have reduced the problem to
$$ \bl(\Dlimind_{\varepsilon>0,\,\delta>0\,R>0}
\field_{U_{\varepsilon,\,\delta,\,R}}\otimes\wC_p\br)\circ 
   \wC_{\Delta_Y}\simeq 0. $$
Consider the projection on the first and third factors
$$ h\cl T^*X\times Y\times Y \lra
       T^*X\times Y \qquad\text{{\em i.e.}
    $((x;\xi),(t',x'),(t'',x'')) \mapsto ((x;\xi),(t'',x''))$.} $$
Then
$$ \bl(\Dlimind_{\varepsilon>0,\,\delta>0\,R>0}
\field_{U_{\varepsilon,\,\delta,\,R}}\otimes\wC_p\br)
    \circ \wC_{\Delta_Y} \simeq
    \Dr h_{!!}\bl(\Dlimind_{\varepsilon>0,\,\delta>0\,R>0}
\field_{V_{\varepsilon,\,\delta,\,R}}\br)\otimes\wC_p,$$
where $V_{\varepsilon,\,\delta,\,R}
=\{R\,t'+\langle \xi_0,x'-x\rangle>\varepsilon|x'-x|,\,0<t'\le\delta,
\,|x'-x''|\le\delta,\,|t'-t''|\le\delta\}$.
This vanishes by the following lemma.
\end{proof}
\begin{sublem}
Let $(t,t',x,y)=(t,t',x_1,\ldots,x_n,y_1,\ldots,y_n)$ be the coordinates
of $\R\times\R\times \R^n\times \R^n$,
and let $h\cl \R\times \R\times \R^n\times \R^n\to \R\times\R^n$
be the projection, $h(t,t',x,y)=(t',y)$. For $\xi\in\R^n\setminus\{0\}$
and $\delta>0$, set
$V_{\delta}=\{(t,t',x,y) ;\,t+\langle\xi_0,x\rangle>|x|,\,
|x-y|\le\delta,0<t\le\delta,\,|t-t'|\le\delta\}$.
Then 
$$\Supp\bl(\Dr h_!\,\field_{V_{\delta}}\br)\not\ni0.$$
\end{sublem}
\begin{proof}
Let us decompose $h$ into
$\R\times \R\times\R^n\times \R^n\lra[h_1]\R\times \R^n\times \R^n
\lra[h_2]\R\times\R^n$,
where $h_1(t,t',x,y)=(t',x,y)$ and $h_2(t',x,y)=(t',y)$.
When $|x-y|\le\delta$, the fiber
$V_{\delta}\cap h_1^{-1}(t',x,y)$
is
$\{t;\max(0,\,|x|-\langle\xi_0,x\rangle)<t\le\min(\delta,t'+\delta),\,
t'-\delta\le t\}$.
Hence, setting
$$
W_{\delta}=\{(t',x,y);
\max(0,\,|x|-\langle\xi_0,x\rangle)<t'-\delta\le\min(\delta,t'+\delta),\,
|x-y|\le\delta\},$$
we have
$\Dr h_1{}_!\,\field_{V_\delta}\simeq
\field_{W_\delta}$.
Since $\Supp(\field_{W_\delta})\subset\{(t',x,y);\delta\le t'\}$,
we obtain
$$\Supp(\Dr h_!\,\field_{V_\delta})\subset \{(t',y);\delta\le t'\}.$$
\end{proof}


        \subsection{Microlocal convolution of kernels}

Let $X$, $Y$ and $Z$ be manifolds,
and let $p_{ij}$ be the $(i,j)$-th projection from
$T^*X\times T^*Y\times T^*Z$.
As usual, denote by $a\cl T^*X\ra T^*X$ the antipodal map. Then define the antipodal
projection $p_{12}^a$ by
$$ p_{12}^a\cl T^*X\times T^*Y\times T^*Z\lra[{p_{12}}]T^*X\times T^*Y 
\lra[{\on{id}\times a}] T^*X\times T^*Y. $$
For $\me{F}\in\Db(\I(\field_{T^*X\times T^*Y}))$ and $\me{G}\in\Db(\I(\field_{T^*Y\times T^*Z}))$,
we set
$$ \me{F}\overset{a}{\circ} \me{G}=\Dr p_{13!!}\big(p_{12}^{a-1}\me{F}
    \otimes p_{23}^{-1}\me{G}\big). $$
In an analogous way, for $S_1\subset T^*X\times T^*Y$ and 
$S_2\subset T^*Y\times T^*Z$, we set
$$ S_1\overset{a}{\ctimes{T^*Y}}S_2=
p_{12}^{a-1}(S_1)\cap p_{23}^{-1}(S_2)\subset T^*X\times T^*Y\times T^*Z. $$
Now we are ready to state the main theorem:
\begin{thm}[Microlocal convolution of kernels]\label{mainthm}
  Let $\me{K}_1\in\Db(\I(\field_{X\times Y}))$ 
and $\me{K}_2\in\Db(\I(\field_{Y\times Z}))$.
  \bi
\item There is a natural morphism
       \eq
&& \mu_{X\times Y}\me{K}_1\overset{a}{\circ}
           \mu_{Y\times Z}\me{K}_2 \lra
          \mu_{X\times Z}(\me{K}_1\circ \me{K}_2). \label{eq:compmor}
\eneq

\item Assume the non-characteristic condition
\eq
&& \SS_0(\me{K}_1)\overset{a}{\ctimes{T^*Y}}\SS_0(\me{K}_2) \cap
           (T^*_XX\times T^*Y\times T^*_ZZ) \subset
              T^*_XX\times T^*_YY\times T^*_ZZ, \label{eq:nonch}
\eneq

 Then \eqref{eq:compmor}
        is an isomorphism outside $$\ol{p_{13}\big(\SS_0(\me{K}_1)
         \overset{a}{\ctimes{T^*Y}}\SS_0(\me{K}_2)\cap T^*X\times T^*_YY
         \times T^*Z\big)}.$$
  \ei
\end{thm}
\begin{proof}
(a)\quad  We shall first construct the morphism. 
Consider the manifolds
$\me{X}_1=X\times Y$, $\me{X}_2=Y\times Z$ and  
$\me{X}=\me{X}_1\times\me{X}_2=X\times Y\times Y\times Z$ 
together with the diagonal embedding
$$ \me{Y}\seteq X\times Y\times Z \overset{j}{\hookrightarrow} \me{X}. $$
Denote by $\me{Z}=X\times Z$, and 
let $q_{13}\cl\me{Y}\ra\me{Z}$ be the projection. The map
$$ T^*\me{Y} \hookrightarrow \me{Y}\ctimes{\me{X}}T^*\me{X} 
\qquad\text{given by}\qquad
   (x,y,z;\xi,\eta,\zeta)\mapsto (x,y,y,z;\xi,-\eta,\eta,\zeta) $$
defines the cartesian square in the following commutative diagram:
$$ \xymatrix{
    {T^*\me{Y}} \ar[d] \ar@<-1pt>@{ (->}[r]\ar@/^25pt/[rr]^{p} 
\ar@/_30pt/[dd]_{q} 
     \ar@{}[rd]|{\square} & \ddb{\me{Y}\ctimes{\me{X}}T^*\me{X}}
     \ar[d]_(.6){j_d}  \ar@<-1pt>@{ (->}[r]^(.6){j_\pi} &  {T^*\me{X}} \\
    \db{\me{Y}\ctimes{\me{Z}}T^*\me{Z}} \ar@<-1pt>@{ (->}[r]_(.6){q_{13}{}_d}
 \ar[d]^(.6){q_{13\pi}} & 
    {T^*\me{Y}} & \\
  {T^*\me{Z}} & &
     } $$
By Proposition \ref{boxcomposition}, we have an isomorphism
$$
\on{K}_{T^*\me{X}}\circ
(\mu_{\me{X}_1}\me{K}_1\boxtimes\mu_{\me{X}_2}\me{K}_2)
\simeq\mu_{\me{X}}(\me{K}_1\boxtimes\me{K}_2). 
$$
By Theorem \ref{immersionthm} we have a morphism
\eq\label{eq:iso1}
&& \Dr j_{d!!}j_{\pi}^{-1}\mu_{\me{X}}(\me{K}_1\boxtimes \me{K}_2) \lra
   \mu_{\me{Y}}(j^{-1}(\me{K}_1\boxtimes\me{K}_2)).
\eneq
Since $q_{13}$ is smooth we also have an isomorphism
by Theorem \ref{properthm} (ii)
$$ \Dr q_{13\pi!!}q_{13d}^{-1}\mu_\me{Y}(j^{-1}(\me{K}_1\boxtimes\me{K}_2))
   \isoto\mu_{\me{Z}}(\Dr q_{13!!}j^{-1}(\me{K}_1\boxtimes \me{K}_2))
  \simeq \mu_{\me{Z}}(\me{K}_1\circ\me{K}_2). $$
Hence we get a morphism
\eq \label{eq:iso2}
&&\Dr q_{!!}p^{-1}
\Bl(\on{K}_{T^*\me{X}}\circ
(\mu_{\me{X}_1}\me{K}_1\boxtimes\mu_{\me{X}_2}\me{K}_2)\Br) \lra
    \mu_{\me{Z}}(\me{K}_1\circ\me{K}_2). 
\eneq
Hence we obtain
\eqn
&&\mu_{\me{X}_1}\me{K}_1\overset{a}{\circ}\mu_{\me{X}_2}\me{K}_2
\simeq\Dr q_{!!}p^{-1}(\mu_{\me{X}_1}\me{K}_1\boxtimes\mu_{\me{X}_2}\me{K}_2)\\
&&\hs{10ex}\lra\Dr q_{!!}p^{-1}
\Bl(\on{K}_{T^*\me{X}}\circ
(\mu_{\me{X}_1}\me{K}_1\boxtimes\mu_{\me{X}_2}\me{K}_2)\Br) \lra
    \mu_{\me{Z}}(\me{K}_1\circ\me{K}_2). 
\eneqn

\medskip
\noindent
(b)\quad
By Theorem \ref{immersionthm},
\eqref{eq:iso1} is an isomorphism under the non-characteristic hypothesis,
and hence
\eqref{eq:iso2} is also an isomorphism under the same hypothesis.

Therefore in order to show (ii), it is enough to show that
\begin{align}
\ba{l}
\mu_{\me{X}_1}\me{K}_1\overset{a}{\circ}
   \mu_{\me{X}_2}\me{K}_2\simeq 
     \Dr q_{!!}p^{-1}\Bl(\on{K}_{T^*\me{X}}\circ
(\mu_{\me{X}_1}\me{K}_1\boxtimes\mu_{\me{X}_2}\me{K}_2)\Br)\\
\qquad\text{outside ${p_{13}\big(\SS_0(\me{K}_1)
         \overset{a}{\ctimes{T^*Y}}\SS_0(\me{K}_2)\cap T^*X\times T^*_YY
         \times T^*Z\big)}$.}\ea
\label{eq:conv}
\end{align}
First note that 
\begin{align*}
\mu_{\X_1}\me{K}_1\overset{a}{\circ} \mu_{\X_2}\me{K}_2&\simeq
(\on{K}_{T^*\X_1}\circ\mu_{\X_1}\me{K}_1)\overset{a}{\circ} 
(\on{K}_{T^*\X_1}\circ\mu_{\X_2}\me{K}_2)\\[1ex]
 &\simeq \big(\on{K}_{T^*\me{X}_1}\comp{T^*Y}\on{K}_{T^*\me{X}_2}\big)
    \circ\big(\mu_{\me{X}_1}\me{K}_1\boxtimes \mu_{\me{X}_2}\me{K}_2\big).
\end{align*}
Consider the diagram
$$ \xymatrix{ 
& {T^*\me{Y}\times T^*\me{X}}  
\ar[dr]^{p'=(p,\on{id})} \ar[dl]_{q'=(q,\on{id})} & \\
   {T^*\me{Z}\times T^*\me{X}} & & {T^*\me{X}\times T^*\me{X}} }$$
Then we have
$$\Dr q_{!!}p^{-1}\Bl(\on{K}_{T^*\me{X}}\circ
(\mu_{\me{X}_1}\me{K}_1\boxtimes\mu_{\me{X}_2}\me{K}_2)\Br)
\simeq
\bl(\Dr q'_{!!}p^{\prime-1}\on{K}_{T^*\me{X}}\br)\circ
(\mu_{\me{X}_1}\me{K}_1\boxtimes\mu_{\me{X}_2}\me{K}_2).$$
Using Proposition \ref{InvIm} and Corollary \ref{DirIm}, we have
$$ \Dr q'_{!!}p^{\prime-1}\on{K}_{T^*\me{X}}\simeq
\kn L_{\sigma}\big(T^*\me{Y},T^*\me{Z}\times T^*\me{X}\big), $$
where $T^*\me{Y}$ is embedded into $T^*\me{Z}\times T^*\me{X}$
by $(q,p)$ and
the section $\sigma$ is given by 
$$\sigma=(\omega_X,\omega_Z,-\omega_X,-\omega_Y,
-\omega_Y,-\omega_Z).$$

In order to see \eqref{eq:conv} under the non-characteristic hypothesis,
it is enough to show that
\eq
&&\parbox{75ex}{
$\on{K}_{T^*\me{X}_1}\comp{T^*Y}\on{K}_{T^*\me{X}_2}
\lra
\kn L_{\sigma}\big(T^*\me{Z}\ctimes{T^*\me{Z}}T^*\me{Y},
T^*\me{Z}\times T^*\me{X}\big)$
is an isomorphism on 
$T^*\me{Z}\times \bl(T^*X\times\bdot{T}^*(Y\times Y)\times T^*Z\br)
\subset T^*\me{Z}\times T^*\me{X}$.}
\eneq
However it is a consequence of Proposition \ref{prop:kernel}
(note that (iii) and (v) in the proposition
fail on $T^*X\times T^*_{Y}Y\times T^*Z$).
\end{proof}

      \subsection{A vanishing theorem for microlocal holomorphic functions}

\begin{thm}\label{th:Odegree0}
  Let $X$ be a complex manifold of dimension $n$. 
Then, $\mu_X(\me{O}_X)\vert_{\bdot{T}^*X}$ is concentrated in degree $-n$.
\end{thm}
\begin{proof}
We may assume $X=\C^n$.
Let $q_1\cl T^*X\times X\to T^*X$ and
$q_2\cl T^*X\times X\to X$ be the projections.
Let $p=(x_0,\xi_0)\in\bdot T^*X$. Then, we have
$$ \mu_X(\me{O}_X)\otimes\wCC_p\simeq
\wCC_p\otimes\Dr q_{1!!} 
\bl(\Dlimind_{\varepsilon,\delta>0} (\C_{\on{F}_{\delta,\varepsilon}}\otimes
      q_2^{-1}\me{O}_X)\br)[2n], $$
where $F_{\delta,\varepsilon}=
\big\{\bl((x,\xi),x'\br);\delta\geqs\langle\xi_0,x'-x\rangle>
\varepsilon|x'-x|\big\}$. Hence 
it is enough to show that
$$ \Dr q_{1!}(\C_{\on{F}_{\delta,\varepsilon}}\otimes
      q_2^{-1}\me{O}_X) $$
is concentrated in degree $n$. We have
\begin{align*}
  \Dr q_{1!}(\C_{\on{F}_{\delta,\varepsilon}}\otimes
      q_2^{-1}\me{O}_X)_{(x_1,\xi_1)}&\simeq
   \Dr\Gamma_c\left(\{x'\in X;\delta\geqs \langle\xi_0,x'-x_1\rangle >
\varepsilon |x'-x_1|\}, \me{O}_X\right).
\end{align*}
The cohomology with compact support of $\me{O}_X$ on the difference of
two convex open subsets is concentrated in degree $n$.
\end{proof}
Now,
$H^{-n}\bl(\mu_X(\me{O}_X)\vert_{\bdot{T}^*X}\br)$ 
has a structure of $\me{E}_X\vert_{\bdot{T}^*X}$-module,
i.e. there exists a canonical ring homomorphism
$\me{E}_X\vert_{\bdot{T}^*X}\to
\me{E}nd\bl(H^{-n}(\mu_X(\me{O}_X)\vert_{\bdot{T}^*X})\br)$.

Indeed, let $p_k\cl X\times X\to X$ be the $k$-th projection,
and $\me{O}_{X\times X}^{(0,n)}\seteq\me{O}_{X\times X}
\otimes_{p_2^{-1}\me{O}_X}p_2^{-1}\me{O}_X^{(n)}$.
We have morphisms
$\Dr p_1{}_!(\me{O}_{X\times X}^{(0,n)}[n]\otimes p_2^{-1}\me{O}_X)\to
\Dr p_1{}_!(\me{O}_{X\times X}^{(0,n)}[n])\to\me{O}_X$
which induce
$\me{O}_{X\times X}^{(0,n)}[n]\to \Dr\Homo(p_2^{-1}\me{O}_X, p_1^!\me{O}_X)$.
Thus we obtain
\begin{align*}
\me{E}_X\to\mu_{\Delta_X}(\me{O}_{X\times X}^{(0,n)}[n])
&\to\mu_{\Delta_X}\bl( \Dr\Homo(p_2^{-1}\me{O}_X, p_1^!\me{O}_X)\br)\\
&\qquad\simeq\mu hom(\me{O}_X,\me{O}_X)\simeq
\Dr\me{H}om\bl(\mu_X(\me{O}_X),\mu_X(\me{O}_X)\br).
\end{align*}
Hence, Theorem \ref{th:Odegree0} implies that
$\mu_X(\me{O}_X)\vert_{\bdot{T}^*X}$ belongs to
$\Db\bl(\Mod\bl(\me{E}_X\vert_{\bdot{T}^*X},\on{I}(\C_{\bdot{T}^*X})\bl)\br)$,
the derived category 
of the abelian category $\Mod\bl(\me{E}_X\vert_{\bdot{T}^*X},
\on{I}(\C_{\bdot{T}^*X})\bl)$ 
of ind-sheaves $\f F$ on $\bdot{T}^*X$ endowed with
a ring homomorphism
$\me{E}_X\vert_{\bdot{T}^*X}\to\me{E}nd\,(\f F)$.
This implies the following theorem.
\begin{thm}
  Let $X$ be a complex manifold. Then 
$\me{F}\mapsto\mu hom(\me{F},\me{O}_X)\vert_{\bdot{T}^*X}$ 
is a well defined functor
from $\Db(\C_X)$ to $\Db(\me{E}_X\vert_{\bdot{T}^*X})$.
\end{thm}


\end{document}